\documentclass[aos]{imsart}

\RequirePackage{amsthm,amsmath,amsfonts,amssymb}
\RequirePackage[authoryear]{natbib}%% uncomment this for author-year citations
\RequirePackage[colorlinks,citecolor=blue,urlcolor=blue]{hyperref}
\RequirePackage{graphicx}

\startlocaldefs
\theoremstyle{plain}

\newtheorem{theorem}{Theorem}[section]
\newtheorem{lemma}[theorem]{Lemma}
\newtheorem{proposition}[theorem]{Proposition}
\newtheorem{corollary}[theorem]{Corollary}
\theoremstyle{remark}
\newtheorem{definition}[theorem]{Definition}
\newtheorem{remark}[theorem]{Remark}

\newtheorem{conjecture}[theorem]{Conjecture}
\newtheorem{example}[theorem]{Example}

\usepackage{amsmath,amsfonts,amssymb,amsthm,graphicx}
\usepackage{dsfont}    % it's a bitmap font (on Volodya's computer)
\usepackage{enumerate}
\usepackage{color}
\usepackage{multirow}\usepackage{graphicx}
\usepackage{amsfonts}
\usepackage{amsmath}
\usepackage{amssymb}

\usepackage{color}
\usepackage{dsfont} 

\usepackage{pgfplots}
\usepackage{subcaption}
\usepgfplotslibrary{fillbetween}
\usetikzlibrary{patterns}
\usepackage{amssymb}
\usepackage{graphicx}
\usepackage{amsmath}
\usepackage{bm}
\usepackage{tikz}

\usepackage{caption}
  
\usetikzlibrary{decorations.pathreplacing}

\usepackage{tkz-graph}
\usepackage{pgfplots}
\usetikzlibrary{calc,patterns}

\usepackage{subcaption}
\usepackage{url}
\usepackage{fancyhdr}
\usepackage{indentfirst}

\usepackage{enumerate}
\usepackage{dsfont}
\usepackage[british]{babel}
\usepackage{amsthm}
\usepackage{color}
\usepackage{natbib}
\usepackage{comment}
\usepackage{tikz-cd}
\usepackage{graphicx}
\usepackage{amsfonts}
\usepackage{amsmath}
\usepackage{amssymb}
\usepackage{url}
\usepackage{bbm}
\usepackage{fancyhdr}
\usepackage{indentfirst}
\usepackage{enumerate}
\usepackage{dsfont}
\usepackage[british]{babel}
\usepackage{cleveref}
\usepackage{amsthm}
\usepackage{color}
\renewcommand{\cite}{\citet}
\usepackage{comment}

\usepackage{algorithm}
\usepackage[noend]{algpseudocode}  % algpseudocode loads algorithmicx
\renewcommand{\d}{\,\mathrm{d}}

\newcommand{\p}{\mathbb{P}}

\usepackage{xcolor}

\newcommand{\E}{\mathbb{E}}    % expectation
\newcommand{\R}{\mathbb{R}}    % real numbers
\newcommand{\N}{\mathbb{N}}    % natural numbers

\newcommand{\I}{\mathcal{I}}   % the Borel sigma-field on [0,1]
 
\usepackage{tikz}
\usetikzlibrary{arrows.meta,positioning}

\newcommand{\bary}{\mathrm{bary}}

\newcommand{\conv}{\mathrm{Conv}}

\newcommand{\ee}{\varepsilon}

\newcommand{\n}[1]{\left\lVert#1\right\rVert}
\newcommand{\B}{\mathcal{B}}

\renewcommand{\P}{\mathcal{P}}

\newcommand{\Ber}{\mathrm{Ber}}

\newcommand{\ba}{\boldsymbol\alpha}
\newcommand{\bb}{\boldsymbol\beta}
\newcommand{\bx}{\mathbf{x}}
\newcommand{\by}{\mathbf{y}}
\newcommand{\bl}{\boldsymbol{\lambda}}

\newcommand{\bmu}{\boldsymbol{\mu}}

\newcommand{\bnu}{\boldsymbol{\nu}}
\newcommand{\T}{\mathcal{T}}
\newcommand{\Tau}{\mathbb{T}}

\def\d{\mathrm{d}}
\DeclareMathOperator\supp{supp}
\DeclareMathOperator\aff{aff}
\DeclareMathOperator\ri{ri}
\def\laweq{\buildrel \mathrm{law} \over =}

\def\lawis{\buildrel \mathrm{law} \over \sim}

\newcommand{\bone}{ {\mathbbm{1}} }
\renewcommand{\S}{\mathscr{S}}
\renewcommand{\P}{\mathcal{P}}

\newcommand{\Q}{\mathcal{Q}}
\newcommand{\F}{\mathcal{F}}
\newcommand{\G}{\mathcal{G}}
\renewcommand{\H}{\mathbb{H}}

\usepackage{mathrsfs}
\newcommand{\M}{\mathcal{M}}
\newcommand{\cM}{\mathcal{M}}
\newcommand{\e}{\mathbf{1}}
\newcommand{\z}{\mathbf{0}}
\newcommand{\Leb}{\mathrm{U_1}}

\newcommand{\K}{\mathcal{K}}

\newcommand{\lcx}{\preceq_{\mathrm{cx}}}
\newcommand{\lst}{\preceq_{\mathrm{st}}}
\newcommand{\gst}{\succeq_{\mathrm{st}}}
\newcommand{\pst}{\prec_{\mathrm{st}}}
\newcommand{\sst}{\succ_{\mathrm{st}}}

\newcommand{\spn}{\mathrm{Span}}
\newcommand{\gcx}{\succeq_{\mathrm{cx}}}
\newcommand{\X}{\mathfrak{X}}
\newcommand{\Y}{\mathfrak{Y}}

\renewcommand{\ge}{\geqslant}
\renewcommand{\le}{\leqslant}
\renewcommand{\geq}{\geqslant}
\renewcommand{\leq}{\leqslant}
\renewcommand{\epsilon}{\varepsilon}

 \pgfdeclarepatternformonly{south east lines}{\pgfqpoint{-0pt}{-0pt}}{\pgfqpoint{3pt}{3pt}}{\pgfqpoint{3pt}{3pt}}{
    \pgfsetlinewidth{0.4pt}
    \pgfpathmoveto{\pgfqpoint{0pt}{3pt}}
    \pgfpathlineto{\pgfqpoint{3pt}{0pt}}
    \pgfpathmoveto{\pgfqpoint{.2pt}{-.2pt}}
    \pgfpathlineto{\pgfqpoint{-.2pt}{.2pt}}
    \pgfpathmoveto{\pgfqpoint{3.2pt}{2.8pt}}
    \pgfpathlineto{\pgfqpoint{2.8pt}{3.2pt}}
    \pgfusepath{stroke}}
\endlocaldefs

\begin{document}

\begin{frontmatter}
\title{On the existence of powerful p-values \\ and e-values for composite hypotheses}
%\title{A sample article title with some additional note\thanksref{t1}}
\runtitle{Existence of powerful p/e-values}
%\thankstext{T1}{A sample additional note to the title.}

\begin{aug}
%%%%%%%%%%%%%%%%%%%%%%%%%%%%%%%%%%%%%%%%%%%%%%%
%% Only one address is permitted per author. %%
%% Only division, organization and e-mail is %%
%% included in the address.                  %%
%% Additional information can be included in %%
\author[A]{\fnms{Zhenyuan}~\snm{Zhang}\ead[label=e1]{zzy@stanford.edu}},
\author[B]{\fnms{Aaditya}~\snm{Ramdas}\ead[label=e2]{aramdas@cmu.edu}}
\and
\author[C]{\fnms{Ruodu}~\snm{Wang}\ead[label=e3]{wang@uwaterloo.ca}}
%%%%%%%%%%%%%%%%%%%%%%%%%%%%%%%%%%%%%%%%%%%%%%
%% Addresses                                %%
%%%%%%%%%%%%%%%%%%%%%%%%%%%%%%%%%%%%%%%%%%%%%%
\address[A]{Department of Mathematics, Stanford University\printead[presep={,\ }]{e1}}

\address[B]{Depts.\ of Statistics \& Data Science, and Machine Learning, Carnegie Mellon Univ.\printead[presep={,\ }]{e2}}
\address[C]{Department of Statistics and Actuarial Science, University of Waterloo\printead[presep={,\ }]{e3}}
\end{aug}

\begin{abstract}
Given a composite null $\P$ and composite alternative $\Q$, when and how can we construct a p-value whose distribution is exactly uniform under the null, and stochastically smaller than uniform under the alternative? Similarly, when and how can we construct an e-value whose expectation exactly equals one under the null, but its expected logarithm under the alternative is positive? 
We answer these basic questions, and other related ones, when $\P$ and $\Q$ are convex polytopes (in the space of probability measures). We prove that such constructions are possible if and only if $\Q$ does not intersect the span of $\P$. If the p-value is allowed to be stochastically larger than uniform under $P\in\P$, and the e-value can have expectation at most one under $P\in\P$, then it is achievable whenever $\P$ and $\Q$ are disjoint. More generally, even when $\P$ and $\Q$ are not polytopes, we characterize the existence of a bounded nontrivial e-variable whose expectation exactly equals one under any $P \in \P$. The proofs utilize recently developed techniques in simultaneous optimal  transport.  A key role is played by coarsening the filtration: sometimes, no such p-value or e-value exists in the richest data filtration, but it does exist in some reduced filtration, and our work provides the first general characterization of this phenomenon. We also provide an iterative construction that explicitly constructs such processes, and under certain conditions it finds the one that grows fastest under a specific alternative $Q$. We discuss implications for the construction of composite nonnegative (super)martingales, and end with some conjectures and open problems.
\end{abstract}

\begin{keyword}[class=MSC]
\kwd[Primary ]{62B15}
\kwd{62G10}
\kwd{49Q22}
\kwd{60G42}
\end{keyword}

\begin{keyword}
\kwd{p-values}
\kwd{e-values}
\kwd{composite hypothesis testing}
\kwd{convex order}
\kwd{simultaneous optimal transport}
\kwd{nonnegative martingale}
\end{keyword}

\end{frontmatter}
%%%%%%%%%%%%%%%%%%%%%%%%%%%%%%%%%%%%%%%%%%%%%%
%% Please use \tableofcontents for articles %%
%% with 50 pages and more                   %%
%%%%%%%%%%%%%%%%%%%%%%%%%%%%%%%%%%%%%%%%%%%%%%
%\tableofcontents

\section{Introduction}
\label{sec:intro}

Consider a universe of distributions $\Pi$ on 
a sample space $(\X,\mathcal F)$, where $\X$ is a Polish space.
%a filtered probability space $(\Omega,\F)$. 
% One draw from any distribution $P$ in $\Pi$ results in a sequence $X_1,X_2,\dots$ of observations.
The data are generated according to some $\mathbb{P} \in \Pi$.
% \com{assuming we can move this to later} 
Let $\P$ and $\Q$ be disjoint  subsets of $\Pi$.
When we say we are testing $\P$, we mean that we are testing the null hypothesis $\mathbb{P} \in \P$.  When we say we are testing $\P$ against $\Q$, we mean  additionally that the alternative hypothesis  is  $\mathbb{P} \in \Q$.

We ask (and answer) several central questions in this paper. The first one is:
\begin{quote}
\normalsize
     \textbf{(Q-exact-p).}  Given a null $\P$ and an alternative $\Q$, when can we find an \emph{exact} p-value for $\P$ that has nontrivial power under $\Q$? To elaborate, we would like to find a $[0,1]$-valued random variable $T$ that is exactly uniform for every $P\in\P$, but is stochastically smaller than uniform under every $Q \in \Q$.
\end{quote}

\noindent The second central question in this paper is the following:
% We focus mainly on the
% one-step version of the above problem, which is formulated as follows:
% with $\mathcal P$ and $\mathcal Q$ representing the distributions of one observation: 
\begin{quote}\normalsize
     \textbf{(Q-exact-e).} Given a null $\P$ and an alternative $\Q$, when does there exist an \emph{exact} e-value for $\P$ that has nontrivial power under $\Q$? To elaborate, we would like to find a nonnegative random variable $X$ such that $\E^P[X] = 1$ for every $P\in\P$, but $\E^Q[\log X] > 0$ (or $\E^Q[X] > 1$) for every $Q \in \Q$.
\end{quote}
% \com{Shall we highlight the word exact here?}
% An answer to the above question will suffice for an answer to the original question. 
We will provide a   complete answer to both questions in this paper, when $\P$ and $\Q$ are convex polytopes in the space of probability measures on $\X$. The solution is surprisingly clean and will be explained soon below.

%: we prove that constructions for \textbf{(Q-exact-p)} and \textbf{(Q-exact-e)} exist if and only if  $\Q$ does not intersect the span of $\P$. \com{This sentence needs to be revised in accordance with the contributions paragraphs.}

We also answer the non-exact versions of both problems, where we only require the p-value $T$ to be stochastically larger than uniform under any $P \in \P$:
\begin{quote}\normalsize
     \textbf{(Q-general-p).} Given a null $\P$ and an alternative $\Q$, when does there exist a  p-value for $\P$ that has nontrivial power against $\Q$? (see Terminology below.)
\end{quote}
Or, for the e-value, we require that $\E^P[X]\leq 1$ for any $P \in \P$:
\begin{quote}\normalsize
     \textbf{(Q-general-e).} Given a null $\P$ and an alternative $\Q$, when does there exist an  e-value for $\P$ that has nontrivial ``e-power'' against $\Q$? (see Terminology below.)
\end{quote}
For these non-exact problems, we can still provide a clean characterization of the existence for both  \textbf{(Q-general-p)} and \textbf{(Q-general-e)}.
%the solutions are still very clean: constructions for \textbf{(Q-general-p)} and \textbf{(Q-general-e)} exist if and only if $\P$ and $\Q$ are disjoint.\com{This sentence needs to be revised in accordance with the contributions paragraphs.}
%  {In the following, we use \textbf{(Q-existence)} to denote the set of the four questions \textbf{(Q-exact-p)}-\textbf{(Q-general-e)} above.}
An immediate follow-up question is:
\begin{quote}\normalsize
     \textbf{(Q-power).}  Suppose that we know the p-values or e-values in the above questions do exist. How can we algorithmically construct  powerful, or even optimal, ones?
\end{quote} 
 This question is important for the application of our ideas in hypothesis testing. 

These appear to be rather fundamental questions, and our answers will be proved using recent techniques in simultaneous optimal transport, combined with classical convex geometric arguments. 
 {A natural motivation for exactness of p-values and e-values comes from the trivial observation that, in the case of a simple null hypothesis,
any non-exact p-value or e-value can be strictly improved. Although this is not necessarily true for composite hypotheses, the existence of such exact p-values and e-values, as well as the trade-off between exactness and power, is useful for the design of tests.
}

Note that in the characterizations for \textbf{(Q-exact-p)} and \textbf{(Q-general-p)}  above, a technical condition of joint non-atomicity will be assumed, which is essentially  equivalent to allowing for external randomization. Our proofs are constructive and yield a simple iterative construction addressing \textbf{(Q-power)}, called SHINE (Separating Hyperplanes Iteration for Nontrivial and Exact e/p-variables), that can in principle explicitly build these objects and calculate their values on a given dataset, but it is only computationally feasible for low-dimensional settings.

Towards the end of the paper, we show how answers to the above two questions help answer a final related question:
\begin{quote}\normalsize
    \textbf{(Q-martingale).} Given a null $\P$ and an alternative $\Q$, can we determine if there is a nonnegative (super)martingale $M$ for $\P^\infty$ that grows to infinity under $\Q^\infty$? In other words, when can we find a process $M$ that is a nonnegative (super)martingale  under $P^\infty$ simultaneously for every $P\in\P$, but it almost surely grows to infinity under $Q^\infty$ for every $Q\in\Q$?
\end{quote}

% In what follows, we call p-values, e-values and nonnegative supermartingales as ``objects'' of statistical evidence.
% In short, such objects (that fulfill their respective criteria) exist if and only if the convex hull of $\Q$ does not intersect the span (or convex hull) of $\P$. 

Before proceeding, we introduce  important terminology used throughout  the paper.

\smallskip

% \subsection{Terminology}
\noindent \textbf{Terminology.}
We define pivotal, exact, and nontrivial e- and p-variables below.
% We collect the most important terminology on p-variables and e-variables used throughout the paper. 

% For simplicity in this paper, we often consider the ``iid setting'', where the distributions in $\P,\Q$ are of the form $P^\infty$, meaning that we see independent and identically distributed observations from $P$. We say $P \in \P$ as shorthand for $P^\infty \in \P$, in which case we say that the null is true, and if $P \in \Q$, we say that the null s false and the alternative is true. 

%We assume that $\X$ is a Polish space. The collections $\P,\Q$ of probability measures on $(\X,\mathcal F)$ represent a  null hypothesis 
%$\mathbb{P} \in \P$, 
%and an alternative hypothesis 
%$\mathbb{P} \in \Q$, respectively, where we recall that $\mathbb P$ generated the data.
% In this case, we say we are testing $\P$ against $\Q$.  

\begin{enumerate} 
\item A random variable
$X$ is \emph{pivotal} for $\P$
if $X$ has the same distribution under all $P\in \P$.

\item A nonnegative random variable 
$X$ is a \emph{e-variable} for $\P$
if $\E^P[X]\le 1$ for all $P\in \P$.
An e-variable $X$ for $\P$ is \emph{exact}  
if $\E^P[X]= 1$ for all $P\in \P$. We say $X$ is \emph{nontrivial} for $\Q$ if $\E^Q[X]>1$ for all $Q\in\Q$.  
 An e-variable $X$ for $\P$ is said to have \emph{nontrivial e-power}  against $\Q$ if 
  for each $Q\in\Q$, $\E^Q[\log X]>0$.

\item A nonnegative random variable 
$X$ is a \emph{p-variable} for $\P$
if $P(X\le \alpha)\le \alpha$ for all  $\alpha \in (0,1)$ and  $P\in \P$, and a p-variable 
$X$ is \emph{exact}
if $P(X\le \alpha)=  \alpha$ for all  $\alpha \in (0,1)$ and  $P\in \P$.
A p-variable $X$ for $\P$ is \emph{nontrivial} (or has \emph{nontrivial power}) against $\Q$ if, for each $Q\in \Q$, 
 $Q(X\le \alpha)\ge \alpha $ for all $\alpha\in (0,1) $ with strict inequality for some $\alpha \in (0,1)$. Without loss of generality, p-variables can be restricted to $[0,1]$ by truncation, without changing their properties.
 \end{enumerate}

Note that an exact p-variable is always pivotal, but not vice versa. An exact e-variable need not be pivotal, and a pivotal e-variable need not be exact. Since $x-1\geq\log x$, an e-variable that has nontrivial e-power against $\Q$ is also nontrivial for $\Q$. 
 We will often omit $\P$ and $\Q$ in our subsequent mentions of p/e-variables when they are clear from the context. 
Realizations of e-variables are called e-values. 
Like many other authors, we do not distinguish these terms when there is no confusion; the same applies to p-values and p-variables.

\begin{remark}
    For the majority of this paper, we suppress the raw data that is observed and used to form the p-values or e-values. One may simply assume that we have observed one data point $Z$ from $\mathbb P$. This $Z$ could itself be a random vector of some size $n \geq 1$ lying in (say) $\mathbb R^d$ for some $d \geq 1$ (which means $\mathbb P$ may be $\mu^n$ for some $\mu$ on $\mathbb R^d$), but we leave all this implicit. Thus our p-values and e-values can be treated as ``single-period'' statistics calculated on a batch of data. We return to the multi-period (sequential) case briefly later in the paper.
\end{remark}

\textbf{Summary of contributions.}
We briefly summarize the main results of this paper below. 
With the help of techniques from simultaneous transport, the existence of  p/e-values for a convex polytope $\P$ and a simple alternative $\Q=\{Q\}$ is fully characterized in Theorems \ref{thm:iff} and \ref{thm:p-variable}: under a natural condition of non-atomicity, we show that
pivotal, exact, and powerful p/e-values exist if and only if
$Q\not \in \spn \P$; 
powerful p/e-values exist if and only if 
$Q\not \in \P$. 
% These two results provide partial answers to \textbf{(Q-existence)} in case $\Q$ is simple.
Theorems \ref{thm:iff2} and \ref{thm:pvariable2} extend these earlier results to the case of composite alternatives that are polytopes:
for convex polytopes $\P$ and $\Q$, similar conclusions as before hold  with the condition  $Q\not \in \spn \P$  being replaced by 
$\spn\P\cap\Q=\emptyset$,
and the condition
 $Q\not \in \P$ being replaced by 
$\P\cap\Q=\emptyset$.
Theorem \ref{thm:infinite case} extends these results to the case of general (non-polytope) infinite $\P,\Q$, where the situation is more complicated: we now additionally need a common reference measure and a closure with respect to the total variation distance. 

For the particular  case of a simple alternative ($\Q=\{Q\}$), we can speak of maximizing the e-power under $Q$ among all exact e-variables. The  exact e-variable with the largest e-power is studied in a series of results including Theorems \ref{thm:maximumexist} and \ref{theorem:cx}, and this finally leads to  the SHINE construction, with maximality of the constructed e-variable shown in Theorem \ref{thm:maximal}, providing an answer to \textbf{(Q-power)}.

Finally, the above results directly give rise to an answer to  \textbf{(Q-martingale)} by obtaining sufficient conditions for the existence of a powerful e-process   (Corollary \ref{coro:testm}).\\
% that leads to an explicit e-variable that has considerably large e-power.  

% \red{Two caveats: 
% \begin{itemize}
%     \item At some point I also call it an exact p-variable if it's $\lst\Leb$ under $\P$ but $=\Leb$ under $\Q$. I think there's no difference  (at least for atomless)
%     \item In the proof of Theorem~\ref{thm:iff2} I decided that the null and alternative is flipped...
% \end{itemize}}

% \paragraph{Our contributions.}[TBD]

\textbf{Related results.} The most directly related work is that of \cite{GrunwaldHK19}, which focuses primarily on e-values, and in particular \textbf{(Q-general-e)}. To paraphrase one of their main results, consider any $\P$ and $\Q$ with a common reference measure, whose convex hulls do not intersect. They show that as long as a particular ``worst case prior'' exists, then one can construct an e-value for $\P$ which maximizes the worst case e-power for $\Q$. This is a topic we return to later in the paper, when we provide a more detailed geometric study of \textbf{(Q-general-p)} and \textbf{(Q-general-e)} together. We need fewer technical conditions to establish our results, but their additional assumptions allow them to handle general $\P,\Q$ that are not polytopes. See also \cite{harremoes2023universal} for a very recent follow-up work by the same group, which relaxes some of the original technical conditions.

A second related work is that of \cite{ramdas2022testing}. Here, the authors work in the sequential setting and ask when nontrivial nonnegative (super)martingales for $\P^\infty:=\{P^\infty: P \in \P\}$ exist. We can paraphrase their geometric solution: assuming a common reference measure, nontrivial nonnegative (super)martingales cannot exist if the ``fork-convex hull'' of $\P^\infty$ intersects $\Q^\infty$.

The above papers hint at a deeper underlying geometric picture, and our work elaborates significantly on this theme, completely characterizing the case of convex polytopes. One key point is that the earlier works did not give a systematic and thorough treatment of what one can accomplish in reduced filtrations, while this is a central aspect of our paper. Informally, we will (optimally) transport $\P$ to a single measure $\mu$, while transporting $\Q$ to a single measure $\nu \neq \mu$, and this collapse of the null and alternative corresponds exactly to working in a coarser $\sigma$-algebra. 

The above idea of transport from multiple measures to specified measures is addressed in the  framework of simultaneous transport studied by \cite{wang2022simultaneous}. We borrow several techniques from their work and build on them significantly to provide answers to  our questions. In particular, our work tightly connects arguably basic testing problems with the modern theory of optimal transport.

A third classical yet fundamental related work is  Kraft's theorem (\citet[Theorem 5]{kraft1955some}), which states that if there is a $\sigma$-finite  reference measure $R$ that dominates every distribution in $\P\cup\Q$, then for each $\ee>0$ there exists a $[0,1]$-valued random variable $X$ with 
    \begin{align}
        \inf_{Q\in\Q}\E^Q[X]\geq\ee+\sup_{P\in\P}\E^P[X]\label{eq:kraft}
    \end{align}if and only if the total variation distance $d_{\mathrm{TV}}(\conv\P,\conv\Q)\geq\ee$. Kraft's theorem serves as a starting point for distinguishing sets of distributions (\cite{hoeffding1958distinguishability}) and impossible inference (\cite{bertanha2020impossible}). In particular, in Remark \ref{remark:kraft} below we will see how Kraft's theorem can answer \textbf{(Q-general-e)} above.

Finally, likelihood ratios play an intimate role throughout our paper, but in rather different ways than classical hypothesis testing results. For general composite nulls and alternatives,  generalized likelihood ratio-based methods require certain regularity conditions  in order for Wilks' theorem~\citep{wilks1938large} to apply, which in turn yields an asymptotically exact p-value. An alternative, recent e-value approach is taken by universal inference~\citep{wasserman2020universal}. Our paper takes a very different approach, designing non-asymptotically exact p-values \textbf{(Q-exact-p)} or non-asymptotically conservative p-values \textbf{(Q-general-p)}, and also doing the same for e-values \textbf{(Q-exact-e, Q-general-e)}. We do not impose the regularity conditions required for Wilks' theorem to hold  (our assumptions are different and quite mild), and we are interested in when such p-values or e-values exist and how one can construct them \textbf{(Q-power)}. As a  rough, but instructive, intuition for how likelihood ratios play a role in our work, when deriving exact p-values or e-values, our method tries to find a transport map that can simultaneously transport the entire composite $\P$ into a single uniform $U$, while simultaneously transporting the composite $\Q$ into some distribution $F \neq U$. Now, having effectively converted the given composite problem into a point null $U$ and a point alternative $F$, one can use simple likelihood ratios to design either the p-values or e-values.\\

\textbf{Background on e-values.}
  E-values  are an alternative to p-values, and they have recently been actively studied in statistical testing by   \cite{wasserman2020universal}, \cite{Shafer:2021}, \cite{Vovk/Wang:2021}, \cite{GrunwaldHK19}, and \cite{howard2021time} under various names.  Tests based on e-values are closely related to nonnegative supermartingale techniques for testing and estimation, which date back to work by Robbins~\citep{Darling/Robbins:1967b,robbins1974expected},
  and they emphasize continuous monitoring, optional stopping or continuation of experiments. 
  The notion of e-processes generalizes that of  likelihood ratios to composite hypotheses~\citep{ramdas2022testing}. 
  Some advantages of testing with e-values are summarized in  \citet[Section 2]{wang2022false}.
  The idea of testing with e-values is intimately connected to game-theoretic probability~\citep{Shafer/Vovk:2001,Shafer/Vovk:2019}.  
   For a recent review on e-values and game-theoretic statistics, see \cite{ramdas2022game}.\\

% and that if $\P$ is ``locally dominated'', then \emph{admissible} e-processes must satisfy 
% \begin{equation}\label{eq:einf}
%   E_t = \inf_{P \in \P} M^P_t
% \end{equation}
% for some test martingale family $(M^P)_{P \in \P}$.
% (Technically, the $\inf$ above is an ``essential infimum''.)

\textbf{Notation.} We collect the notation we use  throughout this paper.

\begin{enumerate}
    \item \emph{Topology}. For a set $A\subseteq\R^d$, $A^\circ$ (resp.~$\overline{A},~\partial A,~A^c,~\conv A$) is the interior (resp.~closure, boundary, complement, convex hull) of $A$ and $\aff A$ is the smallest affine subspace of $\R^d$ containing $A$. For an affine subspace $S\subseteq \R^d$, we denote by $\ri (A;S)$ is the relative interior of $A$ in $S$, that is, the interior of $A$ in the relative topology on $S$. 
    
    \item \emph{Probability and measure}. All measures we consider will be finite and have a finite first moment, i.e., $\int|x|\mu(\d x)< \infty$. For a Polish space $\X$, we let $\mathcal{M}(\X)$ be the set of all finite measures on $\X$ and $\Pi(\X)$ be the set of probability measures on $\X$. For $\mu\in\cM(\R^d)$, we denote its barycenter by $\bary(\mu):=\int_{\R^d} x\,\mu(\d x)/\mu(\R^d)$. For a finite set $\mathcal A$ of random variables or probability measures on the same space, we define $\conv \mathcal A$ and $\spn\mathcal A$ in the usual sense of convex hull and span. We write $X\lawis_P\mu$, or simply $X\lawis\mu$, if the random variable $X$ has distribution $\mu$ under $P$.  We say ``a probability measure $\mu$ is supported on a set $A$'' if $\mu(A)=1$. This does not imply that $A$ is closed or $A=\supp\mu$.   The product measure is denoted by $P\otimes Q$. If $\P=\{P_1,\dots,P_L\}$ and $\Q=\{Q_1,\dots,Q_M\}$ are two sets of probability measures on $\X$, we sometimes denote the tuple $(P_1,\dots,P_L,Q_1,\dots,Q_M)$ by $(\P,\Q)$. For $P,Q\in\cM(\X)$ we write $P\ll Q$ if $P$ is absolutely continuous with respect to $Q$ (sometimes we say $Q$ dominates $P$), and $P\approx Q$ if $P\ll Q\ll P$.

    \item \emph{Stochastic orders}. For $F,G\in\Pi(\R)$, we write $F\lst G$ if $F((-\infty,a])\geq G((-\infty,a])$ for all $a\in\R$. Also, 
$F\pst G$ if $F\lst G$ and $F\neq G$.  For $\mu,\nu\in\cM(\R^d)$, we denote by $\mu\lcx\nu$ if $\int\phi\,\d\mu\leq\int\phi\,\d\nu$ for every convex function $\phi$, in which case we say $\mu$ is smaller than $\nu$ in convex order.\footnote{This is sometimes  called the Choquet order in the mathematical literature, e.g.,  \cite{simon2011convexity}.}  If $\mu,\nu$ are probability measures and $X\lawis\mu$, $Y\lawis\nu$, we sometimes abuse notation and write $X\lcx Y$ instead of $\mu\lcx\nu$. We  write $\mu\leq \nu$ if $\mu(A)\leq \nu(A)$ for every Borel set $A$. 

\item \emph{Other notation}. Bold symbols such as $\bx$ and $ \ba$ will typically denote  vectors.  Write $\e_d=(1,\dots,1)\in\R^d$, $\z_d=(0,\dots,0)\in\R^d$, $\I_d=\{\bx\in\R^d\mid x_1=\dots=x_d\}=\R\e_d$, and $\I_d^+=\{\bx\in\R^d\mid x_1=\dots=x_d\geq 0\}=\R_+\e_d$. When the dimension $d$ is clear, we may omit the subscript $d$ and write $\e,\z,\I,\I^+$ instead. We let  $\Leb$ denote the Lebesgue measure on $[0,1]$. Denote the Euclidean norm by $\n{\cdot}$. 
\end{enumerate}

\textbf{Outline of the paper.} The rest of this paper is organized as follows. Section~\ref{sec:prelim} provides the necessary mathematical background regarding convex order and simultaneous optimal transport. The easier case with a simple alternative ($|\Q|=1$) will be solved first in Section~\ref{sec:simple alt}. Under suitable conditions, we solve the maximization problem of the e-power in Section~\ref{sec:e} and illustrate the SHINE construction for finding a  powerful e-variable in Section~\ref{sec:alg} for a simple alternative, thus answering \textbf{(Q-power)}. We answer \textbf{(Q-exact-p)}, \textbf{(Q-exact-e)}, \textbf{(Q-general-p)} and \textbf{(Q-general-e)} in full  in Section~\ref{sec:composite alt}, where we deal with a  general composite (and even infinite) alternative $\Q$. Finally, an application to composite test (super)martingales related to \textbf{(Q-martingale)}  will be discussed in Section~\ref{sec:testmartingale}, followed by a summary in Section~\ref{sec:summary}.   
Appendix \ref{sec:general} contains a few general results on the existence of p/e-values, followed by proofs of our main results in Appendices~\ref{App:B}-\ref{app:sec6proof}.  Appendix \ref{appendix} contains a few technical results that are used in our proofs. 

\section{Preliminaries on convex order and simultaneous transport}\label{sec:prelim}

In this section, we collect results related to convex order and simultaneous transport for future use. We rely on some results from  \cite{shaked2007stochastic} and \cite{wang2022simultaneous}. 

 {In the setting of classical optimal transport theory, one usually starts with two measures $\mu\in\Pi(\X),\nu\in\Pi(\Y)$ on Polish spaces $\X,\Y$, and a typical goal would be optimizing a certain functional over  $(X,Y)$ with respective marginals $\mu,\nu$ (such $(X,Y)$ are called couplings). The set of such couplings is also referred to as \emph{transport plans}. In certain cases, one is interested in a special class of transport plans where $Y$ is required to be a function of $X$. Such couplings are called \emph{transport maps}. See \cite{santambrogio2015optimal} and \cite{villani2009optimal} for background on optimal transport.}

 {A coupling $(X,Y)$ on $\R^d\times\R^d$ is called a \emph{martingale coupling} if $\E[Y| X]=X$. Given $\mu,\nu\in\Pi(\R^d)$, a \emph{martingale transport (plan)} from $\mu$ to $\nu$ is a martingale coupling $(X,Y)$ such that $X\lawis\mu$ and $Y\lawis\nu$. } We recall from \cite{strassen1965existence} that there exists a martingale transport from $\mu$ to $\nu$ if and only if $\mu\lcx\nu$  {(see point 3 in the notation subsection for a definition)}. This result is called Strassen's theorem. The relation $\lcx$ is  a partial order on  $\Pi(\R^d)$. Given a subset $\mathcal N\subseteq\Pi(\R^d)$,  we say $\mu$ is a (Pareto) \emph{maximal} element in $\mathcal N$ if there exists no $\nu\in\mathcal N$ such that $\nu\neq\mu$ and $\mu\lcx\nu$; we say $\mu$ is the \emph{maximum} element in $\mathcal N$ if $\nu\lcx\mu$ for each $\nu\in\mathcal N$.
% In the following, we collect a few properties of the convex order. 
These next facts can be found in \citet[Section 3.A]{shaked2007stochastic}.

\begin{lemma}\label{lemma:cxproperty}
    The followings hold for all integrable real-valued random variables:
    \begin{enumerate}[(i)]
    \item If $\E[X]=\E[Y]$, then $X\lcx Y$ if and only if 
    $\E[(X-a)_+]\leq\E[(Y-a)_+]\text{ for all }a\in\R.$
        \item If $\{X_n\}$ is a sequence of random variables that converge weakly to $X$ and $\E[|X_n|]\to\E[|X|]$, then 
        $X_n\lcx Y\ \implies\ X\lcx Y.$
    \end{enumerate}
\end{lemma}

 {The recent work of \cite{wang2022simultaneous} proposed the notion of simultaneous optimal transport as an extension of classical optimal transport. As explained above, classical optimal transport theory concerns a coupling between two measures. In the setting of simultaneous optimal transport, one starts from two $d$-tuples of probability measures $\bmu=(\mu_1,\dots,\mu_d)$ on $\X$ and $\bnu=(\nu_1,\dots,\nu_d)$ on $\Y$, and requires that the transport plan (or map) sends $\mu_j$ to $\nu_j$ \emph{simultaneously} for all $1\leq j\leq d$. If $d=1$, this coincides with the classical optimal transport.}

 {Let us give a formal definition.} For $d\geq 1$ and two $\R^d$-valued measures $\bmu,\bnu$ on Polish spaces $\X,\Y$ (denoted by $\bmu\in\mathcal{M}(\X)^d$ and $\bnu\in\mathcal{M}(\Y)^d$) such that $\bmu(\X)=\bnu(\Y)$, let $\K(\bmu,\bnu)$ and $\T(\bmu,\bnu)$  denote the set of all \emph{simultaneous transport plans} and \emph{maps} from $\bmu$ to $\bnu$ respectively, i.e., $\K(\bmu,\bnu)$ is the set of all stochastic kernels $\kappa$ such that
$$\kappa_\#\bmu(\cdot):=\int_\X\kappa(x;\cdot)\bmu(\d x)=\bnu(\cdot),$$
and 
$$\T(\bmu,\bnu)=\{T:\X\to\Y \mid \bmu\circ T^{-1}=\bnu\}.$$
When $d=1$, $\K(\bmu,\bnu)$ is often represented as the set of all joint distributions on $\X \times \Y$ whose marginals are $\bmu$ and $\bnu$ respectively, but for $d>1$, we prefer the above representation. 
 {The mathematical structure of simultaneous optimal transport is very different from classical optimal transport, and the existence of simultaneous transport plans (or maps) is  a non-trivial question.}
To further characterize the existence of simultaneous  transport maps and plans, we need the notion of joint non-atomicity.

\begin{definition}\label{def:ssww}
    Consider a tuple of probability measures $\bmu=(\mu_1,\dots,\mu_d)$  on a Polish space $\X$. We say that  $\bmu$  is \emph{jointly atomless} if there exists   $\mu \gg \sum_{i=1}^d\mu_i$ and  a  random variable $\xi $ such that  under $\mu$, $\xi$ is atomless and independent of $(\d \mu_1 / \d \mu,\dots,\d \mu_d/\d \mu) $. 
\end{definition}

 As a simple example, $(\mu_1\times\Leb,\dots,\mu_d\times\Leb)$ on $\X\times [0,1]$ is jointly atomless for each collection $(\mu_1,\dots,\mu_d)$ on $\X$. We refer to \cite{shen2019distributional} and \cite{wang2022simultaneous} for more discussions on this notion.

In statistical terms, the hypothesis  $\{P_1,\dots,P_L\}$ as a tuple being jointly atomless is equivalent to allowing for additional randomization, i.e., simulating a uniform random variable independent of the Radon--Nikodym derivatives $(\d P_1/\d P,\dots,\d P_L/\d P)$ for some $P\in \Pi(\X)$. 
It suffices if simulating a uniform random variable independent of existing random variables is always allowed. Such an assumption is common in statistical methods based on resampling or data splitting.

\begin{proposition}
\label{prop:ssww}
Consider $\bmu\in\Pi(\X)^d$ and $\bnu \in \Pi(\Y)^d$. Let $\bl\in \R_+^d$ satisfy $\n{\bl}_1=1$,  and define $\mu:=\bl^\top \bmu,$ and $\nu:=\bl^\top\bnu$. Assume that $\mu_j\ll\mu$ and $\nu_j \ll \nu$ for each $1\leq j\leq d$. Then,
\begin{enumerate}[(i)]
\item The set  $\mathcal K(\bmu,\bnu)$ is non-empty if and only if 
 $$\left(\frac{\d \mu_1}{\d \mu},\dots,\frac{\d \mu_d}{\d \mu}\right)\Big|_{\mu}\gcx \left(\frac{\d \nu_1}{\d \nu},\dots,\frac{\d \nu_d}{\d \nu}\right)\Big|_{\nu},$$
 where $X|_P$ means the distribution of a random variable $X$ under a measure $P$.
 %on the left-hand side, the $|_\mu$ means the Radon--Nikodym derivative is considered as a random variable under the measure $\mu$. 
 \item Assume that $\bmu$ is jointly atomless.  
The set $\mathcal T(\bmu,\bnu)$ is non-empty  if and only if $$\left(\frac{\d \mu_1}{\d \mu},\dots,\frac{\d \mu_d}{\d \mu}\right)\Big|_{\mu}\gcx \left(\frac{\d \nu_1}{\d \nu},\dots,\frac{\d \nu_d}{\d \nu}\right)\Big|_{\nu}.$$ 
\end{enumerate}  
\end{proposition}

\begin{proof}
\sloppy Theorem 3.4 of  \cite{wang2022simultaneous} implies that the statements hold with $\bl=(1/d,\dots,1/d)$. The more general case follows from Lemma 3.5 of \cite{shen2019distributional}, in the direction (iii)$\Rightarrow$(ii) there.
\end{proof}

We briefly describe the intuition behind this result, which is crucial for our paper. In the sequel, a coupling $(X,Y)$ is \emph{backward martingale} if $\E[X|Y]=Y$; that is, $(Y,X)$ forms a martingale. It is \emph{Monge} if $Y$ is a measurable function of $X$.  The key observation is that the pushforward $\kappa_\#\bmu$ mixes the ratios between different coordinates of the (vector-valued) masses of $\bmu$ at different places of $\X$; see Figure~\ref{fig:sot}. The ``ratios'' can be recognized as Radon--Nikodym derivatives. The ``mix'' effect can be interpreted as a backward martingale transport, because reversing the transport arrows (or equivalently, looking at the transport in the backward direction)  gives rise to  a martingale coupling of the Radon--Nikodym derivatives.  Strassen's theorem then gives the convex order constraint on the  Radon--Nikodym derivatives. In \cite{wang2022simultaneous}, such an observation leads also to the  MOT-SOT\footnote{Here, MOT stands for martingale optimal transport, and SOT stands for simultaneous optimal transport.} parity that relates the simultaneous transport to the underlying backward martingale transport, which will be useful for our purpose when constructing explicitly an e/p-variable. We state a weak form of the MOT-SOT parity below, which can be proved similarly to Corollary 3 of  \cite{wang2022simultaneous}.

\begin{figure}[!t] 
\centering
\begin{tikzpicture} 

\draw[gray, very thick] (-4.9,0)--(0.9,0);
\draw[gray, very thick] (-4.9,-3.9)--(0.9,-3.9);
\filldraw[blue, very thick] (-4.6,0) rectangle (-4.4,0.3);
%\filldraw[red,ultra thick] (-4.6,0.3) rectangle (-4.4,1);
\filldraw[blue, very thick] (-3.6,0) rectangle (-3.4,0.8);
%\filldraw[red,ultra thick] (-3.6,0.8) rectangle (-3.4,1);
\filldraw[blue, very thick] (-2.6,0) rectangle (-2.4,0.7);
%\filldraw[red,ultra thick] (-2.6,0.7) rectangle (-2.4,1);
\filldraw[blue, very thick] (-1.6,0) rectangle (-1.4,0.1);
%\filldraw[red,ultra thick] (-1.6,0.1) rectangle (-1.4,1);
\filldraw[blue, very thick] (-0.6,0) rectangle (-0.4,0.5);
%\filldraw[red,ultra thick] (-0.6,0.5) rectangle (-0.4,1);
\filldraw[blue, very thick] (0.4,0) rectangle (0.6,0.6);
%\filldraw[red,ultra thick] (0.4,0.6) rectangle (0.6,1);

\filldraw[blue, very thick] (-3.2,-3.9) rectangle (-2.9,-2.8);
%\filldraw[red,ultra thick] (-3.2,-3.9) rectangle (-2.9,-2);
\filldraw[blue, very thick] (-1.1,-3.9) rectangle (-0.8,-2);
%\filldraw[red,ultra thick] (-1.1,-3.1) rectangle (-0.8,-2);

\draw[very thick, ->](-4.5,-0.1)--(-3.3,-1.9);
\draw[very thick, ->](-3.5,-0.1)--(-1.2,-1.9);
\draw[very thick, ->](-2.5,-0.1)--(-3,-1.9);
\draw[very thick, ->](-1.5,-0.1)--(-2.8,-1.9);
\draw[very thick, ->](-0.5,-0.1)--(-0.9,-1.9);
\draw[very thick, ->](0.5,-0.1)--(-0.7,-1.9);

\draw[gray, very thick] (4.1-1.6,0)--(9.9-1.6,0);
\draw[gray, very thick] (4.1-1.6,-3.9)--(9.9-1.6,-3.9);
%\filldraw[blue, very thick] (-4.6+7.4,0) rectangle (-4.4+7.4,0.3);
\filldraw[red,ultra thick] (-4.6+7.4,0) rectangle (-4.4+7.4,0.7);
%\filldraw[blue, very thick] (-3.6+7.4,0) rectangle (-3.4+7.4,0.8);
\filldraw[red,ultra thick] (-3.6+7.4,0) rectangle (-3.4+7.4,0.2);
%\filldraw[blue, very thick] (-2.6+7.4,0) rectangle (-2.4+7.4,0.7);
\filldraw[red,ultra thick] (-2.6+7.4,0) rectangle (-2.4+7.4,0.3);
%\filldraw[blue, very thick] (-1.6+7.4,0) rectangle (-1.4+7.4,0.1);
\filldraw[red,ultra thick] (-1.6+7.4,0) rectangle (-1.4+7.4,0.9);
%\filldraw[blue, very thick] (-0.6+7.4,0) rectangle (-0.4+7.4,0.5);
\filldraw[red,ultra thick] (-0.6+7.4,0) rectangle (-0.4+7.4,0.5);
%\filldraw[blue, very thick] (0.4+7.4,0) rectangle (0.6+7.4,0.6);
\filldraw[red,ultra thick] (0.4+7.4,0) rectangle (0.6+7.4,0.4);

%\filldraw[blue, very thick] (-3.2+7.4,-5) rectangle (-2.9+7.4,-3.9);
\filldraw[red,ultra thick] (-3.2+7.4,-3.9) rectangle (-2.9+7.4,-2);
%\filldraw[blue, very thick] (-1.1+7.4,-5) rectangle (-0.8+7.4,-3.1);
\filldraw[red,ultra thick] (-1.1+7.4,-3.9) rectangle (-0.8+7.4,-2.8);

\draw[very thick, ->](-4.5+7.4,-0.1)--(-3.3+7.4,-1.9);
\draw[very thick, ->](-3.5+7.4,-0.1)--(-1.2+7.4,-1.9);
\draw[very thick, ->](-2.5+7.4,-0.1)--(-3+7.4,-1.9);
\draw[very thick, ->](-1.5+7.4,-0.1)--(-2.8+7.4,-1.9);
\draw[very thick, ->](-0.5+7.4,-0.1)--(-0.9+7.4,-1.9);
\draw[very thick, ->](0.5+7.4,-0.1)--(-0.7+7.4,-1.9);
\end{tikzpicture}

\caption{A showcase of simultaneous transport: here the input vector $\bmu$ is two-dimensional, as is the output vector $\bnu$. The two input distributions are discrete distributions over the same alphabet of size six and are drawn in different colors in the top row, with the height of a bar indicating its mass. The two target distributions are binary, indicated on the bottom row. The \emph{simultaneous} transport requires that the maps that transport from $\mu_1$ to $\nu_1$ (left) and from
$\mu_2$ to $\nu_2$ (right) are identical. This map is achieved by mixing (averaging) the Radon--Nikodym derivatives. Denoting $\bar\nu=\nu_1+\nu_2$ and $\bar \mu = \mu_1+\mu_2$, we have $\frac{\d\nu_1}{\d \bar \nu}(1) = \frac13 \left( \frac{\d\mu_1}{\d \bar \mu}(1) + \frac{\d\mu_1}{\d \bar \mu}(3) + \frac{\d\mu_1}{\d \bar \mu}(4)\right)$, and analogously for the other coordinate. 
% and $3 \frac{\d\mu_1}{\d \bar \mu}(2) = \frac{\d\mu_1}{\d \bar \mu}(2) + \frac{\d\mu_1}{\d \bar \mu}(5) + \frac{\d\mu_1}{\d \bar \mu}(6)$.
% his figure is an  from \citet[Figure 1]{wang2022simultaneous}. 
}
\label{fig:sot}
\end{figure}

\begin{proposition}\label{prop:commute}
    Let $\bmu\in\Pi(\X)^d$ and $\bnu \in \Pi(\Y)^d$ satisfy $\bmu\ll\mu_d$, $\bnu\ll\nu_d$ (where we recall that $\mu_d$ is the $d$-th component of the vector-valued measure $\bmu$), and $\K(\bmu,\bnu)$ non-empty.  Suppose that $\bmu$ is jointly atomless and $(\d\bmu/\d\mu_d)|_{\mu_d}$ is atomless. Then there exists a backward martingale coupling between $(\d\bmu/\d\mu_d)|_{\mu_d}$ and $(\d\bnu/\d\nu_d)|_{\nu_d}$ that is also Monge. Moreover, if we denote by $h$ the map that induces this Monge transport, then there exists a simultaneous transport map $T\in\T(\bmu,\bnu)$ satisfying
    $$\frac{\d\bnu}{\d\nu_d}(T(x))=h\left(\frac{\d\bmu}{\d \mu_d}(x)\right),~x\in\X.$$
\end{proposition}

In the above proposition, we have picked the $d$-th entry $\mu_d,\nu_d$ to evaluate the Radon--Nikodym derivatives. One could as well use $\mu_j,\nu_j$ for any $1\leq j\leq d$, or even $\bar \mu,\bar \nu$. When applying this result, we have in mind that the last entry of $\bmu,\bnu$ will be given by the alternative and the rest by the null, which makes it convenient to   evaluate the Radon--Nikodym derivatives using the $d$-th entry.

Finally, we recall the following  basic fact on Radon--Nikodym derivatives.

\begin{lemma}\label{lemma:RN}
    Let $d\in\N$ and $\tau$ be a probability measure supported on $\R_+^d$ with mean $\e$.  Then there exist probability measures $F_1,\dots,F_d$ supported on $[0,1]$ such that
    $$\left(\frac{\d F_1}{\d \Leb},\dots,\frac{\d F_d}{\d \Leb}\right)\Big|_{\Leb}= \tau.$$
\end{lemma}

\begin{proof}
Since $\Leb$ is atomless, $\T(\Leb,\tau)\neq\emptyset$.\footnote{It is a standard fact in optimal transport that a Monge transport map from $\mu$ to $\nu$ exists if $\mu$ is atomless.} Pick $(f_1,\dots,f_d)\in \T(\Leb,\tau)$, and define $F_i$ by $\d F_i/\d \Leb=f_i$ for $1\leq i\leq d$. This is well-defined since $f_i$ is nonnegative a.e.~and $\E^{\Leb}[f_i]=1$ for each $1\leq i\leq d$.    
\end{proof}

\section{Composite null and simple alternative}\label{sec:simple alt}
In this section, we characterize the existence of exact and pivotal p-variables and e-variables for composite null and simple alternative (singleton).
Although our results in this case are  covered by the more general result for composite alternatives treated in Section~\ref{sec:composite alt}, studying this setting  first helps with building intuition behind our proof techniques. Moreover, the concept of e-power studied in Section~\ref{sec:e} is defined for a single $Q$ in the alternative hypothesis.  
We fix $\P=\{P_1,\dots,P_L\}$ and $\Q=\{Q\}$ in $\Pi(\X)$ and will assume that
\begin{align}\label{JA}
    (\P,\Q)\text{ is jointly atomless,} \tag{JA}
\end{align}unless otherwise stated.
The main results are Theorems~\ref{thm:iff} and~\ref{thm:p-variable} below.  When \begin{align}\label{AC}
    P_1,\dots,P_L\ll Q \tag{AC}
\end{align}holds,  we define the measure $\gamma=(\d P_1/\d Q,\dots,\d P_L/\d Q)|_{Q}$ on $\mathbb R^L$.

% The following lemma gives a sufficient condition for the existence of a (pivotal and exact) e-variable that is nontrivial for $\Q$, that will be used multiple times.

\begin{theorem}\label{thm:iff}
Suppose that we are testing $\P=\{P_1,\dots,P_L\}$ against $\Q=\{Q\}$ and \eqref{JA} holds. The following are equivalent:
    \begin{enumerate}[(a)]\item there exists an exact (hence pivotal) and nontrivial  p-variable;
        \item there exists a pivotal, exact,  bounded e-variable that has nontrivial e-power against $\Q$;
\item there exists an exact  e-variable that is nontrivial against $\Q$;
\item there exists a random variable $X$ that is pivotal for $\mathcal P$ but has a different distribution under $Q$, where the laws of $X$ under both are atomless;
        \item it holds that $Q\not\in\spn(P_1,\dots,P_L)$.
    \end{enumerate}
\end{theorem}
To prove Theorem \ref{thm:iff} we need the following preparation.

\begin{lemma}\label{lemma:balls}
    Suppose that $Q\not\in\spn(P_1,\dots,P_L)$ and \eqref{AC} holds. There exists a disjoint collection of closed balls $B_1,\dots,B_{k}$ in $\R^L$ of positive measure (under $\gamma$) not containing $\e$ such that denoting by $t_j$ the point of $B_j$ closest to $\e$, we have $\e\in\conv(\{t_1,\dots,t_k\})^\circ$. 
\end{lemma}

\begin{proof}
Since $Q\not\in\spn(P_1,\dots,P_L)$, the measure $\gamma$ cannot have support contained in a hyperplane in $\R^L$ by definition. In other words, $\aff\supp \gamma=\R^L.$ By Lemma~\ref{lemma:elem}(ii), $\e=\bary(\gamma)\in (\conv\supp\gamma)^\circ$. Therefore, there exist $s_1,\dots,s_k\in\supp\gamma$ such that $\e\in(\conv\{s_1,\dots,s_k\})^\circ$. Let $B_j$ be the ball centered at $s_j$ with radius $r>0$ for $1\leq j\leq k$. For $r$ small enough, these balls will be disjoint from $\e$, and the closest points $t_1,\dots,t_k$ satisfy $\e\in\conv(\{t_1,\dots,t_k\})^\circ$.
\end{proof}

\begin{proposition}\label{prop:existence}
We have $Q\not\in \spn(P_1,\dots,P_L)$ if and only if there exist probability measures $G\neq F$ such that
$$\K((P_1,\dots,P_L,Q),(F,\dots,F,G))\neq\emptyset.$$
If moreover \eqref{JA} holds, then $Q\not\in \spn(P_1,\dots,P_L)$ if and only if there exist probability measures $G\neq F$ such that
$$\T((P_1,\dots,P_L,Q),(F,\dots,F,G))\neq\emptyset.$$
In addition, in both cases above, we may pick $F=\Leb$ and $G$ atomless.
\end{proposition}

\begin{proof}
The ``if'' is clear since $Q\in\spn(P_1,\dots,P_L)$ would imply $G\in\spn(F)=\{F\}$. For ``only if'', let $F=\Leb$ and consider first the case where \eqref{AC} holds. Then using Proposition~\ref{prop:ssww} with $d=L+1$, it suffices to prove that there exists some $G\gg F$ such that
$$\left(\frac{\d P_1}{\d Q},\dots,\frac{\d P_L}{\d Q},\frac{\d Q}{\d Q}\right)\Big|_{Q}\gcx \left(\frac{\d F}{\d G},\dots,\frac{\d F}{\d G},\frac{\d G}{\d G}\right)\Big|_{G}.$$
Equivalently, we need to show that
\begin{equation}\label{eq:exist-G}\gamma=\left(\frac{\d P_1}{\d Q},\dots,\frac{\d P_L}{\d Q}\right)\Big|_{Q}\gcx \left(\frac{\d F}{\d G},\dots,\frac{\d F}{\d G}\right)\Big|_{G}.\end{equation}

We will first consider a special type of density $\d F/\d G$ which allows us to construct $G$ such that \eqref{eq:exist-G} holds. 
Suppose that 
$$\frac{\d G}{\d F}(x)=\begin{cases}
    1&\text{ if }0\leq x\leq 1-\ee;\\
    1+\ee&\text{ if }1-\ee< x\leq 1-\frac{\ee}{2};\\
    1-\ee&\text{ if }1-\frac{\ee}{2}< x\leq 1,
\end{cases}$$
where $\ee>0$ is a  small number. Clearly, $G$ is atomless. Moreover,  $({\d F}/{\d G})|_{G}$ is concentrated on $[(1+\ee)^{-1},(1-\ee)^{-1}]$ and $\p^G[{\d F}/{\d G}=1]= 1-\ee$. Therefore, the measure $( \d F/\d G ,\dots, \d F/\d G )|_{G}$ is supported on the line segment $\{\bx\in\R^L \mid x_1=\dots=x_L\in[(1+\ee)^{-1},(1-\ee)^{-1}]\}$, with mean $\e$ and $\p^G[ \d F/\d G \neq 1]=\ee$. We will find a measure $( \d F/\d G ,\dots, \d F/\d G )|_{G}$ that satisfies the condition above and also \eqref{eq:exist-G}.

Consider a disjoint collection of closed balls $\{B_j\}_{1\leq j\leq k}$ in $\R^L$ as constructed in Lemma~\ref{lemma:balls}. By Lemma~\ref{lemma:CM}, there is $\delta>0$ and a segment  $\{\bx\in\R^L \mid x_1=\dots=x_L\in[1-\delta,1+\delta]\}$ containing $\e$, such that any  measure  of total mass $\delta$ supported on it will be smaller in extended convex order  than some $\widetilde{\gamma}$ such that $\widetilde{\gamma}\leq\gamma|_{\bigcup_{j=1}^k B_j}$.  We choose $\ee>0$ so that $(1-\ee)^{-1}<1+\delta$. As a result, the measure $G$ constructed in the above paragraph satisfies
$$ \omega:=\left(\left(\frac{\d F}{\d G},\dots,\frac{\d F}{\d G}\right)|_{G}\right)\Big|_{\R^L\setminus\{\e\}}\lcx \widetilde{\gamma}.$$
 The measure $( \d F/\d G ,\dots, \d F/\d G )|_{G}-\omega$ is concentrated at $\e$, which is smaller in convex order than any measure with barycenter $\e$ and the same total mass. Since $\bary(\gamma)=\bary(\widetilde{\gamma})=\e$, we conclude  
$$\left(\frac{\d F}{\d G},\dots,\frac{\d F}{\d G}\right)\Big|_{G}\lcx \gamma.$$

\sloppy If \eqref{AC} does not hold, then we define $Q'=Q/2+(P_1+\dots+P_L)/(2L)$, and repeat the above arguments, so that there is $\kappa$ sending $(P_1,\dots,P_L,Q')$ to some $(F,\dots,F,G')$ where $G'\neq F$. By linearity, $\kappa$ also sends $(P_1,\dots,P_L,Q)$ to $(F,\dots,F,G)$ where $G=2G'-F\neq F$.
\end{proof}

\begin{proof}[Proof of Theorem~\ref{thm:iff}]
The direction (a)$\Rightarrow$(b) is proved as Proposition~\ref{prop:calibration}, (b)$\Rightarrow$(c) is clear from definition, (c)$\Rightarrow$(e) is proved as Proposition~\ref{prop:nec}, and (e)$\Rightarrow$(d) is Proposition~\ref{prop:existence}. To show (d)$\Rightarrow$(a), let $X$ be a random variable that has a common law $F$ under $P\in\P$, and law $G$ under $Q$. Let $\phi$ be given in Lemma~\ref{lemma:p-variable}. It follows immediately that $\phi\circ X$ is an exact p-variable.  
\end{proof}

\begin{theorem}\label{thm:p-variable}
    Suppose that we are testing $\P=\{P_1,\dots,P_L\}$ against $\Q=\{Q\}$ and \eqref{JA} holds. The following are equivalent:
\begin{enumerate}[(a)]

    \item there exists a  nontrivial p-variable;
    \item there exists a bounded e-variable that has nontrivial e-power against $\Q$;
    \item  there exists an  e-variable that is nontrivial for $\Q$;
    \item it holds that $Q\not\in\conv(P_1,\dots,P_L)$.
    \end{enumerate}
\end{theorem}

\begin{remark}
    The directions (c)$\Leftrightarrow$(e) in Theorem~\ref{thm:iff} and (c)$\Leftrightarrow$(d) in Theorem~\ref{thm:p-variable} also hold without \eqref{JA}, in view of Proposition~\ref{prop:no ja}. 
\end{remark}

% \begin{remark}
%     The direction $(d)\Rightarrow(c)$ of Theorem~\ref{thm:p-variable}
% \end{remark}

\begin{example}\label{ex:1}
\begin{enumerate}[(i)]
\item Let $P_1\lawis \Ber(0.1),~P_2\lawis \Ber(0.2),$ and $Q\lawis \Ber(0.3)$. It follows that $Q\in\spn(P_1,P_2)\setminus\conv(P_1,P_2)$. By Theorems~\ref{thm:iff} and~\ref{thm:p-variable}, a nontrivial e-variable (or p-variable) exists, but an exact nontrivial e-variable (or p-variable) does not exist.
\item Let $P_1\lawis  \mathrm{N}(-1,1),~P_2\lawis  \mathrm{N}(1,1)$, and $Q\lawis  \mathrm{N}(0,1)$. By Theorem~\ref{thm:iff}, there exists a pivotal exact nontrivial e-variable (or p-variable).
\end{enumerate}

\end{example}

\begin{remark}
    {When the sample space $\X$ is finite (say $|\X|=d$) and $Q\not\in\spn\P$, it is easy to construct nontrivial exact e-variables. We can associate the distributions with their Radon--Nikodym derivatives, which are just $d$-dimensional vectors, and one can consider an e-variable of the form $1+Y$ where $Y$ is proportional to the orthogonal part of $Q$ relative to the span of $\P$ (so that $Y$ integrates to zero under any $P\in\P$, but has positive expectation under $Q$). 
    % We defer the details to Appendix {\color{red}XXX}. 
    In case $\mathcal P=\{P\}$, taking $P$ as the reference measure, this construction yields $Y=\alpha (\d Q/\d P-1)$ for any $\alpha \in [0,1]$, and the e-variable is precisely $\d Q/\d P$ when $\alpha =1$. For infinite $\X$, such a direct construction exploiting orthogonality is no longer possible because the Radon--Nikodym derivatives do not live in a Hilbert space.}
\end{remark}

Next,  Section~\ref{sec:e} constructs a powerful exact e-variable  by additionally imposing pivotality.

\section{Constructing a powerful exact e-variable}\label{sec:e}

We focus on e-variables in this section. 
Provided the existence, our next step is to maximize the e-power of an e-variable that is pivotal and exact.  
The e-power  of an e-variable $X$ can be measured by $\E^Q[\log X]$, which has long been a popular criterion; see for example~\cite{Kelly56,breiman1961optimal, bell1988game,shafer2011test, GrunwaldHK19,waudby2024estimating}.\footnote{In short, it captures the rate of growth of the test martingale under the alternative $Q$; see Section~\ref{sec:testmartingale}.}  It has been recently called the \emph{e-power} of $X$ (\cite{vovk2022efficiency}), a term we continue to use for simplicity. 
In this section, we will fix $\P=\{P_1,\dots,P_L\}$ and $\Q=\{Q\}$. Our goal is to solve
\begin{align}\label{eq:opt}
   \begin{split}
       \max\quad &\E^Q[\log X],\\ 
    \text{s.t.}\quad &X: ~X\text{ is a  pivotal exact e-variable}.
   \end{split} 
\end{align}
This optimization problem turns out to be a special case of a more general problem that is illustrated by \eqref{eq:optcx} below. Such a connection will be explained in Section~\ref{sec:econv}. We describe an equivalent condition for the existence of a maximal element for \eqref{eq:optcx} in 
Section~\ref{sec:characterize}. A further sufficient condition in the case $L=2$ is illustrated in 
Section~\ref{sec:L=2}. Section~\ref{sec:multiple observations} contains a few discussions regarding batching multiple data points and how it affects the e-power.  
Finally, we provide several examples in Section~\ref{sec:examples}. In this section, we let 
\begin{equation}\label{eq:gamma}
\gamma:= \Big(\frac{\d P_1}{\d Q},\dots,\frac{\d P_L}{\d Q}\Big)\Big|_Q.
\end{equation}
In particular, $\gamma$ is a  probability measure on $\R_+^L$ with mean $\e$.

\subsection{E-power maximization and convex order}\label{sec:econv}

We first recall the maximizer of e-power in the case of a simple null versus a simple alternative, which has an explicit form. This fact is used frequently in the above literature. 
\begin{example}\label{ex:opt}
Let us first illustrate an example with simple null $\P=\{P\}$ ($L=1$) and simple alternative $\Q=\{Q\}$.  Clearly, any e-variable is pivotal. Thus~\eqref{eq:opt} reduces to
\begin{align}
    \label{eq:opt11}\begin{split}
        \max\quad&\E^Q[\log X],\\
        \text{s.t.}\quad &X:~  X\geq 0,\, \E^P[X]=1.
    \end{split}
\end{align}
By Gibbs' inequality, the maximum value is attained by the likelihood ratio, i.e., when $X=\d Q/\d P$ (see~\cite{Shafer:2021} for this simple setting).

Below we illustrate the solution to \eqref{eq:opt11} using our theory, which sheds light on the composite null case. For simplicity, we assume \eqref{JA} and \eqref{AC}. Denote by $\gamma:=({\d P}/{\d Q}  )|_{Q}$. Consider the set $\cM_\gamma$ of probability measures $\mu$ such that $\mu\lcx \gamma$: 
$$
\cM_\gamma := \{\mu \in \Pi(\R): \mu\lcx \gamma\}
$$
Using Lemma~\ref{lemma:RN}, every $\mu\in\cM_\gamma$ (in fact also for $\mu\not\in\cM_\gamma$) corresponds to a probability measure $F$ such that $({\d F}/{\d \Leb})|_{\Leb}=\mu$.  
By Proposition~\ref{prop:ssww}, there exists a random variable $Y$ that has law $F$ under $P$ and law $\Leb$ under $Q$. Next, consider $X$  of the form $X=({\d\Leb}/{\d F})(Y)$, and we optimize $\E^Q[\log X]$ over $F$ satisfying $({\d F}/{\d \Leb})|_{\Leb}\in\cM_\gamma$. It is clear that the constraint 
$$\E^P[X]=\E^P\Big[\frac{\d\Leb}{\d F}(Y)\Big]=\E^F\Big[\frac{\d\Leb}{\d F}\Big]=1$$ 
is satisfied, and the objective in \eqref{eq:opt11} becomes
$$\E^Q[\log X]=\E^Q\Big[\log\Big(\frac{\d\Leb}{\d F}(Y)\Big)\Big]=\E^{\Leb}\Big[\log\frac{\d\Leb}{\d F}\Big]=\E^{\Leb}\Big[-\log\frac{\d F}{\d \Leb}\Big].$$
We have thus arrived at the  optimization problem
\begin{align}\label{eq:opt12}\begin{split}
    \max\quad&\E^{\Leb}\Big[-\log\frac{\d F}{\d \Leb}\Big],\\
    \text{s.t.}\quad &F\in \Pi(\R):~ F\ll\Leb,\, \frac{\d F}{\d \Leb}\Big|_{\Leb}\in\cM_\gamma.
\end{split}
    \end{align}
The value \eqref{eq:opt12} gives a lower bound on \eqref{eq:opt11}. Since the set  $\cM_\gamma$ has a maximum element $\gamma$ in convex order, the problem \eqref{eq:opt12} has a trivial solution $\E^Q[-\log({\d P}/{\d Q})]$. This corresponds to the solution to \eqref{eq:opt11} using Gibbs' inequality. 

The fact that the two values \eqref{eq:opt11} and \eqref{eq:opt12} are the same is not a coincidence and holds more generally  for composite nulls, which we will prove in Theorem~\ref{prop:Xcondition}.  With a composite null, the main difficulty arises from solving \eqref{eq:opt12}, because the set $\cM_\gamma$ has a complicated structure, and may not contain a maximum element in convex order.
\end{example}

\sloppy 
As explained in Example~\ref{ex:opt}, the first step to solving \eqref{eq:opt} is to impose the further condition that $X$ is of the form $({\d G}/{\d F})(Y)$ for some $F,G,Y$.  As a consequence of Gibbs' inequality, this does not affect the optimal value of \eqref{eq:opt}, as shown in the following result.
\begin{theorem}\label{prop:Xcondition}
   Assume \eqref{JA} and \eqref{AC}. There exists a maximizer $X$ to \eqref{eq:opt} of the form $X=({\d G}/{\d F})(Y)$, where $F,G\in\Pi(\R)$, and $Y\in\T((P_1,\dots,P_L,Q),(F,\dots,F,G))$.
\end{theorem}
The fact that the log-optimal pivotal and exact e-variable is a likelihood ratio is quite aesthetically appealing, a phenomenon that is known to be true without the restrictions of pivotality and exactness~\citep{GrunwaldHK19,larsson2024numeraire}, but in this more general case $F$ could be a \emph{sub-}probability distribution. 

% \com{This Gibb's inequality works only for log but not concave!!}

Given $X=({\d G}/{\d F})(Y)$ where $Y\in\T((P_1,\dots,P_L,Q),(F,\dots,F,G))$, we may rewrite
$$\E^Q[\log X]=\E^Q\Big[\log\Big(\frac{\d G}{\d F}(Y)\Big)\Big]=\E^{G}\Big[-\log\frac{\d F}{\d G}\Big].$$
As a consequence of Proposition~\ref{prop:ssww}, the optimization problem \eqref{eq:opt} is equivalent to finding
\begin{align}\label{eq:opt2}
  \begin{split}
      \max\quad &\E^{G}\Big[-\log\frac{\d F}{\d G}\Big],\\
      \text{s.t.}\quad&F,G\in\Pi(\R):~\left(\frac{\d F}{\d G},\dots,\frac{\d F}{\d G}\right)\Big|_{G}\lcx \left(\frac{\d P_1}{\d Q},\dots,\frac{\d P_L}{\d Q}\right)\Big|_Q.
  \end{split}  
\end{align}

More generally, since $x\mapsto -\log x$ is convex on its domain, we may formulate the problem of optimizing $\E^G[\phi(\d F/\d G)]$ for all convex function $\phi:\R_+\to\R$. In other words, let $\gamma$ be the law of $(\d P_1/\d Q,\dots,\d P_L/\d Q)$ under $Q$ and introduce the set $\cM_\gamma$ of probability measures 
 supported on $\I_L^+$ that is smaller than $\gamma$ in convex order, and our goal is  to \begin{align}\label{eq:optcx}
   \begin{split}
        \max \quad&\mu~\text{ in }~\lcx,\\
        \text{s.t.}\quad&\mu\in\cM_\gamma.
   \end{split} 
 \end{align} This will be the goal of the present section. The reader should keep in mind that unfortunately, even if \eqref{eq:optcx} allows a unique maximum element, it does not necessarily solve \eqref{eq:opt} uniquely when the logarithm in \eqref{eq:opt} is replaced by other concave functions. This is because Theorem~\ref{prop:Xcondition} requires Gibbs' inequality, where the logarithm plays a crucial role.

\subsection{Existence of the maximum element in convex order}\label{sec:characterize}

To ease our presentation, we will assume further that 
\begin{align}\label{A}
    \gamma\text{ from~\eqref{eq:gamma} does not give positive  mass to any hyperplane in }\R^L.\tag{N}
\end{align}
That is, for every half-space  $\H\subseteq\R^L,~\gamma(\partial\H)=0$. This is a technical assumption  which greatly simplifies our proofs (as we will explain in Remarks~\ref{rem:atomless} and~\ref{rem:atomless2}), and we expect that analogous results hold without such an assumption.
\begin{proposition}\label{prop:LdimHx}
   Let $\gamma$ be a  probability measure on $\R_+^L$ with mean $\e$. Consider $x\geq 0$. There exists a closed half-space $\H_x$ of $\R^L$ and a  measure $\mu_x$ supported on $\H_x$, such that
\begin{enumerate}[(i)]
\item the positive diagonal $\I^+_L\not\subseteq\H_x$;
\item $-\e\in\H_x$;
    \item $x\e\in \partial\H_x={\H_x}\cap\H_x^c$, where $\H_x^c$ is the closed complement of  $\H_x$;
    \item the measure $\mu_{\H_x^c}:=\gamma-\mu_x$ is supported on $\H_x^c$, and the barycenters of $\mu_x$ and $\mu_{\H_x^c}$ both lie on $\I^+$.
    \end{enumerate}In this case, we call $\partial\H_x$ a separating hyperplane at $x$. 
    Moreover, if \eqref{A} holds, there exists a unique  measure $\mu_x$ satisfying the above conditions, in which case it also holds that $\mu_x=\gamma|_{\H_x}$ and $\mu_{\H_x^c}=\gamma|_{\H_x^c}$.
\end{proposition}

We remark that if $\gamma$ has a strictly positive density on $\R_+^L$, then the above $\H_x$ is unique.

 Recall from \eqref{eq:optcx} that our goal is to find the maximum element in $\cM_\gamma$ in convex order.

\begin{theorem}\label{thm:maximumexist}
Assuming \eqref{A}, the following are equivalent.
\begin{enumerate}[(a)]
    \item There exists a unique maximum element $\mu$ in convex order in $\cM_\gamma$, i.e., 
         $\mu\lcx \gamma$ and  for each $\nu$ supported on $\I^+$ with $\nu\lcx \gamma$, it holds that $\nu\lcx \mu$.
    
    \item The class of measures $\{\mu_x\}_{x\geq 0}$ from Proposition~\ref{prop:LdimHx} is monotone (in the usual order), i.e., for all $x\leq y$, $\mu_x\leq \mu_y$.
\end{enumerate}
\end{theorem}

\begin{example}
    Suppose that $L=1$. It is clear from the proof of Proposition~\ref{prop:LdimHx} that condition (b) in Theorem \ref{thm:maximumexist} is always satisfied. Therefore, the maximum element $\mu$ in $\cM_\gamma$ always exists. This agrees with Example~\ref{ex:opt}, where the likelihood ratio maximizes the e-power.  
\end{example}

\begin{remark}\label{rem:atomless}
The only place we used our assumption \eqref{A} is on the uniqueness of the measure $\mu_x$ in Proposition~\ref{prop:LdimHx}.    When there is no uniqueness, the condition (b) in Theorem~\ref{thm:maximumexist} needs to be replaced by the existence of a monotone selection of measures $\{\mu_x\}_{x\geq 0}$, each of them satisfying the conditions in Proposition~\ref{prop:LdimHx}.
\end{remark}

    The condition (b) in Theorem~\ref{thm:maximumexist} is in general not easy to check, especially in higher dimensions.\footnote{In this paper when we mention ``dimension'' we typically refer to the dimension of the null, but not the dimension of the underlying space $\X$.}  Later, we supply a sufficient condition in Section~\ref{sec:L=2}, and a few examples in Section~\ref{sec:examples}.

\subsection{A sufficient condition  in case  \texorpdfstring{$|\P|=2$}{}}\label{sec:L=2}

When $L=2$, we provide a sufficient condition for the class of measures $\{\mu_x\}_{x\geq 0}$ to satisfy the monotonicity condition $x\leq y\implies \mu_x\leq\mu_y$. In view of Theorem~\ref{thm:maximumexist}, this condition implies the existence of the maximum element $\mu$. We keep the same setting as in Section~\ref{sec:characterize} and assume \eqref{A}, with the exception that $L=2$.

\begin{theorem}\label{theorem:cx}  \sloppy  Assume \eqref{A}, \eqref{JA},  \eqref{AC}. Suppose that there exists a  convex set $\Gamma\subseteq \R^{2}$ such that $\gamma(\partial \Gamma)=1$.\footnote{This assumption is far from being necessary, but might be convenient to verify.} Then there exists a unique maximum element $\mu$ in convex order in $\cM_\gamma$. Moreover, $\mu$ is the unique probability measure on the $\I_2^+$ with $\mu([0,x]^2)=\mu_x(\R^2)$, where $\mu_x$ was given in Proposition~\ref{prop:LdimHx} applied with $L=2$. In particular, there exist distinct measures $F,G\in\Pi(\R)$ such that  $( \d F/\d G , \d F/\d G )|_G=\mu$, attaining the maximum in \eqref{eq:opt2}.
\end{theorem}

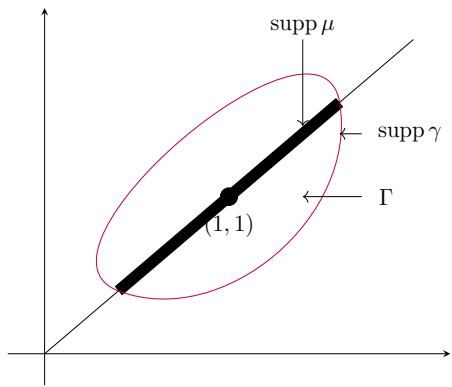
\begin{figure}[!t]\begin{center}\resizebox{0.4\textwidth}{0.33\textwidth}{
\begin{tikzpicture}

\begin{axis}[axis lines=middle,
            enlargelimits,
            ytick=\empty,
            xtick=\empty,]
\addplot[name path=F,thin,black,domain={0:1}] {2*x} node[pos=.9, above]{};
\addplot[name path=G,line width=5pt,black,domain={0.2:0.8}] {2*x} node[pos=.9, above]{};
\node[circle,draw,fill=black, minimum size=0.1pt,scale=0.8,label=below:{$(1,1)$}] (D) at  (500,100){};
\path[-,purple,every node/.style={font=\sffamily\small}] (800,160) edge[bend right=120] node [left] {} (200,40);
\path[-,purple,every node/.style={font=\sffamily\small}] (800,160) edge[bend left=60] node [left] {} (200,40);
\path[->,black](860,140) edge[bend left=0] node [right] {\hspace{0.3cm}$\supp\gamma$} (800,140);
\path[->,black](860,100) edge[bend left=0] node [right] {\hspace{0.6cm}$\Gamma$} (700,100);
\path[->,black](700,200) edge[bend left=0] node [above=0.6cm] {$\supp\mu$} (700,145);
\end{axis}
\end{tikzpicture}
}\end{center}

\caption{Illustration of Theorem~\ref{theorem:cx}. The convex set $\Gamma$ is enclosed by the red contour $\partial\Gamma$ on which $\gamma$ (the law of $({\d P_1}/{\d Q},{\d P_2}/{\d Q})$ under $Q$) is supported. The measure $\mu$ is supported on the thick segment on the diagonal $\I$. }
\label{fig:1}
\end{figure}

\begin{remark}\label{rem:atomless2}
    With essentially the same arguments, we may remove assumption \eqref{A} from Theorem~\ref{theorem:cx}. With the presence of atoms, selecting any monotone collection $\{\mu_x\}_{x\geq 0}$ would be enough; see Remark~\ref{rem:atomless}.
\end{remark}

\subsection{On multiple observations}\label{sec:multiple observations}

Before we proceed, let us discuss the case with multiple data points. Suppose that instead of one data point, we observe $n$ iid data points $Z_1,\dots,Z_n$ in the space $\X$ from the experiment. The e-variable is built based on the $n$ data points together instead of a single data point. In other words, given $\P=\{P_1,\dots,P_L\}$ and $\Q=\{Q\}$, we build an e-variable for $\P^n:=\{P_1^n,\dots,P_L^n\}$ that is pivotal, exact, and has nontrivial e-power against $\Q^n:=\{Q^n\}$. We first see that, as long as $Q\not\in \mathcal P$ and $\P$ is linearly independent, at most two  observations are needed to build a pivotal and exact e-variable based on
 Theorem~\ref{thm:iff}. Furthermore, without linear independence of $\P$, a finite number of observations would suffice when the underlying space $\X$ is Euclidean.  %Note that no joint non-atomicity   needs to be assumed for this result.  

 %$$ \ell_n=\max\{\E^{Q^n}[\log X]\mid X\text{ is an  pivotal exact e-variable that has nontrivial power against $\Q^n$ for }\P^n\}.$$

% \begin{proposition}\label{lemma:spansquared}
%     Suppose that $P_1,\dots,P_L$ are linearly independent probability measures on $\X$, and $Q\in\Pi(\X)$ satisfies  \eqref{AC} and $Q\not\in\P=\{P_1,\dots,P_L\}$. Then either $Q\not\in\spn\P$ or $Q^2\not\in\spn\P^2$ (or both). In particular, either $Q\not\in\conv\P$ or $Q^2\not\in\conv\P^2$ (or both).
% \end{proposition}

\begin{theorem}\label{thm:spansquared}

Suppose that $\X$ is an Euclidean space and $P_1,\dots,P_L$ are  distinct probability measures on $\X$. If $Q\in\Pi(\X)$ satisfies $Q\not\in\P=\{P_1,\dots,P_L\}$, then there exists $k\geq 1$ such that $Q^k\not\in\spn\P^k$ (and in particular $Q^k\not\in\conv\P^k$). Moreover, if we also assume that $Q$ satisfies \eqref{AC} and that $P_1,\dots,P_L$ are linearly independent, then either $Q\not\in\spn\P$ or $Q^2\not\in\spn\P^2$ (or both); in particular, either $Q\not\in\conv\P$ or $Q^2\not\in\conv\P^2$ (or both). 
\end{theorem}

In the last claim above, one can show that neither the linear independence condition nor \eqref{AC} can be removed.   
The proof of Theorem~\ref{thm:spansquared} is  put in Appendix~\ref{app:sec4}, which relies on the following fundamental fact: If $\X$ is an Euclidean space and $P_1,\dots,P_L$ are  distinct probability measures on $\X$, then there exists $k\geq 1$ (possibly large) such that  $P_1^k,\dots,P_L^k$ are linearly independent. This  fact may be  known, but we are not aware of a proof in the literature, and we present it as Lemma~\ref{lem:remark:power} in Appendix~\ref{app:sec4}. The weaker statement that there exists $k$ for which $Q^k\not\in\conv\P^k$ also follows from Lemma 2 of \cite{berger1951uniformly}.

% \begin{remark}\label{remark:power}
%   When $P_1,\dots,P_L$ are not linearly independent, we may take further powers of them to achieve linear independence. More precisely, if $\X$ is an Euclidean space and $P_1,\dots,P_L$ are  distinct probability measures on $\X$, there exists $k\geq 1$ (possibly large) such that  $P_1^k,\dots,P_L^k$ are linearly independent. This fundamental fact may be  known, but we are not aware of a proof in the literature, and we present it as Lemma~\ref{lem:remark:power} in Appendix~\ref{app:sec4}. Combined with Proposition~\ref{lemma:spansquared}, we see that without linear independence, we have the weaker conclusion that for some $k\geq 1$, $Q^k\not\in\spn\P^k$ (and in particular $Q^k\not\in\conv\P^k$), i.e., batching $k$ data points suffices.\footnote{The existence of $k\geq 1$ such that $Q^k\not\in\conv\P^k$ also follows from Lemma 2 of \cite{berger1951uniformly}.} In fact, the absolute continuity assumption \eqref{AC} may also be removed in this case, using a minor modification of the proof presented in Appendix~\ref{app:sec4}.
% \end{remark}

\begin{example}
Suppose that $P_1\lawis \Ber(0.1),~P_2\lawis \Ber(0.2),$ and $Q\lawis \Ber(0.3)$. We explained in Example~\ref{ex:1} that an exact nontrivial e-variable does not exist. Nevertheless, Theorem~\ref{thm:spansquared} implies that $Q^2\not\in \spn(P_1^2,P_2^2)$, and hence an exact and nontrivial e-variable exists for a batch of two data points. An example of such an e-variable $X$ is given by 
\begin{align*}
    X(\omega)\approx\begin{cases}
        1.009&\text{ for }\omega=(0,0);\\
        0.939&\text{ for }\omega=(0,1),(1,0);\\
        1.338&\text{ for }\omega=(1,1).
    \end{cases}
\end{align*}
\end{example}

Let us denote by $ \ell_n$ the maximum e-power with $n$ data points for $\P^n$ against $\Q^n$ using a pivotal and exact e-variable, similarly as in \eqref{eq:opt}.

\begin{proposition}\label{prop:multidata}
    In the setting above, suppose that  \eqref{AC} and \eqref{JA} hold, and $Q\not\in\P=\{P_1,\dots,P_L\}$. For any $n,m\in\N$, $ \ell_{n+m}\geq  \ell_n+ \ell_m$. In particular, $ \ell_n/n$ converges to a positive limit less than or equal to $\min_{P\in \mathcal P}\E^Q[\log (\d Q/\d P)]$.
\end{proposition}

The above proposition formalizes the straightforward observation that constructing an e-value using $m+n$ points is potentially more powerful than multiplying two e-values together that were constructed separately using $m$ and $n$ points respectively. 

\smallskip

\textbf{Is there a loss of e-power  caused by imposing exactness or pivotality?} 
% \red{rewrite--the ``loss" does not make sense} By Proposition~\ref{prop:multidata}, $\ell_{2^n}/2^n$ is increasing in $n$. 
% Let $\ell^*_n $ be the maximum e-power for $\P^n$ against $\Q^n $  achieved without assuming pivotality and exactness. Clearly, $\ell^*_n \le n \min_{P\in \mathcal P}\E^{Q}[\log (\d Q/\d P)]$. 
% Note that the value
% $$ \ell^*_n/n - \ell_n/n 
% $$
% is the average loss of e-power caused by imposing pivotality and exactness, which  is getting weaker as the number of data points grows.
 {The superadditivity property established in Proposition~\ref{prop:multidata} implies in particular that $\ell_{2^n}/2^n$ is increasing in $n$. {The intuitive reason of the increase in the average e-power is partly due to the fact that the pivotality constraint becomes less restrictive for a higher number of observations.} To see this, imagine laws $P\in\P$ with a complicated entangled overlapping structure. To achieve pivotality, we need to send all laws $P$ simultaneously to a single distribution $F$, the ways of which may be quite limited due to the overlapping structure.\footnote{For instance, if $P\in\P$ all have disjoint support (no overlap), the transport map can be picked independently on the disjoint supports to send $P$ to $F$, but this is not possible of $P\in\P$ all have the same support.} On the other hand, with multiple observations, the laws $P^n,\,P\in\P$ have much fewer overlapping parts than $P\in\P$ do (for instance, $P^\infty,\,P\in\P$ are mutually singular), meaning that there are more ways to achieve pivotality. }

%  {To see the above discussion intuitively, let us consider the case where the laws $\P=\{P_1,\dots,P_L\}$ have completely disjoint supports. In this case, picking an (exact) e-value is equivalent to simultaneously picking $L$ (exact) e-values for testing $P_j,\,1\leq j\leq L$ against $Q$. }

It remains an open question whether $\ell_n/n \to \min_{P\in\mathcal P}\E^Q[\log (\d Q/\d P)]$. Note that this conjectural limit can be  different from $\min_{P\in\conv\mathcal P}\E^Q[\log (\d Q/\d P)]$ because the linearity structure is lost after taking powers. 
If the above question is answered in the affirmative, then the loss of e-power  
vanishes asymptotically, by noting that $n \min_{P\in\mathcal P}\E^Q[\log (\d Q/\d P)]$ is an upper bound on the theoretical best e-power for testing $\mathcal P^n$ against $\mathcal Q^n$ (see Example~\ref{ex:opt}).
In Example~\ref{ex:nstep}, we present a   setting of Gaussian distributions in which $\ell_n/n \to \min_{P\in\mathcal P}\E^Q[\log (\d Q/\d P)]$ holds true. 
We conjecture that this limit holds true in general, but we did not find a proof.\footnote{ {The argument in Example~\ref{ex:nstep} is analytical. On the other hand, numerical verification of this conjecture remains a challenging task due to drastic extremal values of the Radon--Nikodym derivatives (in high dimensions, almost all mass of $(\d P^n/\d Q^n)|_{Q^n}$ concentrates near 0 or $\infty$), exponential time complexity, and the slow convergence of $\ell_n/n$. We leave it as an open problem to design a more efficient iterative algorithm for this problem (or more generally, computing numerically the best e-power in high dimensions), or to prove that one cannot exist.} }

\subsection{Examples}\label{sec:examples}

The condition in Theorem~\ref{theorem:cx} that $\gamma=({\d P_1}/{\d Q},{\d P_2}/{\d Q})|_{Q}$ is supported on the boundary of a convex set is not very restrictive.  When $P_1,P_2,Q\in\Pi(\R)$, the vector of density functions $(({\d P_1}/{\d Q})(x),({\d P_2}/{\d Q})(x))$ forms a parameterized curve in $\R^2$ by $x\in\R$. In certain nice cases, such a curve lies on the boundary of a convex set. We illustrate with a few examples below.

\begin{example}\label{ex:d=2}
    Consider $P_1\lawis  \mathrm{N}(-1,1),~P_2\lawis  \mathrm{N}(1,1)$, and $Q\lawis  \mathrm{N}(0,1)$. It follows from a direct computation that
    $$\gamma=\left(\frac{\d P_1}{\d Q},\frac{\d P_2}{\d Q}\right)\Big|_{Q}=\left(e^{-\xi-1/2},e^{\xi-1/2}\right)\Big|_{\xi\lawis \mathrm{N}(0,1)},$$
which is supported on the hyperbola $\{(x_1,x_2)\in \R_+^2\mid x_1x_2=1/e\}$, the boundary of the  convex set $\{(x_1,x_2)\in \R_+^2 \mid x_1x_2\geq 1/e\}$. 
By Theorem~\ref{theorem:cx}, there exists a unique maximal element $\mu$ in $\cM_\gamma$ in convex order.  

Using the notation from Proposition \ref{prop:LdimHx}, it is easy to see that $\H_x=\{(x_1,x_2)\in \R^2 \mid x_1+x_2\leq 2x\}$ and $\mu_x=\gamma|_{\H_x}$. Moreover, Theorem~\ref{theorem:cx} yields that $\mu$ is the unique probability measure on $\I^+$ with 
\begin{align}
    \mu([0,x]^2)=2\Phi\left(\log(\sqrt{e} x+\sqrt{ex^2-1})\right)-1\text{ for }x\geq \frac{1}{\sqrt{e}},\label{eq:density}
\end{align} where $\Phi$ is the Gaussian cumulative density function. 
It can be directly seen from the figure below that points $x,y\in\R$ are  shrunk to a single point precisely when the points $(e^{x-1/2},e^{-x-1/2})$ and $(e^{y-1/2},e^{-y-1/2})$
 are symmetric around $\I$. This happens if and only if $x=-y$. In other words, 
 the most powerful pivotal e-variable is a function of $|Z|$, where $Z$ is the observed data point.
 Using Example~\ref{ex:opt} on testing the simple hypothesis $|Z|\lawis |\xi+1|$ against  $|Z|\lawis |\xi|$,  this e-variable is given by 
 $ X=2 e^{1/2} / (e^Z+ e^{-Z})  =  e^{1/2} \cosh(Z)^{-1} $,
 and the e-power is
 $\E^Q[\log X] \approx 0.125$.
 In the sequential setting where iid observations $Z_1,\dots,Z_n$ are available (treated in the next example), we effectively reduce the filtration generated by $Z_1,\dots,Z_n$ to the one generated by $|Z_1|,\dots,|Z_n|$. This corresponds to the intuition that taking absolute value transports $P_1,P_2$ to the same measure but not for $Q$, and indeed this is the optimal solution to \eqref{eq:opt2}.
\end{example}

\begin{example}\label{ex:nstep}
We consider the setting in Example~\ref{ex:d=2} but instead of one data point, we observe $n$ iid data points 
% \com{I think we should stick to $Z_i$ for data when they are needed} 
$Z_1,\dots,Z_n$ in the experiment. Here, we build an e-variable based on the $n$ data points together  instead of 
building an e-variable for each data point; this allows for more flexibility than Example~\ref{ex:d=2}.  
  In this setting, $P_1=\mathrm{N}(-\mathbf 1_n,I_n)$,
 $P_2=\mathrm{N}(\mathbf 1_n,I_n)$,
 and  $Q=\mathrm{N}(\mathbf 0_n,I_n)$,
 where $\mathbf 1_n=(1,\dots,1)\in \R^n$,  $\mathbf 0_n=(0,\dots,0)\in \R^n$,
 and $I_n$ is the $n\times n$ identity matrix.
 It follows from a direct computation that 
$$\gamma=\left(\frac{\d P_1}{\d Q},\frac{\d P_2}{\d Q}\right)\Big|_{Q}=\left(e^{-\xi-n/2},e^{\xi-n/2}\right)\Big|_{\xi\lawis {\mathrm{N}}(0,n)},$$
    which is very similar to Example~\ref{ex:d=2}. 
%By Theorem~\ref{theorem:cx}, there exists a unique maximal element $\mu$ in $\cM_\gamma$ in convex order.  
Using a similar argument as in Example~\ref{ex:d=2},  the most powerful pivotal e-variable is given by 
$E_n=  e^{n/2}   \cosh({\sum_{i=1}^n Z_i})^{-1} $.
Note that this is different from the sequential one built in Example~\ref{ex:d=2}
which is 
$E_n^*= e^{n/2} \prod_{i=1}^n \cosh(Z_i)^{-1}$.
The contrast between $E_n$ and $E_n^*$ is interesting to discuss. 
On the one hand, $E_n$ has better e-power than $E_n^*$ since $E_n\ge E_n^*$ due to the log-convexity of the  $\cosh$ function. This is intuitive, as $E^*_n$ effectively tests more null hypotheses such as
$\mathrm{N}(\boldsymbol \mu,I_n)$ for $\boldsymbol \mu \in \{-1,1\}^n$ than $E_n$.
On the other hand, $n\mapsto  E^*_n$ is a martingale under both $P_1$ and $P_2$, 
but  we can check that $n\mapsto E_n$ is not a martingale under either $P_1$ or $P_2$. 
In Section~\ref{sec:alg}, we will compare the e-power of the two  approaches numerically, and in Section~\ref{sec:testmartingale}, we further discuss test martingales.
Finally, we note that $\ell_n/n=\E^{Q} [\log E_n]/n \to 1/2 = \min_{i=1,2} \E^{Q}[\log (\d Q/\d P_i)]/n $,
and hence the upper bound in Proposition~\ref{prop:multidata} is sharp. On the other hand, $\E^Q[\log E_n^*]/n =1/2-\E^Q[\log \cosh (Z_1)]\approx 0.125.$
\end{example}

\begin{figure}[!t]\begin{center}\resizebox{0.4\textwidth}{0.33\textwidth}{
\begin{tikzpicture}

\begin{axis}[axis lines=middle,
    enlargelimits,xmin=0,xmax=2,ymin=0,ymax=2,
            ylabel=$x_2$,
            xlabel=$x_1$,samples=300]
\addplot[name path=F,thin,black,domain={0:2}] {x} node[pos=.9, above]{};
\addplot[name path=G,line width=3pt,purple,domain={0.6:2}] {x} node[pos=.9, above]{};
\addplot[name path=G,thin,black,domain={0.18:2}] {1/(e*x)} node[pos=.9, above]{};
\node[circle,draw,fill=black, minimum size=0.1pt,scale=0.4,label=below:{$\e$}] (DD) at  (100,100){};
\node[circle,draw,fill=black, minimum size=0.1pt,scale=0.4,label=below:{}] (DD) at  (80,80){};
\node[circle,draw,fill=black, minimum size=0.1pt,scale=0.4,label=below:{}] (DD) at  (110,110){};

\path[dashed,<-,black](84,76) edge[bend left=0] node [left,yshift=-0.2cm,xshift=-0.11cm] {} (134,26);
\path[dashed,<-,black](76,84) edge[bend left=0] node [left,yshift=-0.2cm,xshift=-0.11cm] {} (26,134);
\path[dashed,<-,black](114,106) edge[bend left=0] node [left,yshift=-0.2cm,xshift=-0.11cm] {} (200,20);
\path[dashed,<-,black](106,114) edge[bend left=0] node [left,yshift=-0.2cm,xshift=-0.11cm] {} (20,200);
\end{axis}
\end{tikzpicture}
}\end{center}

\caption{An illustration of Example~\ref{ex:d=2}: $\gamma$ is supported on the hyperbola $x_1x_2=e^{-1}$, the optimal $\mu$ is supported on the red ray. Dashed arrows indicate the reduction of filtration.}
\label{fig:4}
\end{figure}
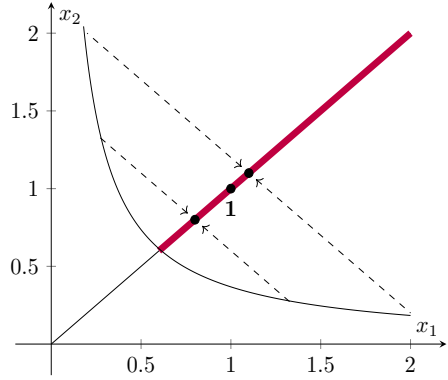

\begin{example}
Let us examine some further  sufficient conditions with $L=2$. Consider  $P_1,P_2,Q\in\Pi(\R)$ such that $P_1,P_2\ll Q$ and $\d P_i/\d Q\in C^2(\R)$ for $i=1,2$. Recall that a simple $C^2$  parameterized curve $(x(t),y(t))$ in $\R^2$ lies on the boundary of a convex set if and only if its curvature 
$$
k=\frac{x'y''-y'x''}{((x')^2+(y')^2)^{3/2}}
$$
is always nonnegative or always nonpositive (Theorem 2.31 of  \cite{kuhnel2015differential}).   Therefore, $(\d P_1/\d Q,\d P_2/\d Q)$ lies on the boundary of a convex set if  
$$\Big(\frac{\d P_1}{\d Q}\Big)'\Big(\frac{\d P_2}{\d Q}\Big)''-\Big(\frac{\d P_1}{\d Q}\Big)''\Big(\frac{\d P_2}{\d Q}\Big)'$$remains of a constant sign. As a simple example, this is the case if $P_1,P_2,Q$ are Gaussian distributions on $\R$ with different means but the same variance, or with the same mean but different variances. In particular, this recovers Example~\ref{ex:d=2}.

More generally, suppose that $P_1,P_2,Q$ have densities $p_1,p_2,q\in C^2(\R)$ where $q$ is strictly positive, and denote by $W(f_1,\dots,f_n)$ the  Wronskian of $f_1,\dots,f_n$. Then we have the further  sufficient condition that  
$$W\left(({p_1}/{q})',({p_2}/{q})'\right)\neq 0\text{ everywhere},$$or equivalently, $W(p_1,p_2,q)(x)\neq 0$ for all $x\in\R$.  By the Abel-Liouville identity (\cite{teschl2012ordinary}), this is the case if $p_1,p_2,q$ form a fundamental system of solutions of the ODE 
$$y^{(3)}=a_2(x)y^{(2)}+a_1(x)y^{(1)}+a_0(x)y$$for some continuous functions $a_i:\R\to\R,~0\leq i\leq 2$. 
\end{example}

In higher dimensions with more than two nulls, Theorem~\ref{theorem:cx} is often not applicable. Nevertheless, given enough symmetry, we may directly compute $\{\mu_x\}$ from Theorem~\ref{thm:maximumexist} and prove that they are monotone. Surprisingly, many intuitively straightforward tests are suboptimal. 

\begin{example} 
\sloppy Consider probability measures $P_1\lawis  \mathrm{N}((0,1),I),~P_2\lawis  \mathrm{N}((-\sqrt{3}/2,-1/2),I),$ $P_3\lawis  \mathrm{N}((\sqrt{3}/2,-1/2),I)$, and $Q\lawis  \mathrm{N}((0,0),I)$. Note that Theorem~\ref{theorem:cx} is not directly applicable here. It is natural to guess from Example~\ref{ex:d=2} that the optimal solution is the Euclidean norm, i.e., the distance from $0$ in $\R^2$. On the contrary, we show this is not the case.
A routine computation gives that
\begin{align*}
    \gamma&=\left(\frac{\d P_1}{\d Q},\frac{\d P_2}{\d Q},\frac{\d P_3}{\d Q}\right)\Big|_{Q}\\
    &=\Big(e^{\xi_1-1/2},e^{({-\xi_1-\sqrt{3}\xi_2-1})/{2}},e^{({-\xi_1+\sqrt{3}\xi_2-1})/{2}}\Big)\Big|_{\xi_1,\xi_2\lawis \mathrm{N}((0,0),I)\text{ independent}}.
\end{align*}
Note that this forms an exchangeable random vector and that the support of $\gamma$ is contained in $\{(x,y,z)\in\R^3_+ \mid xyz=e^{-3/2}\}$. By symmetry and exchangeability, the unique optimal solution shrinks the set $R_a:=\{(x,y,z) \mid x+y+z=a,xyz=e^{-3/2}\}$ into a single point in $\R^3$ for every $a>0$. For $(\xi_1,\xi_2)\in\R^2$, the Radon--Nikodym derivative $(\d P_1/\d Q,\d P_2/\d Q,\d P_3/\d Q)(\xi_1,\xi_2)\in R_a$ if and only if  \begin{align}
    h(\xi_1,\xi_2):=e^{{\xi_1}}+e^{({-\xi_1-\sqrt{3}\xi_2})/{2}}+e^{({-\xi_1+\sqrt{3}\xi_2})/{2}}=\sqrt{e}\,a.\label{h}
\end{align}
 In other words, we reduce the filtration generated by the sequence of observations $Z_1,\dots,Z_n$ to the one generated by $h(Z_1),\dots,h(Z_n)$. 
It is clear that \eqref{h} does not agree  with $\xi_1^2+\xi_2^2=a'$ for any  $a'\in\R$, so  taking the Euclidean  distance from $0$ instead of $h$ is suboptimal.  
Using a general technique of constructing the most powerful e-variable in Section~\ref{sec:shine2} below, one can show that the e-variable takes the form
$X=3(e^{Z^{(1)}-1/2}+e^{(-Z^{(1)}-\sqrt{3}Z^{(2)}-1)/2}+e^{(-Z^{(1)}+\sqrt{3}Z^{(2)}-1)/2})^{-1}=3\sqrt{e}h(Z)^{-1}$, 
where $Z=(Z^{(1)},Z^{(2)})$ forms a single  observation. 
This example also generalizes to more than three nulls (Gaussian with the same variance) whose means form a regular polygon centered at $0$. 

\end{example}

Finally, we supply the following two examples illustrating an explicit calculation of $\{\mu_x\}_{x\geq 0}$ with the presence of atoms in $\gamma$. The following example also shows that without pivotality,  the maximum value of \eqref{eq:opt}  increases.

\begin{example}
Let $a,b,c,d\in(0,1)$ such that $\max(a+c,b+d)<1$. On the probability space $\Omega=[0,3]$, we define  measures $P_1,P_2,Q$ where
\begin{enumerate}[(i)]
    \item $P_1$ has density $a\bone_{[0,1]}+c\bone_{[1,2]}+(1-a-c)\bone_{[2,3]}$;
    \item $P_2$ has density $b\bone_{[0,1]}+d\bone_{[1,2]}+(1-b-d)\bone_{[2,3]}$;
    \item $Q$ is uniformly distributed.
\end{enumerate}
It is clear that $P_1,P_2\ll Q$ and  $(P_1,P_2,Q)$ is jointly atomless. The measure $\gamma=({\d P_1}/{\d Q},{\d P_2}/{\d Q})|_{Q}$ has the form
$$\gamma=\frac{1}{3}\left(\delta_{(3a,3b)}+\delta_{(3c,3d)}+\delta_{(3(1-a-c),3(1-b-d))}\right).$$
Note that assumption \eqref{A} is not satisfied. In the following, we specify the choices of $a=0.2,~b=0.3,~c=0.5,~d=0.6$. The triangle connecting the points $(0.6,0.9),(1.5,1.8),(0.9,0.3)$ intersects with $\I^+$ at the points $(0.7,0.7)$ and $(1.3,1.3)$. Note that the measures
$$\frac{1}{3}\delta_{(3a,3b)}+\frac{1}{6}\delta_{(3(1-a-c),3(1-b-d))}\quad\text{ and }\quad\frac{1}{3}\delta_{(3c,3d)}+\frac{1}{6}\delta_{(3(1-a-c),3(1-b-d))}$$
have barycenters equal to $(0.7,0.7)$ and $(1.3,1.3)$ respectively. In particular, we may pick 
\begin{align*}
    \mu_x=\begin{cases}
        0&\text{ for }x<0.7;\\
        \frac{1}{3}\delta_{(3a,3b)}+\frac{1}{6}\delta_{(3(1-a-c),3(1-b-d))}&\text{ for }0.7\leq x<1.3;\\
        \gamma&\text{ otherwise}.
    \end{cases}
\end{align*}so that $\{\mu_x\}_{x\geq 0}$ is monotone and  satisfies the four conditions in Proposition~\ref{prop:LdimHx}. 

In view of Theorem~\ref{theorem:cx} and Remark~\ref{rem:atomless2}, the maximum is attained in \eqref{eq:opt2} and hence in \eqref{eq:opt}, by the choice $X(\omega)=(10/7)\bone_{\{\omega\in [0,1]\cup[2,2.5]\}}+(10/13)\bone_{\{\omega\in[1,2]\cup[2.5,3]\}}$.\footnote{The way that $X$ is obtained will be explained in detail in Section~\ref{sec:shine2} below.} The optimal value is $\approx 0.047$. 
On the other hand, if we remove the constraint that $X$ is pivotal, then with $X\approx 1.63\delta_{[0,1]}+0.66\delta_{[1,2]}+1.15\delta_{[2,3]}$, we have $\E^Q[\log X]\approx 0.07$, showing that the maximum in \eqref{eq:opt} increases.

%a standard convex optimization yields 1.63006789 1.14540028 0.66073267
%$$\sup\left\{\E^Q[\log X]:~X\geq 0,~\E^{P_1}[X]=\E^{P_2}[X]=1\right\}\approx 0.07.$$

\end{example}

 {In the case where $\gamma$ is not supported on the boundary of a convex set in $\R^2$, the following example shows that the conclusion of Theorem \ref{theorem:cx} may not hold.}

\begin{example}
    Let $\gamma$ be a probability measure on $\R^2$ centered at $\e$ that is supported on
    $$\{(x,y)\mid x+y\leq 1\}\cup \{(x,y)\mid x\leq y,~x\leq x_0\}\cup[1,\infty)^2,$$
    where $x_0\approx 0.903$ is the unique real solution to $8x^3=6x^2+1$.  Assume that $\d\gamma/\d(\Leb\otimes\Leb)=1$ on the set $\{(x,y)\mid x+y\leq 1\}\cup \{(x,y)\mid x\leq y,~x\leq x_0\}$. A routine computation shows that the separating hyperplanes at $x=1/2,x_0$ are $\H_{1/2}=\{(x,y)\mid y=1-x\}$ and $\H_{x_0}=\{(x,y)\mid x=x_0\}$. In particular, $\mu_{1/2}\leq \mu_{x_0}$ does not hold. Theorem \ref{thm:maximumexist} then implies that there is no maximum element $\mu$ in convex order in $\M_\gamma$.  
\end{example}

\section{The SHINE construction}\label{sec:alg}

The current section develops the SHINE construction (Separating Hyperplanes Iteration for Nontrivial and Exact e/p-variables), that effectively produces a pivotal nontrivial exact e/p-variable via separating hyperplanes (see Proposition~\ref{prop:LdimHx}, which is the key to our construction). Unless otherwise stated, we follow the setup of Section~\ref{sec:e}. 

The first goal of the SHINE construction is to solve the optimization problem \eqref{eq:optcx}. In the case where the condition in Theorem~\ref{thm:maximumexist} is satisfied, the construction outputs the maximum element. When the maximum element $\mu$ does not exist or when the condition (b) in Theorem~\ref{thm:maximumexist} is hard to check, we provide  a reasonable \emph{maximal} element $\mu$ in convex order. 
In the second part of the SHINE construction, we recover the corresponding e/p-variable from the output $\mu$ in the first part. The two parts are respectively illustrated in Sections~\ref{sec:shine1} and~\ref{sec:shine2}. We end this section by providing examples and simulation results in Section~\ref{sec:simulation}.

\subsection{Description of the SHINE construction} \label{sec:shine1}

\sloppy Start with $\mu^{(0)}=\delta_{\e}$, $x^{(0)}_1=\e$, and $\mu^{(0)}_1=\gamma$ from~\eqref{eq:gamma}. At step  $s\geq 0$, we are given $\mu^{(s)},~\{x^{(s)}_k\}_{1\leq k\leq 2^s},$ and $\{\mu^{(s)}_k\}_{1\leq k\leq 2^s}$. For each $k$, we apply Proposition~\ref{prop:LdimHx} to the sub-probability measure $\mu^{(s)}_k$ at the point $x^{(s)}_k$. This yields a unique decomposition of $\mu^{(s)}_k$ into two measures, each having a barycenter on $\I^+$. Denote them by $\mu^{(s+1)}_{2k-1}$ and $\mu^{(s+1)}_{2k}$. For $1\leq k\leq 2^{s+1}$, define $x^{(s+1)}_k=\bary(\mu^{(s+1)}_{k})$. Finally, let $\mu^{(s+1)}$ be the probability measure having mass $\mu^{(s+1)}_{k}(\R^L)$ on $x^{(s+1)}_k$ for every $k$, i.e.,
\begin{align}
    \mu^{(s+1)}:=\sum_{k=1}^{2^{s+1}}\mu^{(s+1)}_{k}(\R^L)\delta_{x^{(s+1)}_k}.\label{eq:mun}
\end{align}
 The output of the SHINE construction at step $s$ is 
the measure $\mu^{(s)}$.

It is easy to see that each $\mu^{(s)}$ is centered at $\e$ and supported on $\I^+$. Moreover,  $\mu^{(s)}\lcx\gamma$ by Strassen's theorem because $\mu^{(s)}$ is the aggregation of barycenters of different components in the decomposition of $\gamma$. By Markov's inequality, the sequence $\{\mu^{(s)}\}$ is tight and allows a weak limit. In fact, an even stronger assertion can be made. Define $\{X_s\}$ as the coupling of the first coordinate of $\{\mu^{(s)}\}$ such that $X_0=1$ and at each $s\geq 0$, for $j=2k-1,2k$, 
\begin{equation}
\p\left[X_{s+1}=x^{(s+1)}_{j} \mid X_s=x^{(s)}_k\right]=\frac{\mu^{(s+1)}_{j}(\R^L)}{\mu^{(s+1)}_{2k-1}(\R^L)+\mu^{(s+1)}_{2k}(\R^L)}.\label{eq:Xn}
\end{equation}
By construction,
$$x^{(s+1)}_{2k-1}\mu^{(s+1)}_{2k-1}(\R^L)+x^{(s+1)}_{2k}\mu^{(s+1)}_{2k}(\R^L)=x^{(s)}_k\mu^{(s)}_k(\R^L)=x^{(s)}_k(\mu^{(s+1)}_{2k-1}(\R^L)+\mu^{(s+1)}_{2k}(\R^L)).$$
It can thus be checked by direct calculation that $\E[X_{s+1} \mid  X_s = x^{(s)}_k] = x^{(s)}_k$, meaning that $\{X_s\}$ forms a nonnegative martingale, and hence converges a.s.~to some $X_\infty$ by the martingale convergence theorem.  {We call $\{X_s\}$ the SHINE martingale (associated with $\gamma$).} Denote by $\mu$ the law of $X_\infty \e=(X_\infty,\dots,X_\infty)$. Then  $\mu\lcx\gamma$ by Lemma~\ref{lemma:cxproperty}(ii).

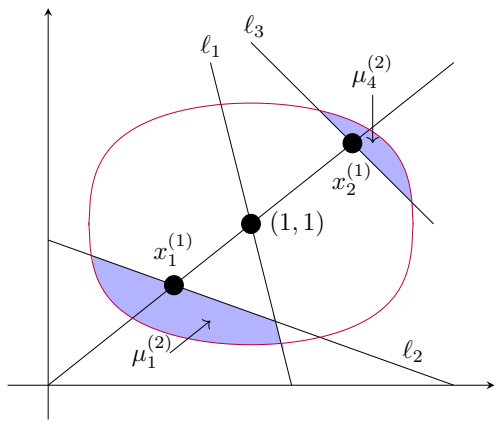
\begin{figure}[!t]\begin{center}\resizebox{0.44\textwidth}{0.36\textwidth}{
\begin{tikzpicture}

\begin{axis}[axis lines=middle,
            enlargelimits,
            ytick=\empty,
            xtick=\empty,samples=300]
\addplot[name path=F,thin,black,domain={0:1}] {2*x} node[pos=.9, above]{};
\addplot[name path=E,thin,black,domain={0.5:0.95}] {1.5-2.5*(x-0.75)} node[pos=.02, above]{$\ell_3$};
\addplot[name path=B,thin,black,domain={0:1}] {0.9-0.9*x} node[pos=.9, above]{$\ell_2$};
\addplot[name path=D,thin,black,domain={0.4:0.6}] {1-10*(x-0.5)} node[pos=.01, above]{$\ell_1$};
%\addplot[name path=G,line width=5pt,black,domain={0.2:0.8}] {2*x} node[pos=.9, above]{};
\node[circle,draw,fill=black, minimum size=0.1pt,scale=0.8,label=right:{$(1,1)$}] (DD) at  (500,100){};
\node[circle,draw,fill=black, minimum size=0.1pt,scale=0.8,label=above:{$x^{(1)}_1$}] (DD) at  (310,62){};
\node[circle,draw,fill=black, minimum size=0.1pt,scale=0.8,label=below:{$x^{(1)}_2$}] (DD) at  (750,150){};
\addplot[name path=A,thin,purple,domain={0:1}] {1.2-1.5*(0.16-(x-0.5)^2)^0.25} node[pos=.9, above]{};
\addplot[name path=C,thin,purple,domain={0:1}] {0.8+1.5*(0.16-(x-0.5)^2)^0.25} node[pos=.9, above]{};

\path[->,black](300,20) edge[bend left=0] node [left,yshift=-0.2cm,xshift=-0.11cm] {$\mu^{(2)}_1$~} (400,40);
\path[->,black](800,180) edge[bend left=0] node [above=0.3cm] {$\mu^{(2)}_4$} (800,150);

\addplot[blue!30]fill between[of=A and B, soft clip={domain=0.11:0.555}];
\addplot[blue!30]fill between[of=A and D, soft clip={domain=0.555:0.58}];
\addplot[blue!30]fill between[of=C and E, soft clip={domain=0.67:0.9}];
\end{axis}
\end{tikzpicture}
}\end{center}

\caption{An illustration of the 
SHINE construction  in dimension  $L=2$. Suppose that the measure $\gamma$ is supported on the region enclosed by the red contour, where $\bary(\gamma)=(1,1)$. In the first step of the SHINE construction, we use Proposition~\ref{prop:LdimHx} to find a line $\ell_1$ through $(1,1)$ that partitions $\gamma$ into two parts $\mu^{(1)}_1$ and $\mu^{(1)}_2$, each of whose barycenters lies on the diagonal. In the second step, we find a line $\ell_2$ through $x^{(1)}_1=\bary(\mu^{(1)}_1)$ that partitions $\mu^{(1)}_1$ into two measures $\mu^{(2)}_1$ and $\mu^{(2)}_2$, each of whose barycenters lies on the diagonal, and similarly a line $\ell_3$. We then proceed iteratively. }
\label{fig:3}
\end{figure}

\begin{remark}
    The first step of the construction, i.e., after finishing step $s=0$, already contains a proof of  Proposition~\ref{prop:existence}, because $\delta_{\e}\neq \mu^{(1)}\lcx \gamma$.  Nevertheless, the ideas behind the original proof of  Proposition~\ref{prop:existence} extend to the composite alternative scenario.
\end{remark}

\begin{example}\label{ex:simons}
    Suppose that $L=1$, i.e., we have simple null $P$ versus simple alternative $Q$, where $P\approx Q$. In this case, Proposition~\ref{prop:LdimHx} applies trivially: for each $s\geq 0$ and $1\leq k\leq 2^s$, the measure $\mu^{(s)}_k$ is decomposed into 
    $$\mu^{(s)}_k=\mu^{(s+1)}_{2k-1}+\mu^{(s+1)}_{2k}:=\mu^{(s)}_k\big|_{[0,\bary(\mu^{(s)}_k))}+\mu^{(s)}_k\big|_{[\bary(\mu^{(s)}_k),\infty)}.$$
   % Note that by assumption \eqref{A}, each $\mu^{(s)}_k$ is atomless. 
    As in \eqref{eq:mun}, this results in a sequence of laws $\{\mu^{(s)}\}_{s\geq 0}$ on $\R$ that are increasing and smaller than $\gamma$ in convex order. This is closely related to a martingale decomposition theorem by \cite{simons1970martingale}: if we denote by $\{Z_s\}_{s\geq 0}$ the natural martingale coupling of $\{\mu^{(s)}\}_{s\geq 0}$, then $Z_s\to Z$ a.s.~for some $Z$ that has law $\gamma$. In other words, the e-variable  obtained from the SHINE construction converges to  $  \d Q/\d P$ a.s.~under both $P$ and $Q$.
    \end{example}

\begin{theorem}\label{thm:maximal}Assume \eqref{JA} and \eqref{AC}. 
For any $s$, we have $\mu^{(s)} \lcx \mu^{(s+1)}$, and if $\mu^{(s)}\neq \mu^{(s+1)}$, then the inequality is strict, meaning that the above SHINE construction makes progress at each step. 
Further, assuming \eqref{A}, it produces a sequence of measures that converges almost surely to a maximal element $\mu$ in convex order in $\cM_\gamma$. In this case, if there exists a maximum element $\mu_0$, then the output of our SHINE construction converges to $\mu_0$.
\end{theorem}

% If $\mu^{(s)}\neq \mu^{(s+1)}$, we always have $\mu^{(s)}$ is strictly smaller in convex order than $\mu^{(s+1)}$ and hence the SHINE construction progresses each step before it stops (meaning that $\mu^{(s)}=\mu^{(s+1)}=\cdots$).
When we apply the construction in practice, we need to stop at finitely many steps, so we will not always obtain an exactly maximal element.  {Later in Section \ref{sec:shineconv}, we show that the e-power from the $k$-th step in SHINE converges exponentially to the optimal value produced by SHINE, with a rate that can be made explicit given mild moment conditions.}

We note in particular that Lemma~\ref{lemma:nuatomless} together with Theorem~\ref{thm:maximal} yield that the construction always gives an atomless measure $\mu$ in the limit. 

% When the maximum element exists, a maximal element must be maximum. This has the following consequence.
% \begin{corollary}
% Assume \eqref{A}, \eqref{JA} and \eqref{AC}.    If there exists a maximum element $\mu_0$, then our SHINE construction produces $\mu_0$.
% \end{corollary}

With the presence of atoms, the decomposition given by Proposition~\ref{prop:LdimHx} is not necessarily unique when applied to our construction. The degree of freedom of each $\mu_x$ is the measure on the hyperplane $\partial\H_x$. To describe a well-defined construction, we need to specify $\mu_x|_{\partial \H_x}$   uniquely for each $x$. 
Analyzing the maximality of the output remains a technical task, which we do not discuss in this paper. 

\subsection{Recovering explicitly an e/p-variable}\label{sec:shine2}
We aim first to recover our e-variable $X$, which we recall from Theorem~\ref{prop:Xcondition}  is of the form $X=({\d G}/{\d F})(Y)$, where $Y\in\T((P_1,\dots,P_L,Q),(F,\dots,F,G))$ and $F,G$ come from our SHINE construction. 
We have seen from \eqref{eq:mun} and \eqref{eq:Xn} that  at the $s$-th step, our construction leads to a canonical martingale coupling of $\mu^{(s)}$ and $\gamma$ that couples the mass $\mu^{(s+1)}_{k}(\R^L)\delta_{x^{(s+1)}_k}$ with $\mu^{(s+1)}_{k}$. We denote the martingale coupling  by $(\Lambda_s,\Gamma_s)$, {which is a random vector of dimension $2L$}. Under assumption \eqref{A}, we know further that the measures $\{\mu^{(s)}_k\}_{1\leq k\leq 2^s}$ are mutually singular, and hence $(\Lambda_s,\Gamma_s)$ is backward Monge, i.e., in the backward direction we have $\Lambda_s=h(\Gamma_s)$ for some $h$. Since $(P_1,\dots,P_L,Q)$ is jointly atomless, we may apply Proposition~\ref{prop:commute} to find a simultaneous transport map $Y\in\T((P_1,\dots,P_L,Q),(F,\dots,F,G))$ such that  for each $x\in\X$,$$\frac{\d F}{\d G}(Y(x))\times\e=h\left(\frac{\d P_1}{\d Q}(x),\dots,\frac{\d P_L}{\d Q}(x)\right).$$
This leads to
$$(X(x))^{-1}\times\e=h\left(\frac{\d P_1}{\d Q}(x),\dots,\frac{\d P_L}{\d Q}(x)\right),~x\in\X.$$
For example, the $s$-th step of the construction gives explicitly
\begin{align}
(X(x))^{-1}\times\e=h(x^{(s+1)}_k)~~\text{ if }~~\left(\frac{\d P_1}{\d Q}(x),\dots,\frac{\d P_L}{\d Q}(x)\right)\in \supp\mu^{(s+1)}_{k},~x\in\X.\label{eq:evariableform}
\end{align}

 Note that the measures $F,G$ can meanwhile be reconstructed from Lemma~\ref{lemma:RN}, and further Lemma~\ref{lemma:quantile}(i) if one requires $F=\Leb$. In this case, $Y$ is the valid p-variable as desired, which can be effectively described by the MOT-SOT parity of \cite{wang2022simultaneous}.

\begin{example}
    Suppose that we are in the setting of Example~\ref{ex:d=2}, with $P_1\lawis  \mathrm{N}(-1,1),~P_2\lawis  \mathrm{N}(1,1)$, and $Q\lawis  \mathrm{N}(0,1)$. Recall that
    $$\gamma=\left(\frac{\d P_1}{\d Q},\frac{\d P_2}{\d Q}\right)\Big|_{Q}=\left(e^{-Z-1/2},e^{Z-1/2}\right)\Big|_{Z\lawis \mathrm{N}(0,1)}.$$
By symmetry of $\gamma$, it is clear that the separating hyperplanes $\H_x$ in the SHINE construction are given by $\H_x=\{(a,b):a+b\leq 2x\}$.  In the first step of the construction, we locate the barycenters of the measures $\gamma|_{\H_1}$ and $\gamma|_{\H_1^c}$. By direct calculation, we obtain $\bary(\gamma|_{\H_1})\approx 0.713\times\e$ and $\bary(\gamma|_{\H_1^c})\approx 1.743\times\e$.  Using \eqref{eq:evariableform}, the corresponding e-variable has the  form 
\begin{align*}
    X(x)=\begin{cases}
         1.403&\text{ if }|x|\leq\log(\sqrt{e}+\sqrt{e-1});\\
         0.574&\text{ if }|x|>\log(\sqrt{e}+\sqrt{e-1}).
    \end{cases}
\end{align*}
The resulting e-power $\E^Q[\log X]$ is approximately 0.089. (One may compare this to the maximum e-power 0.12543, which can be directly computed from \eqref{eq:density}.)  In general, we may construct $X$ in multiple steps. 
\end{example}

\subsection{Convergence rate of SHINE}\label{sec:shineconv}
We complement Theorem \ref{thm:maximal} with the following result on the convergence rate of the e-power given by the SHINE construction. Recall from \eqref{eq:opt2} that the e-power is given by $\E^Q[-\log X_k]$, where $\{X_k\}$ is the SHINE martingale. 

\begin{theorem}
    \label{thm:e-power convergence2}
    
Assume the same conditions as in Theorem \ref{thm:maximal}. Suppose that there exists $\ee>0$ such that
\begin{align}
    \E^Q\Big[\Big(\frac{\d P_j}{\d Q}\Big)^{2+\ee}\Big]<\infty, \quad\text{for some }j,\label{eq:moment1}
\end{align}
and  
\begin{align}
    \E^Q\Big[\Big(\frac{\d P_{j'}}{\d Q}\Big)^{-2}\Big]<\infty, \quad\text{for some }j'.\label{eq:moment2}
\end{align}
Consider the e-power $\mathrm{EP}_k:=\E^Q[-\log X_k]$ where $\{X_k\}$ is the SHINE martingale. 
    Then there exist $r\in(0,1)$ and $C>0$ such that
    $$\mathrm{EP}_\infty-\mathrm{EP}_k\leq C r^k,$$
    where $\mathrm{EP}_\infty=\E[-\log X_\infty]$ and $X_k\to X_\infty$ a.s.
\end{theorem}

Our result relies on a particular feature of the SHINE martingale $\{X_k\}$ produced by \eqref{eq:Xn}. Intuitively, the martingale $\{X_k\}$ has a binary tree representation, and the legs in the tree never intersect with other legs at all levels. In this way, one gains control of the fluctuations of the martingale from its values at previous times. The key step to proving Theorem \ref{thm:e-power convergence2} is the following convergence rate of the $L^2$ Wasserstein distance. In particular, this exponential convergence applies to the Simons martingale introduced by \cite{simons1970martingale}; see also Example \ref{ex:simons}.

\begin{lemma}
    \label{thm:general convergence} 

Suppose that the SHINE martingale $\{X_k\}_{k\geq 0}$ satisfies $\E[|X_\infty|^{2+\ee}]<\infty$ for some $\ee>0$ where $X_k\to X_\infty$ a.s. Then there exist $r=r(\ee)<1$ and a constant $C>0$ (where $C$ may depend on $\ee$, the law of $|X_0-X_1|$, and $\E[|X_\infty|^{2+\ee}]$) such that
$$\E[(X_k-X_\infty)^2]\leq Cr^k.$$
If $X$ is uniformly bounded, one can pick $r<0.827$.
\end{lemma}

{Lemma \ref{thm:general convergence} immediately implies $\mathrm{EP}_\infty-\mathrm{EP}_k\leq C r^k$  
for some $r<0.827$ 
if $ \d P_j/\d Q $ for each $j$ is  bounded above and bounded away from $0$.} 
The proofs of Theorem \ref{thm:e-power convergence2} and Lemma \ref{thm:general convergence} are collected in Appendix \ref{sec:rate}.

\subsection{Simulation results}\label{sec:simulation}

% \begin{itemize}
%     \item I have now algorithm working for $L=2$. The way I found the separating hyperplane is by bisecting, i.e. starting from two slopes $0.999,1.001$ and rotate the line until the barycenter is close to the diagonal. If I use other ways to find the line, it typically requires a good guess and could give bad results. This is the reason I don't clearly see how this works in higher dimensions.
%     \item I recovered the pmf for $X$ under $P_1,P_2,Q$ at $k$-th step, but not the explicit form of $X$ as described above. This would need analytic expressions.
%     \item I can plot the $\E^Q[\log X_k]$ at $k$-th step. \com{theoretical value?}
%     \item I can plot the separating hyperplanes but a big ugly. Can also calculate the explicit expressions for those lines.
%     \item In $L=2$, given enough samples, I believe the code works for any distributions.
% \end{itemize}

% \red{Example: Gaussian case: null is -2 and 2, alternative is 0, Bernoulli case: null is 0.3, 0.7 and alternative is 0.5.
% By our theorem, there does not exist exact p/e for just one observation. But there exist such choices if we have two observations.
% }

\begin{figure}[ht!]
\centering
\begin{subfigure}{.5\textwidth}
  \centering
  \includegraphics[width=.99\linewidth]{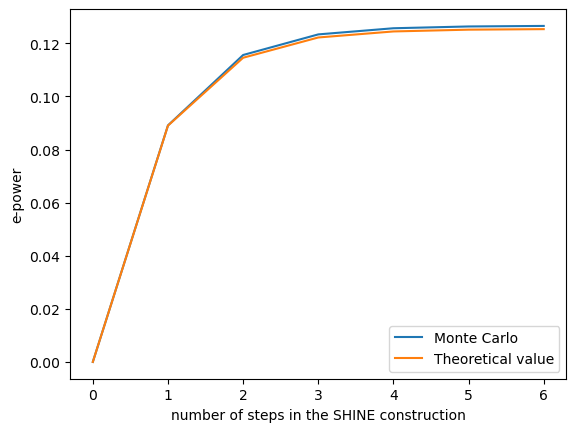}
  \caption{growth of the e-power}
\end{subfigure}%
\begin{subfigure}{.5\textwidth}
  \centering
  \includegraphics[width=.99\linewidth]{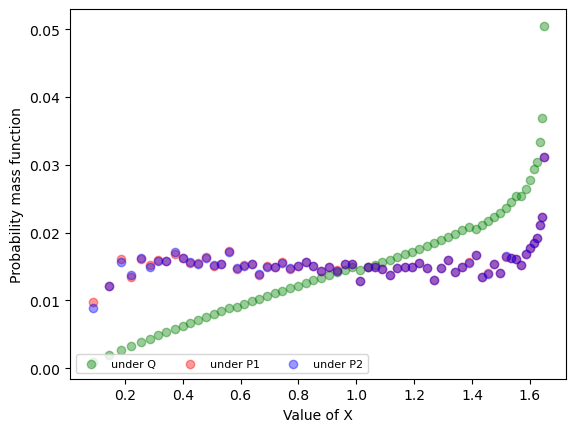}
  \caption{laws of the e-variable $X$ under $P_1,\,P_2,\,Q$}
\end{subfigure}
\caption{The SHINE construction for Example~\ref{ex:d=2}.}
\label{fig:d=2}
\end{figure}

We first consider the setting of Example~\ref{ex:d=2}, where we recall that $P_1\lawis  \mathrm{N}(-1,1),~P_2\lawis  \mathrm{N}(1,1)$, and $Q\lawis  \mathrm{N}(0,1)$. In Figure~\ref{fig:d=2}, we provide two figures illustrating the e-power at each step in the SHINE construction and the corresponding laws of the e-variable under $P_1,~P_2$, and $Q$. In the left panel, we compute the e-power in two ways: from the analytic expression \eqref{eq:density} and by Monte Carlo simulations.  {In the Monte Carlo simulations, we simulate an empirical measure $Q^{(N)}$ of $Q$ and approximate the law $\gamma$ by $\gamma^{(N)}=(\d P_1/\d Q,\d P_2/\d Q)|_{Q^{(N)}}$. After this, we perform the SHINE construction on $\gamma^{(N)}$.} The e-powers are reasonably close with only $N=2\times 10^4$ samples and converge quickly to their limits, where it is straightforward to compute from \eqref{eq:density}  that the theoretical maximum e-power is  approximately 0.12543. In the right panel, we show the distributions of the e-variable under $P_1,~P_2$, and $Q$ at step $s=5$ of the SHINE construction, again by simulating $2\times 10^4$ samples of each distribution. The pivotality of the e-variable implies that the laws of $X$ under $P_1$ and $P_2$ are the same, while the marginal errors shown by the figure are due to our Monte Carlo simulation. Note that within finitely many steps, the SHINE construction always yields a discrete e-variable. With Monte Carlo, our e-variable is \emph{approximately} pivotal since the measure $\gamma$ is atomic thus violating Assumption \eqref{A}.

\begin{figure}[ht!]
\centering
\begin{subfigure}{.5\textwidth}
  \centering
  \includegraphics[width=.99\linewidth]{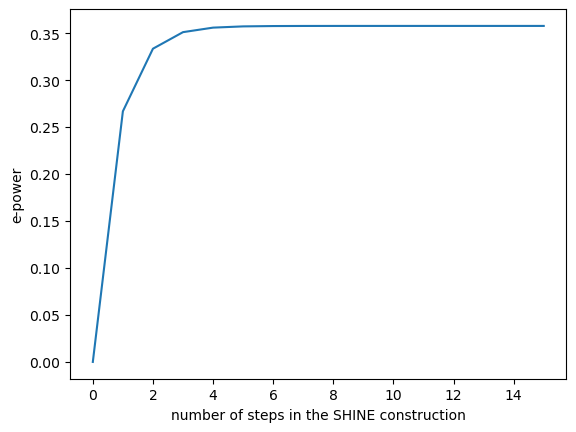}
  \caption{growth of e-power with two data points}
  \label{fig:sub1}
\end{subfigure}%
\begin{subfigure}{.5\textwidth}
  \centering
  \includegraphics[width=.99\linewidth]{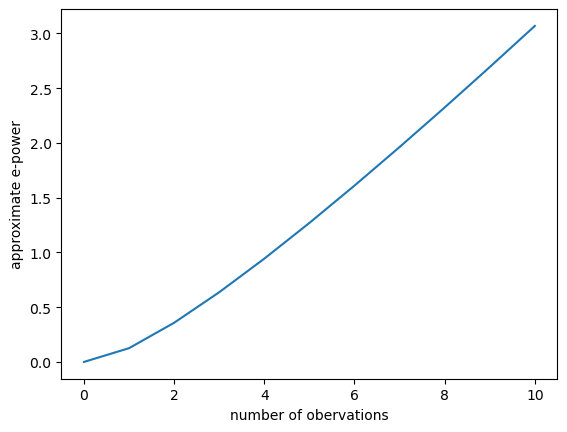}
  \caption{e-power versus number of observations $n$}
  \label{fig:sub2}
\end{subfigure}
\caption{Maximum e-power with multiple data points for Example~\ref{ex:nstep}.}
\label{fig:test}
\end{figure}

In Figure~\ref{fig:test}, we complement the discussions in  Example~\ref{ex:nstep} regarding multiple data points. Recall that $P_1=\mathrm{N}(-\mathbf 1_n,I_n)$,
 $P_2=\mathrm{N}(\mathbf 1_n,I_n)$,
 and  $Q=\mathrm{N}(\mathbf 0_n,I_n)$, where $n\in\N$. Panel (a) computes the theoretical e-power developed after a number of steps with two data points $(n=2)$, which is approximately 0.35775, significantly higher than 0.25086, which is twice the e-power with a single data point. Panel (b) plots the theoretical e-power at the  step $s=7$ of the SHINE construction, for various numbers of observations $n$. Observe that the curve is convex and tends to be linear, reflecting the fact that taking multiple data points increases the average e-power, while the normalized e-power  {$\ell_n/n$} converges as shown in Proposition~\ref{prop:multidata}.

The implementation of the SHINE construction in dimensions greater than two has the obstacle that it is difficult in general to find the separating hyperplanes. We leave this to future work, as well as  generalizations of the SHINE construction when \eqref{A} does not hold.

\section{Composite null and composite alternative}\label{sec:composite alt}
Our goal in this section is to extend Theorems~\ref{thm:iff} and~\ref{thm:p-variable} to composite alternative, i.e., when $|\P|,|\Q|>1$. A full characterization of the existence of (exact and pivotal) nontrivial p/e-variables is provided in the case where both $\P$ and $\Q$ are finite. We also discuss the general case where $\P,\Q$ are infinite, including a few open problems.

\subsection{Existence of an exact and pivotal p/e-variable for the finite case}

We start with the case where $\P,\Q$ are both finite. That is,    given  $\P=\{P_1,\dots,P_L\}$ and $\Q=\{Q_1,\dots,Q_M\}$ such that \eqref{JA} holds, we characterize equivalent conditions for the existence of an (exact and) nontrivial e-variable (or p-variable).

\begin{theorem}\label{thm:iff2}
Assume \eqref{JA}. Suppose that we are testing $\P=\{P_1,\dots,P_L\}$ against $\Q=\{Q_1,\dots,Q_M\}$. The following are equivalent:
    \begin{enumerate}[(a)]
    \item there exists an exact (hence pivotal) and nontrivial p-variable;
        \item there exists a pivotal, exact,  bounded e-variable that has nontrivial e-power against $\Q$;
\item there exists an exact  e-variable that is nontrivial for $\Q$;
\item there exists a random variable $X$ that is pivotal for $\mathcal P$ and satisfies $F\not\in\conv(G_1,\dots,G_M)$, where $F$ is the law of $X$ under every $P\in\P$ and $G_j$ is the law of $X$ under $Q_j$ for $1\leq j\leq M$;
        \item it holds that $\spn(P_1,\dots,P_L)\cap\conv(Q_1,\dots,Q_M)=\emptyset$.
    \end{enumerate}
    Furthermore, the equivalence of (c) and (e) does not require \eqref{JA}.
\end{theorem}

\begin{theorem}\label{thm:pvariable2}
  Assume \eqref{JA}.   Suppose that we are testing $\P=\{P_1,\dots,P_L\}$ against $\Q=\{Q_1,\dots,Q_M\}$. The following are equivalent:
\begin{enumerate}[(a)]

    \item there exists a  nontrivial p-variable;
    \item there exists a bounded e-variable that has nontrivial e-power against $\Q$;
    \item  there exists an e-variable that is  nontrivial for $\Q$;
    \item it holds that $\conv(P_1,\dots,P_L)\cap\conv(Q_1,\dots,Q_M)=\emptyset$.
    \end{enumerate}
    Furthermore, the equivalence of (c) and (d) does not require \eqref{JA} or finiteness of $\P$ and $\Q$.
\end{theorem}

% \begin{remark}\label{rem:ja}By Proposition~\ref{prop:no ja}, the directions (c)$\Leftrightarrow$(e) in Theorem~\ref{thm:iff2} and (c)$\Leftrightarrow$(d) in Theorem~\ref{thm:pvariable2} also hold without condition \eqref{JA}.
% \end{remark}

\begin{remark}
\label{remark:kraft}    The equivalence of $(c)$ and $(d)$ in Theorem~\ref{thm:pvariable2} is a special case of Kraft's theorem, which we recall from \eqref{eq:kraft}. Note that here we do not require that $\P$ and $\Q$ are finite, but only the existence of a reference measure $R$ dominating $\P\cup\Q$. To see that Kraft's theorem  implies the equivalence of $(c)$ and $(d)$ in Theorem~\ref{thm:pvariable2} in case $\P$ or $\Q$ is infinite, suppose that $(c)$ holds. It follows that $\E^Q[X]>1\geq \E^P[X]$ for all $P\in\P$ and $Q\in\Q$.  Kraft's theorem implies the existence of some $\ee>0$ such that $d_{\mathrm{TV}}(\conv\P,\conv\Q)\geq\ee$, and in particular, $(d)$ holds. On the other hand, if $(d)$ is true, then Kraft's theorem yields $\ee>0$ and $X$ satisfying \eqref{eq:kraft}. A suitable linear transformation $Y$ of $X$ then satisfies     $\E^Q[Y]>1\geq \E^P[Y]$ for all $P\in\P$ and $Q\in\Q$ (we may assume $Y$ is positive since $X$ is bounded by construction), and the rest follows from Proposition~\ref{prop:e-var}.
    
    % \com{this needs elaboration: how did e-variables come into Kraft's theorem? maybe reproduce Kraft's theorem, potentially informally, and make this point clearer? also does Kraft's theorem work beyond the finite case?}
    
\end{remark}

% The above characterization extends also to finitely generated measures, which by definition are the  nonnegative mixtures of a finite set of probability  measures, as follows. Given a finite set $\sP$ of probability measures on $\X$, we define 
% $$\M(\sP)=\left\{\sum_{P\in\sP}\lambda_PP~\Big|~\lambda_P\geq 0,~\sum_{P\in\sP}\lambda_P=1\right\}.$$
% We have the following immediate consequence of Theorem~\ref{thm:iff2}.
% \red{start with this, and call it a convex polytope}
\begin{corollary}\label{coro:iff3}  Suppose that we are testing $\P$ against $\Q$, where $\P$ and $\Q$ are convex polytopes in $\Pi$. Denote by $\{P_1,\dots,P_L\}$ (resp.~$\{Q_1,\dots,Q_M\}$) the vertices of the polytope $\P$ (resp.~$\Q$) and assume that $(P_1,\dots,P_L,Q_1,\dots,Q_M)$ is jointly atomless. Precisely the same conclusions in Theorems~\ref{thm:iff2} and~\ref{thm:pvariable2} hold.
\end{corollary}
\begin{proof}
    This follows immediately from Proposition~\ref{prop:polytope}.
\end{proof}

\begin{corollary}\label{coro:nobounded}
    There exists a (pivotal and exact)  e-variable  nontrivial for $\Q$ if and only if  there exists a (pivotal and exact)   e-variable  that has nontrivial e-power against $\Q$.
\end{corollary}

\begin{proof}
This is a direct consequence of     Theorems~\ref{thm:iff2} and~\ref{thm:pvariable2}, and Proposition~\ref{prop:calibration}. 
\end{proof}

\begin{example}
    Fix $0<q_1<q_2<1$ and let $\P=\{\Ber(q_1)\}$ and $\Q=\{\Ber(p) \mid q_2\leq p\leq 1\}$. Corollary~\ref{coro:iff3} then provides an exact nontrivial e-variable (or p-variable). Nevertheless, such an exact nontrivial e-variable (or p-variable) would not exist if we replace $\P$ by $\{\Ber(p) \mid 0\leq p\leq q_1\}$.
\end{example}

Due to the complication of convex order in higher dimensions, it remains a challenging task how to generalize Theorem~\ref{theorem:cx} and the SHINE construction to the composite alternative case.

\subsection{Infinite null and alternative}\label{sec:infinite}

We first state a weaker version of Theorem~\ref{thm:iff2} when both  $\P$ and $\Q$ may be infinite but allow a common reference measure.

\begin{theorem}\label{thm:infinite case}
Assume that there exists a common reference measure $R\in \Pi(\X)$ such that $P\ll R$  for  $P\in\P$ and  $Q\ll R$ for $Q\in\Q$. There exists an exact  bounded e-variable $X$ for $\P$ against $\Q$ satisfying $\inf_{Q\in\Q}\E^Q[\log X]>0$ if and only if $0\not\in\overline{\overline{\spn}\P+\overline{\conv}\Q}$,  where the closure is taken  with respect to  the total variation distance.   If $\Q$ is tight, then we have the further equivalence to $\overline{\spn}\P\cap \overline{\conv}\Q=\emptyset$.
\end{theorem}

Note that we have put a stronger assumption on the e-variable $X$ ($\inf_{Q\in\Q}\E^Q[\log X]>0$) than having nontrivial e-power against $\Q$ (for all $Q\in\Q,~\E^Q[\log X]>0$). Theorem~\ref{thm:infinite case} can thus be seen as a sufficient condition for the existence of an exact  e-variable that has nontrivial e-power against $\Q$. Dealing with pivotal p-variables appears beyond the techniques of this paper.

% \begin{proposition}
% Assume that there exist $\theta_0\in\Theta_0$ and  $\theta_1\in\Theta_1$ such that $P_\theta\ll P_{\theta_0}$  for  $\theta\in\Theta_0$,  $Q_{\theta}\ll Q_{\theta_1}$ for $\theta\in\Theta_1$, and that $\{\P,\Q\}$ is jointly atomless.\footnote{The notion of joint non-atomicity for an infinite collection of measures directly generalizes Definition~\ref{def:ssww}.}  
% Suppose that there exists an exact and nontrivial \red{that may not be pivotal} bounded e-variable for $\P$ against $\Q$. Then $\overline{\spn}\{P_\theta:\theta\in\Theta_0\}\cap \overline{\conv}\{Q_\theta:\theta\in\Theta_1\}=\emptyset$ where the closure is taken with respect to the total variation distance.
% \end{proposition}

% \com{For equivalence of existence of e-variable probably don't need jointly atomless??}

We pose the open problem of characterizing the existence of pivotal, exact, and nontrivial p/e-variables with $\P,\Q$ infinite. For instance, in a very close direction, we pose the 
following conjecture, strengthening Theorem~\ref{thm:infinite case}. We expect that the theory of simultaneous transport between infinite collections of measures will be helpful.

\begin{conjecture}\label{conj}
Suppose that $\P $ and $\Q $ are collections of probability measures on $\X$ with a common reference measure. Assume also that $(\P,\Q)$ is jointly atomless.\footnote{If $\P$ or $\Q$ is infinite, this can be defined in a natural way as in Definition~\ref{def:ssww}.} There exists a pivotal and exact e-variable $X$ satisfying $\inf_{Q\in\Q}\E^Q[\log X]>0$  if and only if $0\not\in\overline{\overline{\spn}\P+\overline{\conv}\Q}$, where the closure is taken  with respect to  the total variation distance. %Moreover, there exists an exact and nontrivial p-variable if and only if $\overline{\spn}\P\cap\conv\Q=\emptyset$.
\end{conjecture}

Our next result shows that surprisingly,  even in  simple settings where $\P$ and $\Q$ are seemingly distant, an exact  e-variable may not exist.

\begin{proposition}\label{prop:id}
Let $P$ be an infinitely divisible distribution on $\R^d$ with a density $p$.  Consider $\P := \{P_\theta\}_{\theta\in\R^d}$ that are the shifts of the measure $P$, where $P_\theta$ has density  $p(x-\theta)$. Let $Q$ be any distribution on $\R^d$ with a density $q$. Then for each $\Q$ that contains $Q$,  there exists no exact  e-variable for $\P$ that is  nontrivial for $\Q$.
\end{proposition}

Note that here we have reached a slightly stronger conclusion than the forward direction of  Theorem~\ref{thm:infinite case}, that even an unbounded e-variable would not exist. The absolute continuity of $Q$ cannot be removed. For instance, if $Q$ has a mass at $x\in\R^d$,  $X=1+\delta_x$ would be an exact e-variable that is nontrivial for $\{Q\}$. 

A particular instance of interest is when $Q$ is Gaussian. In this case, \cite{gangrade2023sequential} proved that for the set of all Gaussians (of all means and all covariances), there does not exist an  e-variable with nontrivial e-power, even non-exact. Thus, our result is stronger in that it allows for a much smaller $\P$ that just includes all translations of any single Gaussian, but it is weaker in that it only shows that an \emph{exact}  e-variable with nontrivial e-power does not exist.

% The same statement holds true if we replace $\Q=\{Q\}$ by $\Q=\{Q_\theta\}_{\theta\in\Theta_1}$ where for some $\theta\in\Theta_1$, $Q_\theta$ is absolutely continuous.
 % \com{I thought you had convex hull not span?}

{We conclude this section with the following example that shows pivotal and exact p/e-values exist for a classic statistical problem. Technically, the construction below does not require any of our previous results, but it leads to a SOT of infinite dimensions.}
 \begin{example}\label{ex:symmetric}
      {Let $\P$ be the class of all symmetric distributions on $\R$ with no mass at $0$, and $\Q$ the class of distributions $Q$ on $\R$ satisfying $Q(\R_+)>1/2$. 
A typical case in applications is to test whether the difference $Y-X$ of pre-treatment measurement $X$
and post-treatment measurement $Y$ 
is symmetric about $0$. 
     Many possible pivotal e-values for $n$ observations can be built based on the signs, the ranks, and the sizes of the data; see \cite{ramdas2020admissible, vovk2022efficiency}.  
     For instance, with one observation, a simple e-value is
     $$X(\omega)=\begin{cases}
         3/2&\quad\text{if }\omega>0;\\
         1/2&\quad\text{if }\omega<0,
     \end{cases}$$
     which is exact. Note that $X$ is also pivotal since $X$ simultaneously maps $\P$ to the uniform distribution on $\{1/2,3/2\}$.
     If one allows for additional randomization using a uniform distribution on $[0,1]$, then 
     $$X(\omega,u)=\begin{cases}
         {u/2}&\quad\text{if }\omega>0;\\
         {(u+1)/2}&\quad\text{if }\omega<0,
     \end{cases}~\quad~u\in[0,1],$$
     is a nontrivial exact p-value.}
 \end{example}

\section{On the existence of nontrivial test (super)martingales}\label{sec:testmartingale}

From here on, for $t\in \{1,2,\dots\}$, let $Z^t$ denote $(Z_1,\dots,Z_t)$, which represents data on $\X^t$, and let $\F$ by default represent the data filtration, meaning that $\F_t=\sigma(Z^t)$.
% \com{I interpret that each data point lives in $\X$. Is that the right interpretation?}

A sequence of random variables $Y \equiv (Y_t)_{t\geq 0}$ is called a \textit{process} if it is adapted to $\F$, that is, if $Y_t$ is measurable  with respect to  $\F_t$ for every $t$. However, $Y$ may also be adapted to a coarser filtration $\G$; for example, $\sigma(Y^t)$ could be strictly smaller than $\F_t$. Such situations will be of special interest to us.
%\com{Let me check my understanding: the filteration becomes coarser because our simultaneous transport may not be injective, right? Yes, that's right.}
Henceforth, $\F$ will always denote the data filtration, and $\G$ will denote a generic subfiltration (which could equal $\F$, or be coarser). An $\F$-stopping time $\tau$ is a nonnegative integer-valued random variable such that $\{\tau \leq t\} \in \F_t$ for each $t \geq 0$. Denote by $\Tau_\F$ the set of all $\F$-stopping times, excluding the constant $0$ and including ones that may never stop. Note that if $\G \subseteq \F$, then $\Tau_\G \subseteq \Tau_\F$.
In this section, $\mathcal P$ is a set of measures on  the sample space $\X^\infty$.\\

% Sometimes, there may not be a single reference measure that locally dominates every distribution in $\P$. This will often be the case for nonparametrically specified $\P$.

\textbf{Test (super)martingales.}
An integrable process $M$ is a \emph{martingale} for $P$ with respect to  $\G$ 
% \com{$\F$ here and below? $\G$ could be $\F$ or smaller, as mentioned just above}  
if 
\begin{equation}\label{eq:martingale}
    \E^P[M_t \mid \G_{t-1}] = M_{t-1} 
\end{equation}
for all $t\geq 1$. $M$ is a \emph{supermartingale} for $P$ if it satisfies \eqref{eq:martingale} with ``$=$'' relaxed to ``$\leq$''.  A (super)martingale is called a \emph{test (super)martingale} if it is nonnegative and $M_0=1$. A process $M$ is called a test (super)martingale for $\P$ if it is a test (super)martingale for every $P \in \P$. The process $M$ is then called a \emph{composite test (super)martingale}.
We say that 
% $M$ is \emph{nontrivial for $\Q$} if $\E^Q[M_\tau]>1$ under all $Q\in \mathcal Q$
% and for all finite stopping times $\tau\in\Tau_\F$, and
 $M$ has \emph{power one against $\Q$} if $\E^Q[\log M_t] \to \infty$ under all $Q\in\Q$.

%We say that $M$ is \emph{nontrivial for $\Q$} if $-M$ is a supermartingale but not a martingale for all $Q\in\Q$, and that $M$ has \emph{nontrivial power against $\Q$} if $\lim_t E^Q[\log M_t]/t > 0$ for all $Q\in\Q$.  These definitions are motivated by the nontriviality of e-variables.

It is easy to construct test martingales for singletons $\P=\{P\}$: we can pick any $Q \ll P$, and then the likelihood ratio process $(\d Q/\d P)(X^t)$ is a test martingale for $P$ (and its reciprocal is a test martingale for $Q$). In fact, every test martingale for $P$ takes the same form, for some $Q$.  

Composite test martingales $M$ are simultaneous likelihood ratios, meaning that they take the form of a   likelihood ratio simultaneously for every element of $\P$. Formally, for every $P \in \P$, there exists a distribution $Q^P\ll P$ and satisfies $M_t = (\d Q^P/\d P)(X^t)$. 
Trivially, the constant process $M_t=1$ is a test martingale for each $\P$, and any decreasing process taking values in $[0,1]$
% \com{What's a decreasing process? a process that almost surely stays constant or goes down}
is a test supermartingale for each $\P$. We call a test (super)martingale \emph{nondegenerate} if it is not always a constant (or decreasing) process.  Nondegenerate test supermartingales do not always exist: their existence depends on the richness of $\P$.\\
% \red{are you omitting the word composite?}

\textbf{On the existence of nondegenerate test (super)martingales.} If $\P$ is too large, there may be no nondegenerate test martingales  with respect to  $\F$.
To explain the situation, suppose that
$\P$ contains only measures of iid sequences with marginal distributions   in a set $\P^{\rm mar}\subseteq \Pi(\X)$.
Examples of the non-existence phenomenon  include the case when $\P^{\rm mar}$ is the set of all mean-zero subGaussian distributions \citep{ramdas2020admissible}, all log-concave distributions~\citep{gangrade2023sequential}, or all Bernoulli distributions~\citep{ramdas2022testing}.   
In all these cases, nondegenerate test martingales have been proven to not exist, at least in the original filtration $\F$.
Sometimes,  nondegenerate test \emph{super}martingales may still exist, as in the subGaussian case. But if $\P^{\rm mar}$  is too large or rich (as in the exchangeable and log-concave cases), even nondegenerate test supermartingales do not exist. 

However, the situation is subtle: in the above situations, there could still exist nondegenerate (or power one) test (super)martingales in some $\G \subseteq \F$. Indeed, for the exchangeable setting described above, \cite{vovk2021testing} constructs exactly such a test martingale in a reduced filtration.  It is a priori not obvious exactly when shrinking the filtration allows for nontrivial test (super)martingales to emerge, and how exactly one should shrink $\F$ (the relevant filtration $\G$ is not evident at the outset).

Our results for (exact) e-variables have direct implications for the existence of test (super)martingales. For simplicity, consider the iid  case, where each $Z_i \lawis  P$ for some $P\in \P^{\rm mar}$ or $P \in \Q^{\rm mar}$; that is, $\P=\{P^\infty \mid  P \in  \P^{\rm mar}\}$ and 
$\Q=\{P^\infty \mid  P \in  \Q^{\rm mar}\}$.

\begin{corollary}\label{coro:testm}
    Let $\P^{\rm mar}$ and $\Q^{\rm mar}$ be subsets of $\Pi(\X)$ allowing for a common reference measure $R\in\Pi(\X)$. If $\overline{\conv}\P^{\rm mar}\cap\overline{\conv}\Q^{\rm mar}=\emptyset$, then there exists a test supermartingale for $\P$ that has power one against  $\Q$. If $0\not\in\overline{(\overline{\spn}\P^{\rm mar}+\overline{\conv}\Q^{\rm mar})}$,  then there exists a test martingale for $\P$ that has power one against $\Q$.  
\end{corollary}

The proof is immediate from Kraft's theorem (see Remark~\ref{remark:kraft}) and Theorem~\ref{thm:infinite case}, and does  not require the joint non-atomicity condition \eqref{JA}. The conditions on $\P$ and $\Q$ imply that an (exact) e-variable (based on $t$ sample points for any $t$) exists for $\P$ that is powerful against $\Q$ by Corollary~\ref{coro:iff3}. We can form our (super)martingale by simply multiplying these e-values for $t=1$ (thus constructively proving the corollary).

We conjecture that the converse direction in the above corollary also holds, perhaps with some additional conditions; in other words, we conjecture that if a test martingale for $\P$ has power one against $\Q$, then the span of $\P^{\rm mar}$ does not intersect $\Q^{\rm mar}$.
(To explain why we cannot directly invoke the reverse directions of our theorems, it is possible that the construction of the e-variable at step $t$ can use information about the distribution gained in the first $t-1$ steps. In short, there (of course) exist test (super)martingales that are not simply the products of independent e-values, and ruling those out requires further arguments, for example, presented in the subGaussian setting by~\cite{ramdas2020admissible}.)

The first (supermartingale) part of Corollary~\ref{coro:testm} is closely related to the main result by~\cite{GrunwaldHK19}, albeit they require some extra technical conditions in their theorem statement while relaxing the polytope requirement. The second (martingale) part is new to the best of our knowledge, and is a key addition to the emerging literature on game-theoretic statistics~\citep{ramdas2022game}.

 \begin{remark}
 Let $ \P^{\rm mar}=\conv(\{P_1,\dots,P_L\})$ with $L$ finite and suppose $Q\in \spn \P^{\rm mar}$ but $Q\not \in \P^{\rm mar}$.
By Theorem~\ref{thm:iff},
there does not exist a nontrivial test martingale for $\P$ against $\{Q^\infty\}$  with respect to  the original filtration. On the other hand, if \eqref{AC} holds, then by   Theorem~\ref{thm:spansquared}, 
there exists a reduced filtration 
 --- in particular, formed by combining data points --- with respect to which a nontrivial test martingale exists.
 \end{remark}

\section{Summary}\label{sec:summary}

This paper uses tools from convex geometry and simultaneous optimal transport to shed light on some fundamental questions in statistics: when can one construct an exact p/e-value  for a composite null, which is nontrivially powerful against a composite alternative? The answer, in the case where the null and alternative hypotheses are convex polytopes in the space of probability measures, is cleanly characterized by convex hulls and spans of the null and alternative sets of distributions. Several other related properties, like pivotality under the null, end up being central.
% where the technical property of joint non-atomicity is assumed in some of our results. 
For general null and alternative hypotheses (which are not polytopes) that allow a common reference measure, we provide a further characterization of the existence of an exact bounded e-variable that has a uniformly positive e-power. 

Our proofs are constructive when the alternative is simple,  and in simple cases, we provide corroborating empirical evidence of the correctness of our theory. A key role is played by the shrinking of the data filtration (accomplished by the  transport map which maps the composite null to a single uniform). Implications for the existence of composite test (super)martingales are also briefly discussed.

We mention some open problems along the way (see Conjecture \ref{conj} and  Sections \ref{sec:multiple observations} and \ref{sec:testmartingale}).
% and anticipate that our results can be generalized beyond convex polytopes with the development of newer technical tools. 
For instance, it is of great interest to extend the SHINE construction to the composite alternative setting.
% and extend our results to general convex subsets of probability measures that are not polytopes.
\begin{acks}[Acknowledgments]
We thank Peter Grunwald,  Martin Larsson, and Johannes Wiesel for helpful discussions. 
 {We are also grateful to three referees for their thoughtful comments and for pointing towards the convergence rate of the SHINE construction.} 
Codes used to generate simulation and numerical results can be found at  \href{https://github.com/Hungryzzy/SHINE}{https://github.com/Hungryzzy/SHINE}.
\end{acks}

%% if your bibliography is in bibtex format, uncomment commands:
\bibliographystyle{imsart-nameyear.bst} % Style BST file (imsart-number.bst or imsart-nameyear.bst)
\bibliography{Composite.bib}       % Bibliography file (usually '*.bib')

%% or include bibliography directly:
\newpage

%%%%%%%%%%%%%%%%%%%%%%%%%%%%%%%%%%%%%%%%%%%%%%
%% Example with single Appendix:            %%
%%%%%%%%%%%%%%%%%%%%%%%%%%%%%%%%%%%%%%%%%%%%%%
\begin{appendix}

\section{General relations on the existence of p- and e-variables} \label{sec:general}
For convex polytopes $\P$ and $\Q$ in $\Pi$, we may write $\P=\conv \widetilde{\P}$ and $\Q=\conv\widetilde{\Q}$ where $\widetilde{\P}$ and $\widetilde{\Q}$ are finite. The following result  helps us to reduce the problems to the case where $\P$ and $\Q$ are finite.
\begin{proposition}
    \label{prop:polytope}
    Suppose that $\P=\conv \widetilde{\P}$ and $\Q=\conv\widetilde{\Q}$. 
    \begin{enumerate}[(i)]
        \item There exists an (exact) nontrivial p-variable for $\P$ and $\Q$ if and only if the same exists for $\widetilde{\P}$ and $\widetilde{\Q}$.\footnote{Here and later, we mean that the statement holds regardless of whether the bracketed constraint exists, i.e., the current sentence contains two (different) statements.}
        \item There exists a (pivotal, exact, bounded) e-variable that is nontrivial for (or has nontrivial e-power against) $\Q$ for $\P$ and $\Q$ if and only if the same exists for $\widetilde{\P}$ and $\widetilde{\Q}$. 
    \end{enumerate}
\end{proposition}

\begin{proof}
   This is clear from definitions of p/e-variables.
\end{proof}

As a result of the above proposition, in what follows, we can concern ourselves, without loss of generality, with the case where $\P$ and $\Q$ are finite subsets of $\Pi(\X)$ (except for Section~\ref{sec:infinite}). Recall again from Section~\ref{sec:intro} the difference between a nontrivial e-variable, and one with nontrivial e-power.

\begin{proposition}
    \label{prop:e-var}
  Let $X$ be a (pivotal and exact) bounded e-variable for $\P$ that satisfies $\inf_{Q\in \Q}\E^Q[X]>1$. Then there exists a (pivotal and exact) bounded  e-variable for $\P$ that satisfies $\inf_{Q\in \Q}\E^Q[\log X]>0$. In particular, suppose that  $\P$ and $\Q$ are both finite and $X$ is a (pivotal and exact) bounded e-variable for $\P$ that is nontrivial for $\Q$. Then there exists a (pivotal and exact) bounded  e-variable for $\P$ that has nontrivial e-power against $\Q$. 
\end{proposition}

\begin{proof}
 The following fact is crucial: by the Taylor expansion of the log function, for every $\ee>0$, there exists $\delta>0$, such that
 for each $x\in[1-\delta,1+\delta]$, $(1-\ee)(x-1)\leq \log x\leq (1+\ee)(x-1)$. 
 Note that each $Y=(1-b)+bX$ with $b>0$  is clearly an e-variable. On the other hand, since $X$ is bounded, the range of $Y$ can be chosen arbitrarily close to $1$ by picking $b$ small enough. Using $\inf_{Q\in\Q}\E^Q[X]>1$ we see that with $b$ small enough, $Y$ is an e-variable that has nontrivial e-power against $\Q$. Note that pivotality and exactness are preserved under this transformation.
\end{proof}

\begin{remark}
Assuming $(\P,\Q)$ is jointly atomless,    Proposition~\ref{prop:e-var} also holds true without the boundedness assumption on $X$; see Corollary~\ref{coro:nobounded} below. However, we are not aware of a simpler proof of this fact.
\end{remark}

In the sequel, when the equivalence of the existences is clear, we may write ``there exists a nontrivial e-variable'' instead.  When $\Q$ is infinite, these two definitions are in general different, as shown by the following example.

\begin{example}
    Let $Z_\mu$ denote the law $\mathrm{N}(\mu,1)$ for $\mu\in\R$, and consider  $\P=\{Z_0\}$ and $\Q=\{Z_\mu\mid\mu>0\}$. Clearly, $X(\omega)=1/2+\bone_{\{\omega>0\}}$ is a bounded e-variable that is nontrivial for $\Q$. Suppose for contradiction that $Y$ is a bounded e-variable that has nontrivial e-power against $\Q$. Since $Y$ cannot be a constant, $\E^{Z_0}[\log Y]<\E^{Z_0}[Y]-1=0$. Since $Z_\mu\to Z_0$ in total variation as $\mu\to 0$, we have for $\mu>0$ small enough that $\E^{Z_\mu}[\log Y]<0$, contradicting $\E^Q[\log Y]>0$ for all $Q\in\Q$. 
\end{example}
%truncate Y from \delta>0 small is enough. we omitted this step here.

The following calibration result is in place to help us construct an e-variable based on  a p-variable.

\begin{proposition}\label{prop:calibration}
    Suppose that $\Q$ is finite.  
    \begin{enumerate}[(i)]
        \item If there exists an exact (hence pivotal) and nontrivial p-variable, then there exists a pivotal, exact, and bounded e-variable with nontrivial e-power against $\Q$.
        \item If there exists a nontrivial p-variable, then there exists a bounded e-variable with nontrivial e-power against $\Q$.
    \end{enumerate}
\end{proposition}

\begin{proof}
(i) Suppose that $X$ is an exact nontrivial p-variable. It follows that $E:=2-2X$ is a pivotal, exact, and bounded e-variable, and $\E^Q[E]>1$ for each $Q\in\Q$. Proposition~\ref{prop:e-var} then finishes the proof. (ii) is similar, where we recall that without loss of generality, a p-variable takes values in $[0,1]$. 
\end{proof}

Our next simple result provides general necessary conditions for the existence of p/e-variables, hence answering the trivial parts of \textbf{(Q-existence)}.
\begin{proposition}\label{prop:nec}
Suppose that $\P$ and $\Q$ are arbitrary subsets of $\Pi(\X)$.
    \begin{enumerate}[(i)]
        \item If there exists a nontrivial e-variable  for $\Q$, then $\conv\P\cap\conv\Q=\emptyset$.
        \item If there exists an exact and nontrivial e-variable  for $\Q$, then $\spn\P\cap\conv\Q=\emptyset$.
    \end{enumerate}
\end{proposition}

\begin{proof}
    For (i), suppose that $R\in \conv\P\cap\conv\Q$, then since $\E^P[X]\le 1$ for all $P\in \P$, we have $\E^R[X]\leq 1$. But $\E^Q[X]>1$ for all $Q\in\Q$ implies $\E^R[X]>1$, yielding a contradiction. For   (ii), suppose that $R\in \spn\P\cap\conv\Q$, then  $\E^P[X]= 1$ for all $P\in \P$ gives that $\E^R[X]= 1$. But $\E^Q[X]>1$ for all $Q\in\Q$ gives $\E^R[X]>1$, yielding a contradiction.
\end{proof}

Let us end this section by incorporating the following important result, which sometimes helps us remove the jointly atomless condition when pivotality is not involved.
\begin{proposition}\label{prop:no ja}Fix $\P$ and $\Q$. If there exists  an (exact) e-variable $($defined on $(\X\times[0,1],\F\otimes\B([0,1])))$ that is nontrivial for $\{Q\times\Leb\}_{Q\in\Q}$ with null $\{P\times\Leb\}_{P\in\P}$, then there exists  an (exact) e-variable  that is nontrivial for $\Q$ with null $\P$.
\end{proposition}

\begin{proof}
Let $Y$ be an exact  e-variable  that is nontrivial for $\{Q\times\Leb\}_{Q\in\Q}$ with null $\{P\times\Leb\}_{P\in\P}$. Define $X=\E^\Leb[Y]$ by taking the expectation of $Y$ over the second coordinate. Then $\E^{P_i}[X]=\E^{P_i\times\Leb}[Y]=1$ and $\E^{Q_j}[X]=\E^{Q_j\times\Leb}[Y]>1$, meaning that $X$ is an exact e-variable nontrivial for $\Q$. The non-exact case is similar.    
\end{proof}

\section{Proof of results from Section \ref{sec:simple alt}}\label{App:B}

\begin{proof}[Proof of Theorem~\ref{thm:p-variable}]That (a)$\Rightarrow$(b) is precisely Proposition~\ref{prop:calibration}; (b)$\Rightarrow$(c) is clear, and (c)$\Rightarrow$(d) is immediate from Proposition~\ref{prop:nec}.  
%The ``if" follows from Hahn-Banach: let $\mu$ be a dominating measure for $P_1,\dots,P_L,Q$. Since $\frac{\d Q}{\d \mu}\not\in\conv(\frac{\d P_1}{\d \mu},\dots,\frac{\d P_L}{\d \mu})$, there exists a bounded random variable $Y$ with $\E^{P_i}[Y]=\E^\mu[\frac{\d P_i}{\d \mu}Y]\leq 0<\E^\mu[\frac{\d Q}{\d \mu}Y]=\E^Q[Y]$ for all $i$. Then let $X=1+Y/\sup{|Y|}$. To see that $\E^Q[\log X]>0$, we simply replace $X$ by $(1-b)+bX$ with $b>0$ chosen  small enough, because our choice of $X$ is bounded.
 For (d)$\Rightarrow$(a), we define the set $D=D_L=(-\infty,1)^L\cup (1,\infty)^L\cup\{\e\}$. We claim that it suffices to find a measure $\mu\lcx\gamma$ that is supported on $D$ and not equal to $\delta_{\e}$. Given such $\mu$, we apply Lemma~\ref{lemma:RN} with $d=L$ to find measures $F_1,\dots,F_L$ such that $({\d F_1}/{\d \Leb},\dots,{\d F_L}/{\d \Leb})|_{\Leb}= \mu.$ Since $\mu$ is supported on $D$, we may without loss assume that there is a threshold $\beta\in(0,1)$ such that for each $1\leq i\leq L$, $\d F_i/\d \Leb\leq 1$ on $[0,\beta)$ and $\d F_i/\d \Leb\geq 1$ on $(\beta,1]$. In particular, $F_i\sst\Leb$. Proposition~\ref{prop:ssww} then yields a random variable $X\in\T((P_1,\dots,P_L,Q),(F_1,\dots,F_L,\Leb))$. Let $\Psi$ be as given in Lemma~\ref{lemma:quantile}(ii). By definition, $\Psi\circ X$ is a nontrivial p-variable.

To  find a measure $\mu\lcx\gamma$ that is supported on $D$ and not equal to $\delta_{\e}$, for simplicity we translate $D$ by $\e$, and from now on $D=(-\infty,0)^L\cup(0,\infty)^L\cup\{\z\}$ and $\gamma$ has mean $\z$. Our goal is to find a measure $\mu$ supported on $D$ such that $\mu\lcx\gamma$ and $\mu\neq\delta_{\z}$. We apply induction on $L$. Suppose that $L=2$. Then since $Q\not\in\conv(P_1,P_2)$, the measure $\gamma$ is not supported on any line that has a negative slope and contains $\z$. There are two cases.
\begin{itemize}
    \item If $\gamma$ is not supported on any line (hyperplane in $\R^2$), then $Q\not\in\spn(P_1,P_2)$. By Theorem~\ref{thm:iff}, a  nontrivial p-variable exists.
    \item If $\gamma$ is supported on a line, then such a line must contain $\z$ and have a positive slope, and hence is contained in   $D$. 
\end{itemize}

Now suppose that $L> 2$. We say a set $K\subseteq\R^L$ is a linear cone if it is the union of a convex cone and its symmetric image around $\z$ in $\R^L$. Clearly, $D$ is a linear cone, and
\begin{enumerate}[(i)]
    \item the intersection of a subspace and a linear cone is a linear cone;
    \item if $S$ is a subspace of $\R^L$ and $K$ is a linear cone, then $\{\z\}\subsetneq S\cap K$ if and only if there exists a one-dimensional subspace $T$ of $S$ such that $T\subseteq K$.
\end{enumerate}
Since $Q\not\in\conv(P_1,\dots,P_L)$, the measure $\gamma$ is not supported on any hyperplane that is contained in $D^c\cup\{\z\}$ and contains $\z$. If $\gamma$ is not supported on any hyperplane, then using Theorem~\ref{thm:iff}, a  nontrivial p-variable exists. Thus we may assume that $\gamma$ is supported on some hyperplane $S$ such that $\{\z\}\subsetneq S\cap D$. We then lower the dimension by one and identify $S=\R^{L-1}$. There are two cases.
\begin{itemize}
    \item If $\gamma$ is not supported on any hyperplane in $S$, then there exist points $x_1,\dots,x_{L-1}\in\supp\gamma$ such that $\ri\aff\{x_1,\dots,x_{L-1}\}=S$. Let $T$ be a one-dimensional subspace of $S$ such that $T\subseteq D$. Then $T\cap\conv\{x_1,\dots,x_{L-1}\}$ is nonempty.
    Using Lemma~\ref{lemma:CM}, we may find a measure $\mu$ supported on the bounded set  $T\cap\conv\{x_1,\dots,x_{L-1}\}$ (and thus supported on $D$) such that $\mu\lcx\gamma$ and $\mu\neq\delta_{\z}$. It follows that a  nontrivial p-variable exists.
    \item If $\gamma$ is supported on a hyperplane $S'$ of $S$, then by Mazur's separation theorem  (\citet[Corollary 3.4]{conway1990course}) and since $\gamma$ is not supported on any hyperplane that intersects with $D$ only at $\z$, we must have  $\{\z\}\subsetneq S'\cap D$. In this case, we have reduced the dimension by one. Thus induction works for this case.
\end{itemize}
By reducing the problem iteratively in the above manner, we eventually   arrive at the problem with $L=2$, which we already showed above. 
\end{proof}

\section{Proof of results from Section \ref{sec:e}}\label{app:sec4}

\begin{proof}[Proof of Theorem \ref{prop:Xcondition}]
Suppose that $Z$ is a maximizer to \eqref{eq:opt}. Since $Z$ is a pivotal e-variable, we denote by $F'$ as the common distribution of $Z$ under $P_i,~1\leq i\leq L$, and $G'$ the distribution of $Z$ under $Q$. Let $\widetilde{Z}$ be the identity random variable on $\R$, we have $\E^{F'}[\widetilde{Z}]=1$. By Gibbs' inequality,
$$\E^Q[\log Z]=\E^{G'}[\log\widetilde{Z}]\leq \E^{G'}\Big[\log\frac{\d G'}{\d F'}\Big]=\E^Q\Big[\log\Big(\frac{\d G'}{\d F'}(Z)\Big)\Big].$$
Thus, $X=({\d G'}/{\d F'})(Z)$ must also be a maximizer to \eqref{eq:opt11}. 
\end{proof}

\begin{proof}[Proof of Proposition \ref{prop:LdimHx}]
We induct on $L$. The base case is $L=1$, where the claims follow simply by picking $\H_x=(-\infty,x]$.

% We first illustrate the case $L=2$.  In the following, we denote by $\H_+=\{(x,y)\in\R^2|y\geq x\}$ and $\H_-=\{(x,y)\in\R^2|x\geq y\}$. 
% Let $x\geq 0$. It follows that the $\bary(\gamma|_{\H_+})\in\H_+$ and  $\bary(\gamma|_{\H_-})\in\H_-$. By continuity, there exist $\H_x$ and $\mu_x$ satisfying (iii) and (iv) above. 
% Suppose that $\I^+\subseteq \H_x$. Then $\gamma$ is concentrated on $\I^+$. In this case, we may then replace $\H_x$ by $(-\infty,x]\times\R$ and $\mu_x$ by $\mu_{(-\infty,x]\times\R}$. Thus (i) holds.
% (ii) can then be assumed true without loss of generality.
% The measure $\mu_x$ is now  determined up to masses on $\H_x\cap\H_x^c$, but we have assumed that $\gamma(\H_x\cap\H_x^c)=0$. Thus the measure $\mu_x$ is unique.

% The general case of $L\geq 2$ can be proved via induction on $L$.  We induct on $L$ to produce $\H_x$ and $\mu_x$ that satisfy (iii) and (iv) above, while the rest follows from similar arguments as in the base case of $L=2$.

Fix an arbitrary $L\geq 2$ and $x\geq 0$. Consider the plane $\P_L=\{\bx\in\R^L \mid \,x_1=x_2\}\subseteq\R^L$, so that $\I^+\subseteq \P_L$. The collection of lines in $\P_L$ through $x\e$ will be denoted by $\mathcal L_\theta,~\theta\in[0,2\pi)$, where $\I^+\subseteq \mathcal L_0$. For each $\mathcal L_\theta$, consider the projection of $\gamma$ on the hyperplane $P_\theta$ to which $\mathcal L_\theta$ is normal. It follows from our induction hypothesis that there is some half-space $\H_{\theta,x}$ of $\R^L$ on which some measure $\mu_{\theta,x}\leq\gamma$ is supported, such that $\bary(\mu_{\theta,x})\in\P_L$ and $\bary(\gamma-\mu_{\theta,x})\in\P_L$, as well as $\mathcal L_\theta\subseteq \partial\H_{\theta,x} $.

Suppose that (i) does not hold. Then $\gamma$ is supported on a hyperplane $S$ in $\R^L$ containing $\I_L$. By the induction hypothesis, we may find a closed half-space $\H_x'$ of $S$ satisfying the conditions (i)-(iv). Clearly, any closed half-space $\H_x$ of $\R^L$ containing $\H_x'$ also satisfies the same conditions.

Therefore, we may assume  (i)  and  that $\gamma$ is not supported on any hyperplane in $\R^L$. In particular, $\mu_{0,x}$ and $\mu_{\pi,x}$ are non-zero. In this case, $\bary (\mu_{0,x})$ and $\bary(\mu_{\pi,x})$ lie in the two different half-planes in $\P_L$ separated by $\I^+$. By continuity of the measure, there exists some $\theta_x$ such that $\bary(\mu_{\theta_x,x})\in\I^+$.  This establishes (iii) and (iv). 
Finally, by replacing $\H_x$ by $\H_x^c$, we may assume that (ii) holds as well.

Suppose that \eqref{A} holds and $\mu_x,\nu_x$ are distinct measures satisfying the above conditions. Then $\mu_x-\nu_x$ is a nontrivial signed measure supported on a hyperplane in $\R^L$, contradicting \eqref{A}.
\end{proof}

The proof of Theorem Theorem~\ref{thm:maximumexist} requires the following lemma.
\begin{lemma}\label{lemma:3}
    Let $\nu$ be a probability measure on $\I^+_L$ such that $\nu\lcx\gamma$. For $\H_x$ and $\mu_x$ defined in Proposition~\ref{prop:LdimHx}, denote by   $b_x\e$ the barycenter of $\mu_{x}$ and $\xi$ distributed as the first marginal of $\nu$. Then $\E[(\xi-x)_-]\leq (x-b_x)\mu_{x}(\R^L)$ for all $x\geq 0$.
    Moreover, equality holds for $x$ if and only if for every martingale coupling $(X,Y)$ such that $X\lawis\nu$ and $Y\lawis\gamma$, it holds $(Y \mid X\leq x)\lawis \mu_x/\mu_x(\R^L)$.
\end{lemma}

\begin{proof}
Let $v_x$ be a unit normal vector to $\partial\H_x$, such that the angle $\theta_x$ between the vectors $v_x$ and $\e$ satisfies $0\leq \theta_x<\pi/2$. For $y\in\R^L$ we write $a_y=\langle y,v_x\rangle$ with Euclidean inner product. Define $\phi_x:\R^L\to\R$ by $$\phi_x(y):=
\begin{cases}{a_y}/{\cos\theta_x}&\text{ if }a_y\leq 0;\\ 0&\text{ if }a_y>0.
\end{cases}$$
Since $\phi_x$ is concave and $\nu\lcx\gamma$, it follows that
$$(b_x-x)\mu_{x}(\R^L)=\int\phi_x\d\gamma\leq\int\phi_x\d\nu=-\E[(\xi-x)_-].$$This completes the proof. The rest is clear.
\end{proof}

\begin{proof}[Proof of Theorem~\ref{thm:maximumexist}] We first prove (b)$\Rightarrow$(a).  We first characterize the measure $\mu$ by the cumulative density function of its first marginal (recall that $\mu$ is supported on the nonnegative diagonal $\I^+$). For $x\geq 0$, pick $\mu_x$ as in Proposition~\ref{prop:LdimHx}. Note that $x\mapsto \mu_x(\R^L)$ is nondecreasing in $x$ and continuous by (b).
Define $\mu$ by the unique probability measure on $\I^+$ such that $\mu([0,x]^L)=\mu_x(\R^L)$. 

We next show that $\mu\lcx\gamma$. By Strassen's theorem, it suffices to find a martingale coupling $(X,Y)$ such that $X\lawis\mu$ and $Y\lawis\gamma$. Let us fix $X\lawis\mu$ and let $$(Y \mid x<X\leq x')\lawis \frac{\mu_{x'}-\mu_x}{\mu_{x'}(\R^L)-\mu_x(\R^L)},$$ where we identify the random variable $X$ supported on $\I$ with its first coordinate. This defines a coupling $(X,Y)$ since $x\leq x'\implies \mu_x\leq \mu_{x'}$. 
Let $a\geq 0$ be arbitrary. On the event $\{X\leq a\}$, $Y$ is distributed as $\mu_a/\mu_a(\R^L)$. By Proposition~\ref{prop:LdimHx}, we have $\E[Y\bone_{\{X\leq a\}}]\times \e\in \I^+$. In addition, 
$$\E\left[X\bone_{\{X\leq a\}}\right]\times \e=\left(\int_0^a x\,\d(\mu_x(\R^L))\right)\times \e,$$
\sloppy which is exactly  $\bary(\mu_a)$ projected to $\I$. Therefore, we must have $\E[X\bone_{\{X\leq a\}}]=\E[Y\bone_{\{X\leq a\}}]$, so that $(X,Y)$ is indeed a martingale. Thus $\mu\lcx\gamma$.

Now by Lemma~\ref{lemma:cxproperty}(i),  it suffices to show that for each $\nu\lcx\gamma$ and $\xi_\nu,\xi_\mu$ denoting the first marginals of $\nu,\mu$, it holds that $\E[(\xi_\nu-x)_-]\leq \E[(\xi_\mu-x)_-]$ for all $x\geq 0$. This is indeed a consequence of Lemma~\ref{lemma:3}, since 
$$\E[(\xi_\mu-x)_-]=x\mu_{x}(\R^L)-\int x\,\d(\mu_x(\R^L))=(x-b_x)\mu_{x}(\R^L).$$

We next show (a)$\Rightarrow$(b). Suppose that $x<y$ but $\mu_x\not\leq\mu_y$.  In particular, $\supp\mu_x\not\subseteq\supp\mu_y$.  We may assume that $\mu_x(\R^L)$ and $\mu_y(\R^L)$ are positive.  Suppose for contradiction that (a) holds with a maximum element $\mu$. We define the measures
$$\nu_x=\mu_x(\R^L)\delta_{\bary(\mu_x)}+(1-\mu_x(\R^L))\delta_{\bary(\gamma-\mu_x)}$$
and similarly $\nu_y$. 
Then with the usual notation, 
\begin{align}
    \E[(\xi_\mu-x)_-]\geq \E[(\xi_{\nu_x}-x)_-]=(x-b_x)\mu_{x}(\R^L)\label{eq:nux}
\end{align}
and
\begin{align}
    \E[(\xi_\mu-y)_-]\geq \E[(\xi_{\nu_y}-y)_-]=(y-b_y)\mu_{y}(\R^L).\label{eq:nuy}
\end{align}
By Lemma~\ref{lemma:3}, equalities hold for \eqref{eq:nux} and \eqref{eq:nuy}, and for every martingale coupling $(X,Y)$ such that $X\lawis\mu$ and $Y\lawis\gamma$, it holds $(Y \mid X\leq x)\lawis \mu_x/\mu_x(\R^L)$ and $(Y \mid X\leq y)\lawis \mu_y/\mu_y(\R^L)$. This contradicts $\supp\mu_x\not\subseteq\supp\mu_y$.
\end{proof}

The following supporting lemma will prove useful in proving Theorem \ref{theorem:cx}.

\begin{lemma}\label{lemma:2}
Suppose that \eqref{A} holds and there exists a convex set $\Gamma\subseteq \R^{2}$ such that $\gamma$ is supported on   $\partial\Gamma$. For $x\geq 0$, let $\H_x$ and $\mu_x$ be defined as  in Proposition~\ref{prop:LdimHx}. Then for $0\leq x\leq x',~\mu_x\leq \mu_{x'}$. 
\end{lemma}

\begin{proof}
Fix $0\leq x<x'$. Define  $C_1=\H_{x}\setminus \H_{x'}^\circ$  and $C_2=\H_{x'}\setminus \H_x^\circ$.  It follows that the positive part of $\mu_x-\mu_{x'}$ is supported on $C_1$. Let us define $$S=\R\e+\partial\H_x\cap\partial\H_{x'}.$$ The line $S$ separates $\R^2$ into two (closed) half-spaces, and we denote by $\H'$ the one that does not contain $C_1$. Since the barycenters of $\mu_x$ and $\mu_{x'}$ lie on $\I^+$, it suffices to show that $\gamma(C_2\setminus (C_1\cup\H'))=0$. Suppose not. Then there exist $z_1\in \partial\Gamma\cap C_1$ and $z_2\in\partial\Gamma\cap C_2\setminus (C_1\cup\H')$. Since $\Gamma$ is  convex, it cannot hold that both $x\e$ and $x'\e$ belong to $\Gamma^\circ$. Suppose that $x\e\not\in\Gamma^\circ$. Then by convexity of $\Gamma$ and our assumption \eqref{A}, we have $\mu_x=0$, thus $\mu_x\leq \mu_{x'}$ holds trivially. The case $x'\e\not\in\Gamma^\circ$ is similar. 
\end{proof}

\begin{proof}[Proof of Theorem~\ref{theorem:cx}]
    The first claim follows from Theorem~\ref{thm:maximumexist} and Lemma~\ref{lemma:2}. The existence of $F,G$ follows from Lemma~\ref{lemma:RN} with $d=L$, by setting $G=\Leb$. 
\end{proof}

For the proof of Theorem \ref{thm:spansquared}, we need the following lemma on independence of powers of probability measures.

\begin{lemma}\label{lem:remark:power}
For distinct measures $P_1,\dots,P_L$ on an Euclidean space $\X$, there exists $k\geq 1$ such that $P_1^k,\dots,P_L^k$ are linearly independent probability measures on $\X^k$.
\end{lemma}
\begin{proof}
%[Proof of Lemma~\ref{lem:remark:power}]  
For notational brevity, we assume $\X=\R$, where the general case of higher-dimensional spaces follows from  essentially the same argument. By assumption, the characteristic functions of the probability laws $P_1,\dots,P_{L}$ are distinct, and hence there exists  $\bl=(\lambda_1,\dots,\lambda_{L-1})\in\R^{L-1}$ such that   the expectations of the $2(L-1)$-dimensional random vector
$$\psi(\bl,X):=(\cos(\lambda_1 X),\sin(\lambda_1X),\dots,\cos(\lambda_{L-1} X),\sin(\lambda_{L-1}X))$$are distinct under the laws $X\lawis P_\ell,~1\leq \ell\leq L$.\footnote{Here, we are using the elementary principle that for $L$ distinct functions on a common domain, there exist $L-1$ points in the domain  such that no two functions agree on all of the $L-1$ points. This can be proved using induction.}
Let us consider a number $k\geq 1$ that will be eventually picked large enough. Define the following  (simultaneous) transport map 
$$T:\R^k\to\R^{2(L-1)k},~(x_1,\dots,x_k)\mapsto (\psi(\bl,x_1),\dots,\psi(\bl,x_k)),$$and denote by $Q_\ell$ the pushforward of $(P_\ell)^k$ under $T$, for $1\leq \ell\leq L$.\footnote{Our notation omits the dependence of $T$ and $Q_\ell$ on $k$.} Note that each $Q_\ell$ is a law on $\R^{2(L-1)k}$, whose coordinates are denoted by $\bx=(x_{i,j})_{1\leq i\leq 2(L-1),1\leq j\leq k}$. By law of large numbers, the law of $Q_\ell$ is concentrated near the set 
\begin{align*}
    R_\ell&:=\Bigg\{\bx:\Big(\frac{1}{k}\sum_{j=1}^kx_{2i-1,j},\frac{1}{k}\sum_{j=1}^kx_{2i,j}\Big)\\
    &\hspace{2cm}=\left(\E^{P_\ell}[\cos(\lambda_iX)],\E^{P_\ell}[\sin(\lambda_iX)]\right),~\forall 1\leq i\leq L-1\Bigg\} .
\end{align*}
Since $\E^{P_\ell}[\psi(\bl,X)]$ are distinct for $1\leq\ell\leq L$, the sets $R_\ell$ are disjoint. 
More precisely, using a quantitative version of the  central limit theorem (e.g., Theorem 1.3 of \cite{talagrand1996new}), we know that for any $\ee>0$, there exists $C=C(\ee,\P)>0$ such that for $k\geq C$, there are disjoint sets $\{A_\ell\}_{1\leq\ell\leq L}$ (each $A_\ell$ contains $Q_\ell$) such that each $Q_\ell$ is concentrated in $A_\ell$, in the sense that
\begin{itemize}
    \item for each $\ell$, $Q_\ell(A_\ell)>1/2$;
    \item for each $\ell'\neq\ell$, $Q_\ell(A_{\ell'})<\ee$.
\end{itemize}
For example, we may take 
\begin{align*}
    A_\ell:=&\Bigg\{\bx:\Big(\frac{1}{k}\sum_{j=1}^kx_{2i-1,j},\frac{1}{k}\sum_{j=1}^kx_{2i,j}\Big)\in\big(\E^{P_\ell}[\cos(\lambda_iX)]-k^{-1/3},\E^{P_\ell}[\cos(\lambda_iX)]\\
    &+k^{-1/3}\big)\times\left(\E^{P_\ell}[\sin(\lambda_iX)]-k^{-1/3},\E^{P_\ell}[\sin(\lambda_iX)]+k^{-1/3}\right),~\forall 1\leq i\leq L-1\Bigg\}.
\end{align*}
We next show that the probability measures $\{Q_\ell\}_{1\leq \ell\leq L}$ are linearly independent for $\ee=1/(2L)$. Indeed, suppose that $\alpha_1Q_1+\dots+\alpha_LQ_L=0$. Then for each $1\leq i\leq L$, by the triangle inequality,
$$\frac{1}{2}|\alpha_i|<|\alpha_i|Q_i(A_i)=\Big|\sum_{\ell\neq i}\alpha_\ell Q_\ell(A_i)\Big|\leq \sum_{\ell\neq i}|\alpha_\ell| Q_\ell(A_i)<\ee\sum_{\ell=1}^L|\alpha_\ell|.$$
Summing the above terms  over $1\leq i\leq L$ yields $2L\ee>1$, a contradiction. In conclusion, with $\ee=1/(2L)$ we found $k\geq 1$ such that $Q_1,\dots,Q_L$ are linearly independent. 
    Since each $Q_\ell$ is the pushforward of $(P_{\ell})^k$ under $T$, the probability measures $(P_{1})^k,\dots,(P_{L})^k$ must be linearly independent as well.
\end{proof}

\begin{proof}[Proof of Theorem~\ref{thm:spansquared}]
The first claim follows immediately from Lemma \ref{lem:remark:power} applied to the distinct measures $P_1,\dots,P_L,Q$. Next, we focus on the case where $Q$ satisfies \eqref{AC} and $\P=\{P_1,\dots,P_L\}$ consists of independent probability measures.
Suppose that $Q\in\spn{\P}$ and $Q^2\in\spn\P^2$.
%For each $1\leq j\leq L$, we decompose $P_j$ into the absolutely continuous part $P_{j}^{\mathrm{ac}}$ and  singular part $P_{j}^{\mathrm{s}}$ with respect to $Q$. 
Denote by $f_j=\d P_{j}/\d Q$, so that by our assumption,  $f_1,\dots,f_L$ are linearly independent as functions in $L^1(Q)$. By construction, there exists a unique tuple of  nonzero numbers $(a_1,\dots,a_L)$ such that $\sum_{j=1}^L a_jP_j=Q$. In particular, $\sum_{j=1}^L a_j=1$ and \begin{align}
    1=a_1f_1+\dots+a_L f_L ~~\text{ $Q$-a.e.}\label{eq:ae1}
\end{align} 
%Integrating \eqref{eq:ae1} against $Q$ yields that for all $1\leq j\leq \ell$, $\int f_j\d Q=1$, and thus $P_j\ll Q$ for all $j$.
Since $Q^2\in\spn\P^2$, there exist $b_1,\dots,b_L$ such that for any set $A\in\F$,
$$\int_{A\times A}1Q(\d x)Q(\d y)=\int_{A\times A}\sum_{j=1}^Lb_jf_j(x)f_j(y)Q(\d x)Q(\d y).$$
By symmetry of the integrand with respect to $x,y$, we must have for any $A,B\in\F$,
$$\int_{A\times B}1Q(\d x)Q(\d y)=\int_{A\times B}\sum_{j=1}^Lb_jf_j(x)f_j(y)Q(\d x)Q(\d y).$$
By Carath\'{e}odory's extension theorem, it holds 
\begin{align}
\sum_{j=1}^Lb_jf_j(x)f_j(y)=1 ~~\text{ $Q^2$-a.e.}\label{eq:ae2}
\end{align}
By considering the a.e.~set  $y\in\X$ where \eqref{eq:ae2} holds and comparing 
    \eqref{eq:ae1} and \eqref{eq:ae2}, we have for any $1\leq j\leq L$, $f_j$ is $Q$-a.e.~constant. This implies $P_1=Q$.  Therefore, $Q\in\P$.
\end{proof}

\begin{proof}[Proof of  Proposition~\ref{prop:multidata}]
    We first show that for any $n,m\in\N$, $ \ell_{n+m}\geq  \ell_n+ \ell_m$. Suppose that $X^{(n)}$ attains the maximum e-power among pivotal and exact e-variables against $\Q^n$ for $\P^n$, and $X^{(m)}$ against $\Q^m$ for $\P^m$. Define $X^{(n+m)}(\omega_1,\omega_2)=X^{(n)}(\omega_1)X^{(m)}(\omega_2)$, where $\omega_1\in \X^n$ and $\omega_2\in \X^m$.  Clearly, $X^{(n+m)}$ is pivotal and exact against $\Q^{n+m}$ for $\P^{n+m}$. Its  e-power is $\E^{Q^{n+m}}[\log X^{(n+m)}]=\E^{Q^n}[\log X^{(n)}]+\E^{Q^m}[\log X^{(m)}]$,  {thus by construction of $X^{(n)}$ and $X^{(m)}$}, $ \ell_{n+m}\geq  \ell_n+ \ell_m$ holds. It then follows from Fekete's lemma that $ \ell_n/n$ converges to some limit in $\R\cup\{\infty\}$. Since this e-power is bounded from above by the e-power for $ \{P_1^n\}$ against $\{Q^n\}$, 
    we have for all $1\leq i\leq L$ that  $\ell_n/n \le \E^{Q^n}[\log (\d Q^n /\d {P^n_i}) ]/n=\E^Q[\log (\d Q /\d {P_i})]$ as we see in Example~\ref{ex:opt}. 
    To see that the limit is positive, it suffices to prove $\ell_k>0$ for some $k>0$ using the property $ \ell_{n+m}\geq  \ell_n+ \ell_m$. Nevertheless, that $\ell_k>0$ follows  directly from Theorem~\ref{thm:spansquared} and Theorem~\ref{thm:iff}.
\end{proof}

\section{Proof of results from Section \ref{sec:alg}}\label{sec:rate}

We start the proof of Theorem~\ref{thm:maximal} with a few simple observations.

\begin{lemma}\label{lemma:denseinsupp}
Suppose that $\rho$ is a finite measure on $\R$, and $I$ is a nonempty bounded open interval. Assume that there exists a sequence of decreasing intervals $I_n\downarrow I$, such that $\bary(\rho|_{I_n})\not\in I$ for every $n$ where the barycenter is well-defined. Then $I\subseteq (\supp\rho)^c$. 
    
\end{lemma}

\begin{proof}
This is a direct consequence of continuity of measure. We omit the details.
\end{proof}

\begin{lemma}\label{lemma:nuatomless}
Assuming \eqref{A},    any maximal element $\nu$ in $\cM_\gamma$  is atomless.
\end{lemma}

\begin{proof}Suppose that $\nu$ has an atom at $x_0\e$. Then a martingale coupling of $\nu$ and $ \gamma$ transports the mass at $x_0\e$ to some measure $\gamma'$ on $\R^L$. In particular, $\bary(\gamma')=x_0\e$ and $\gamma'\leq\gamma$. 
    By assumption \eqref{A}, $\supp\gamma'$ is not contained in any hyperplane. Using a similar argument in the proof of Proposition~\ref{prop:existence}, we conclude that there exists  a measure $\nu'\lcx \gamma'$ supported on $\I^+$ satisfying $\nu'\neq \gamma'(\R^L)\delta_{x_0}$. The new measure $\nu-\gamma'(\R^L)\delta_{x_0}+\nu'$ is then larger than $\nu$ in convex order,  contradicting the maximality of $\nu$.
\end{proof}

\begin{proof}[Proof of Theorem~\ref{thm:maximal}] Suppose that $\mu$ is the measure  from the construction, $\mu\lcx\nu$, and $\nu\lcx\gamma$ with $\nu$ supported on $\I^+$. Our goal is to show $\mu=\nu$. Let $\xi_\mu,\xi_\nu$ denote the first coordinate of $\mu,\nu$. Consider the collection $\mathcal X$ of the first coordinates of all  points $x^{(s)}_k,~s\geq 0,\,1\leq k\leq 2^{s+1}$ defined in the middle of the construction. We first show that for each $x\in\mathcal X$,
\begin{align}\label{eq:munu}
    \p[\xi_\mu\geq x]=\p[\xi_\nu\geq x]\quad\text{ and }\quad \E[\xi_\mu \mid \xi_\mu\geq x]=\E[\xi_\nu \mid \xi_\nu\geq x].
\end{align}
Note that given the first equality, the second equality in \eqref{eq:munu} is equivalent to $\E[(\xi_\mu-x)_+]=\E[(\xi_\nu-x)_+]$. The proof is similar to the ``$\Rightarrow$'' direction of Theorem~\ref{thm:maximumexist}. Let $\pi_\mu$ be any martingale coupling of $(\mu,\gamma)$ and $\pi_\nu$ be any martingale coupling of $(\nu,\gamma)$. For $x=x^{(0)}_1=\e$, by (a symmetric version of) Lemma~\ref{lemma:3}, $\xi_\mu$ attains the maximum value of $\E[(\xi_\mu-x)_+]$, and thus $\E[(\xi_\mu-x)_+]=\E[(\xi_\nu-x)_+]$.   Lemma~\ref{lemma:3} further implies that\footnote{By our assumption on $\gamma$, $\mu$ cannot have an atom at $x$.} $$\pi_\mu([0,x)\times \H_x)=\pi_\nu([0,x)\times \H_x)=1-\pi_\mu([x,\infty)\times \H_x^c)=1-\pi_\nu([x,\infty)\times \H_x^c).$$
In particular, $\p[\xi_\mu\geq x]=\p[\xi_\nu\geq x]$, proving \eqref{eq:munu}.
In the general case, consider $x=x^{(s)}_k$. There exists an interval $J$ whose endpoints are the two neighbor points of $x^{(s)}_k$ in $\{x^{(m)}_k\}_{m<s,\,1\leq k\leq 2^{m+1}}\cup\{0,\infty\}$. By definition, $\nu$ maximizes $\E[(\xi_\mu-x)_+]$, and $\mu$ maximizes $\E[(\xi_\mu-x)_+\bone_{\{\xi_\mu\in J\}}]$. Our induction hypothesis \eqref{eq:munu} applied to the right endpoint of $J$ meanwhile implies  that the two optimization problems are the same. Thus, there exists a similar block decomposition of  the supports of $\pi_\mu,\pi_\nu$ where   the total masses coincide on the blocks, and \eqref{eq:munu} holds for $x=x^{(s)}_k$. We leave the details to the reader. 

% The general case is quite similar, with an inductive block decomposition of the supports of $\pi_\mu,\pi_\nu$ where at each step, the total masses coincide on the blocks. The key is that for $x=x^{(s)}_k$, optimizing $\E[(\xi_\mu-x)_+]$ is equivalent to optimizing $\E[(\xi_\mu-x)_+|\xi_\mu\in J]$, where $J$ is the interval whose endpoints are the two neighbor points of $x^{(s)}_k$ in $\{x^{(m)}_k\}_{m<n,\,1\leq k\leq 2^{m+1}}$. 

We now finish the proof given \eqref{eq:munu}. We claim that the set $\mathcal X$ is dense in $\supp\nu$. Indeed, suppose that $I$ is an open connected component of the open set $\R\setminus \overline{\mathcal X}$. By construction and \eqref{eq:munu}, each $x^{(s)}_k$ is the barycenter of $\nu$ restricted to the interval formed by two neighbor points of $x^{(s)}_k$ in $\{x^{(m)}_k\}_{m<s,\,1\leq k\leq 2^{m+1}}$. In particular, there exist intervals $I_s\downarrow I$ where the endpoints of each $I_s$ belong to $\mathcal X$ and $\bary(\nu|_{I_s})\not\in I$. By Lemma~\ref{lemma:denseinsupp}, $I\subseteq (\supp\nu)^c$, establishing the claim.

Therefore, the distribution functions of $\mu$ and $\nu$ coincide on a dense subset $\mathcal X$ of the support of the atomless measure $\nu$. This implies $\mu=\nu$.     
\end{proof}

We next prepare for the proof of Theorem \ref{thm:e-power convergence2}. We say that a martingale $\{X_k\}_{k\geq 0}$ satisfies the \emph{separated tree condition} if 
\begin{enumerate}[(i)]
    \item 
there exists an array of real numbers $\{x^{(s)}_k\}_{s\geq 0, \,1\leq k\leq 2^s}$ such that for each $s\geq 1$, 
    % \begin{align*}
    %     x^{(s)}_1\leq x^{(s-1)}_1\leq x^{(s)}_2\leq x^{(s)}_3\leq x^{(s-1)}_2\leq\dots\leq x^{(s-1)}_{2^{s-1}}\leq x^{(s)}_{2^s}%\label{eq:tree}
    % \end{align*}
    \begin{align}
        x^{(s)}_1\leq x^{(s)}_2\leq\dots\leq x^{(s)}_{2^s}\quad\text{and}\quad x^{(s-1)}_\ell\in[x^{(s)}_{2\ell-1},x^{(s)}_{2\ell}],~1\leq \ell\leq 2^{s-1}.\label{eq:tree}
    \end{align}
    and if $A_s$ denotes the multi-set consisting of values of $\{x^{(r)}_k\}_{0\leq r< s, \,1\leq k\leq 2^r}$, 
    $$\begin{cases}
        x^{(s)}_k\geq \max\{x\in A_s\mid x\leq x^{(s-1)}_{(k+1)/2},\,x\neq x^{(s-1)}_{(k+1)/2}\}\quad&\text{if}\quad k\text{ is odd;}\\
        x^{(s)}_k\leq \min\{x\in A_s\mid x\geq x^{(s-1)}_{k/2},\,x\neq x^{(s-1)}_{k/2}\}\quad&\text{if}\quad k\text{ is even};
    \end{cases}$$
     \item $X_0$ is a constant, and for $s\geq 1$,
    $$\supp{(X_{s}\mid X_{s-1}=x^{(s-1)}_{k})}=\{x^{(s)}_{2k-1},x^{(s)}_{2k}\}.$$
\end{enumerate}
In particular, if a martingale satisfies the separated tree condition and the inequalities in \eqref{eq:tree} are strict (i.e.,~$x^{(s)}_1< x^{(s)}_2<\dots< x^{(s)}_{2^s}$ for all $s\geq 1$), the martingale is \emph{backward deterministic}, meaning that $\{X_j\}_{j=1}^n$ is $\sigma(X_n)$-measurable for all $n\in\N$ (Section 3.2 of \cite{nutz2022martingale}).

Intuitively, the separated tree condition asserts that the martingale can be represented using a binary tree on $\R$ (with $k$ corresponding to the depth of the tree), and the branches from all different levels, when projected to the real line as intervals, are either disjoint or have containment relationship.

\begin{example}
    The Simons martingale introduced by \cite{simons1970martingale} satisfies the separated tree condition; see Example \ref{ex:simons}. More generally, the SHINE martingale defined by \eqref{eq:Xn} also satisfies the separated tree condition.
\end{example}

\begin{proof}[Proof of Lemma \ref{thm:general convergence}]
The martingale property  implies \begin{equation}\label{eq:mart-seq} \E[(X_\infty-X_{k})^2] 
=\E[X_\infty^2-X_k^2]
=\sum_{j=k}^\infty \E[ X_{j+1}^2- X_j^2] 
=\sum_{j=k}^\infty \E[(X_{j+1}-X_j)^2], \end{equation} and it suffices to bound $\E[(X_k-X_{k+1})^2]$ for each $k$.  
By conditioning on $X_k$, we have
    \begin{align}
        \E[(X_k-X_{k+1})^2]\leq \sum_{j=1}^{2^{k+1}}p_jd_j^2,\label{eq:diff1}
    \end{align}
    where the index $j$ refers to the $2^{j+1}$ \emph{legs} from the support of $X_k$ to the support of $X_{k+1}$, and $p_j,d_j$ are the probability and the displacement (in absolute value) carried by the leg $j$. 

  The next key step is to discover an upper bound for each $p_j$, possibly in terms of $d_j$, using the following Lemma \ref{lemma:2seqs}. The probability $p_j$ can be written naturally as a product of conditional probabilities along a spine (or a directed path from the root of the binary tree) in the binary tree representation of $\{X_k\}_{k\geq 0}$. We label the (absolute value of the) displacement along the spine from level $k-1$ to $k$ by $b_k$, and (absolute value of) the other displacement by $a_k$. In other words, if $\{x_k\}_{k\geq 0}$ denotes the nodes (values) on the spine and each $x_k$ has descendants $x_{k+1}$ and $x_{k+1}'$, then $b_k=|x_{k-1}-x_{k}|$  and $a_k=|x_{k-1}-x_{k}'|$. By the martingale property, the total probability on the spine $\{x_k\}_{k\geq 0}$ until level $n$ is given by $\prod_{k=1}^na_k/(a_k+b_k)$.

Suppose that the martingale is uniformly bounded by $M$. By the separated tree condition, $\sum_j d_j\leq 2M$. Next, it is elementary to check that by the separated tree condition, $b_{k-2}\geq\min\{a_k,a_{k-1}\}$ and $b_{k-1}\geq\min\{a_k,b_k\}$. 
 By Lemma \ref{lemma:2seqs} applied with $n=k+1$, we conclude that $p_j\leq \sqrt{2M/d_j}\,r^{k-1}$, where $r<0.827$. Therefore,
    \begin{align}
        \sum_{j=1}^{2^{k+1}}p_jd_j^2\leq r^{k-1}\sqrt{2M}\sum_{j=1}^{2^{k+1}}d_j^{3/2}\leq 4M^2 r^{k-1}.\label{eq:diff2}
    \end{align}
%    \red{can improve this for unbounded? $\sum_{j=1}^{2^{k+1}}d_j^{3/2}$??}
    Combining \eqref{eq:mart-seq}, \eqref{eq:diff1} and \eqref{eq:diff2},  we get
    \begin{align*}
        \E[(X_k-X_\infty)^2]&\leq CM^2r^k,
        % \sum_{k=n}^\infty |\E[f(X_k)]-\E[f(X_{k+1})]|\\
        % &\leq \sum_{k=n}^\infty 2M^2\sup |f^{(2)}|r^{k-1}=\frac{2M^2\sup |f^{(2)}|r^{n-1}}{1-r},
    \end{align*}
    as desired.

 Let us now consider the general case where $X$ is not uniformly bounded. Take $M>|X_0|$.   Let $\tau_M$ be the first hitting time to $\{x:|x|\geq M\}$ of  $\{X_k\}_{k\ge 0}$.  
We use the same notation as   in \eqref{eq:diff1} for a fixed $k$.
Let $J_M\subseteq \{1,\dots,2^{k+1}\}$ be the collection of all indices of possible paths (up to time $k$) with $\tau_M=\infty$.
Note that 
$$
 \E[\bone_{\{\tau_M=\infty\}} ( X_{k+1}^2- X _k ^2)]  = \sum_{j \in J_M} p_j d_j^2.
$$
Moreover, $\sum_{j\in J_M} d_j\le 2M$.
The same argument as in \eqref{eq:diff2} gives
$$
 \sum_{j \in J_M} p_j d_j^2 \le 4 M^2 r^{k-1}. 
$$
  This implies $\E[\bone_{\{\tau_M=\infty\}} ( X_\infty ^2- X _k ^2)]  \le CM^2 r^k$ for some $C>0$. Hence, 
\begin{align*}
\E[(X_\infty-X_k)^2] & =\E[X_\infty^2-X_k^2] %\\& \le  \E[X_\infty^2-  (X^M_k)^2]  
\\ & =  \E[\bone_{\{\tau_M<\infty\}} ( X_\infty  ^2-  X _k ^2)] + \E[\bone_{\{\tau_M=\infty\}} ( X_\infty ^2- X _k ^2)]  
%\\ &=    \E[\bone_{\{\tau_M<\infty\}} ( X_\infty  ^2-  X _k ^2)] + \E[\bone_{\{\tau_M=\infty\}} ((X_\infty^M)^2-(X^M_k)^2)]    
\\& \le \E[\bone_{\{\tau_M<\infty\}} X_\infty^2]   + C M^2 r^{k}. 
\end{align*}

Let $\tau_{M}^+$ be the first hitting time to $\{x:x\ge M\}$ 
and $\tau_{M}^-$ be the first hitting time to $\{x:x\le -M\}$. 
Note that by the separated tree condition, if $M>X_0$, then 
$\tau_{M}^+<\infty$ implies $X_k \ge X_0$ for all $k\ge 0$ and $X_1>X_0$. Therefore, conditional on the event $\tau_{M}^+<\infty$, $X$ is a martingale that is bounded from below by $X_0$. 
  By Ville's inequality, for $M>X_0$,
\begin{align*}\p(\tau_{M}^+<\infty)
&\le \p(\tau_{M}^+<\infty\mid X_1>X_0)
\\&= \p\left(\sup_{k\ge 0} X_k\geq M\mid X_1>X_0\right)\le \frac{\E[X_\infty\mid X_1>X_0]-X_0} {M-X_0}=
O(M^{-1}).
\end{align*}
Similarly,  the same analysis holds for $\tau_M^-$. Hence,  we can conclude   $\p(\tau_{M}<\infty)=O(M^{-1}).$

Since $X_\infty\in L^{2+\ee}$ for $\ee>0$, $ \E[\bone_{\{\tau_M<\infty\}} X_\infty^2]= O(M^{-\delta'})$ for some $\delta' >0$ by H\"older's inequality.
Taking $M=r^{-k/(2+\delta')}$ and $q=r^{\delta'/(2+\delta')}\in(0,1)$ yields that for some $C>0$ (that may vary from line to line),
$$
\E[(X_\infty-X_k)^2]\leq C( M^{-\delta'} +  M^2 r^{k})  \leq C q^{-k}.
$$
This completes the proof in the general case.
\end{proof}

\begin{lemma}\label{lemma:2seqs}
    Suppose that non-negative numbers $a_1,\dots,a_n,b_1,\dots,b_n$ satisfy:
    \begin{itemize}
        \item $b_{k-2}\geq\min\{a_k,a_{k-1}\}$;
        \item $b_{k-1}\geq\min\{a_k,b_k\}$.
    \end{itemize}
    Then there exists $r<0.827$  such that\footnote{We set $0/0=0$.} 
    \begin{align}
        \prod_{k=1}^n \frac{a_k}{a_k+b_k}\leq \sqrt{\frac{b_1}{b_n}}\,r^{n-2}.\label{eq:produc}
    \end{align}
\end{lemma}

\begin{proof}
Consider a number $K>0$ to be determined. Define
\begin{align}
    c_k=\Big(\frac{a_k}{a_k+b_k}\Big)^2\frac{b_k}{b_{k-1}}.\label{eq:ck}
\end{align}
For $k\geq 2$, we consider the following four different cases.
\begin{enumerate}[(i)]
    \item If $b_k> b_{k-1}$, then $a_k\leq b_{k-1}< b_k$, so that
    $$\frac{a_k}{a_k+b_k}\leq \frac{b_{k-1}}{b_k+b_{k-1}}\leq \frac{1}{2}\sqrt{\frac{b_{k-1}}{b_k}}.$$
    In this case, $c_k\leq 1/4$.
    \item If $b_k\leq b_{k-1}$ and $b_{k-1}> b_{k-2}$, then similarly as in case (i), 
    $$\frac{a_{k-1}}{a_{k-1}+b_{k-1}}\leq \frac{b_{k-2}}{b_{k-1}+b_{k-2}}\leq \frac{1}{2}\sqrt{\frac{b_{k-2}}{b_{k-1}}},$$
    so that $c_{k-1}\leq 1/4$. On the other hand, $c_k\leq 1$. Therefore, $c_kc_{k-1}\leq 1/4$.
    \item If $b_k\leq b_{k-1}\leq b_{k-2}$, and $a_{k-1}\leq Kb_{k-1}$ or $a_k\leq Kb_k$, then either $c_{k-1}\leq K^2/(1+K)^2$ and $c_k\leq 1$, or $c_{k-1}\leq 1$ and $c_k\leq K^2/(1+K)^2$. In this case, $c_kc_{k-1}\leq K^2/(1+K)^2$.
    \item If $b_k\leq b_{k-1}\leq b_{k-2}$, $a_{k-1}>Kb_{k-1}$, and $a_k> Kb_k$, then
    $$b_{k-2}\geq \min\{a_k,a_{k-1}\}\geq K\min\{b_k,b_{k-1}\}\geq Kb_k.$$
    In this case, $c_kc_{k-1}\leq 1/K$.
\end{enumerate}
Upon decomposing the sequence $\{b_k\}$ into increasing and decreasing parts, the above four cases together yield that 
\begin{align}
    \prod_{k=2}^n c_k\leq \max\Big\{\frac{1}{2},\frac{K}{1+K},\frac{1}{\sqrt{K}}\Big\}^{n-2}.\label{eq:ckbound}
\end{align}
%The maximum on the right-hand side is minimized by $L\approx 2.148$, and t
To optimize our bound, we consider the quantity
$$\inf_{K>0}\max\Big\{\frac{1}{2},\frac{K}{1+K},\frac{1}{\sqrt{K}}\Big\}.$$
The minimax value is denoted by $r^2<0.683$. By \eqref{eq:ck} and \eqref{eq:ckbound}, we obtain
$$\frac{b_n}{b_1}\prod_{k=2}^n \Big(\frac{a_k}{a_k+b_k}\Big)^2\leq r^{2(n-2)}.$$
Rearranging gives \eqref{eq:produc}.
\end{proof}

Next, we apply Lemma \ref{thm:general convergence} to deduce the exponential convergence of the e-power. Suppose $f\in C^2(\R)$ has $\sup|f''|<\infty$. Then it would follow from Taylor's theorem and Lemma \ref{thm:general convergence} that
$$|\E[f(X_k)]-\E[f(X_\infty)]|\leq C\sup|f''|r^k.$$
However, in our case of main interest (see e.g., \eqref{eq:opt2}), $f(x)=-\log x$. To overcome the difficulty arising from $\sup|f^{(2)}|=\infty$, we will assume an anti-concentration bound of $X_\infty$ at $0$.

\begin{corollary}\label{thm:2nd}
Suppose that $\{X_k\}_{k\geq 0}$ is a $\R_+$-valued martingale satisfying the separated tree condition with $X_k\to X_\infty$ a.s.~and $\E[ X_\infty ^{2+\epsilon}]<\infty$ for some $\epsilon>0$ and $\E[  X_\infty ^{-2}]<\infty$. Then there exist $r <1$ and a constant $C>0$ such that
$$|\E[\log X_k]-\E[\log X_\infty]|\leq Cr^k.$$ 
\end{corollary}

\begin{proof}
For each $k\in \mathbb N$, by the mean-value theorem and H\"older's inequality,
\begin{align*}
|\E[\log(X_k)]-\E[\log(X_{\infty})]|
&\leq \E\left[|X_{\infty}-X_k| \max\left\{X_k^{-1},  X_{\infty}^{-1} \right\}\right]
\\&\le\E\left[(X_{\infty}-X_k)^2\right]^{1/2}\E \left[\max\left\{  X_k^{-2},  X_{\infty}^{-2}\right\}\right]^{1/2}.
\end{align*}
Note that $ \E [(X_{\infty}-X_k)^2 ] \le c q^{-k}$ for some $q\in (0,1)$ and $c>0$, as implied by Lemma \ref{thm:general convergence}. 
Moreover, $(X_k)_{k\ge 0}$ is a nonnegative martingale, which implies  
$
\E [   X_k^{-2} ]
\le \E[X_{\infty}^{-2}]$.
Hence, 
 $$\E \left[\max\left\{  X_k^{-2},  X_{\infty}^{-2}\right\}\right] 
\le \E \left[   X_k^{-2} + X_{\infty}^{-2} \right]
\le 2\E[X_\infty^{-2}]<\infty.$$ 
Therefore, we get 
 $
|\E[\log(X_k)]-\E[\log(X_{\infty})]| \le Cr^{- k}
 $ for some $C>0$ and $r=q^{1/2}\in (0,1)$. 
\end{proof}

 \begin{proof}[Proof of Theorem \ref{thm:e-power convergence2}] By property of the SHINE construction (Theorem \ref{thm:maximal}), $$X_\infty \e\lcx \left(\frac{\d P_1}{\d Q},\dots,\frac{\d P_L}{\d Q}\right)\Big|_{Q}.$$Therefore,  for each $j\in\{1,\dots,L\}$, we have by linearity of the coordinate map that $(\d P_j/\d Q)|_Q\gcx X_\infty$. In particular, using the definition of the convex order we have
 $$\E^Q\Big[\Big(\frac{\d P_j}{\d Q}\Big)^{2+\ee}\Big]<\infty\quad\implies\quad\E[ X_\infty ^{2+\epsilon}]<\infty$$
 and
    $$\E^Q\Big[\Big(\frac{\d P_{j'}}{\d Q}\Big)^{-2}\Big]<\infty\quad\implies\quad \E[X_\infty^{-2}]<\infty. $$
 These verify the assumptions in Corollary \ref{thm:2nd}. The conclusion of Theorem \ref{thm:e-power convergence2} then follows directly from that of Corollary \ref{thm:2nd}.     
 \end{proof}

\begin{example}
    Suppose that $P_j\sim\mathrm{N}({\mathbf m}_j,1),\,1\leq j\leq L$, and $Q\sim\mathrm{N}({\mathbf m}_Q,1)$ are nondegenerate Gaussian distributions. The laws of the Radon--Nikodym derivatives are
    $$\frac{\d P_j}{\d Q}\Big|_{Q}=e^{\boldsymbol{\xi}\cdot({\mathbf m}_j-{\mathbf m}_Q)+\frac{1}{2}(\n{{\mathbf m}_Q}^2-\n{{\mathbf m}_j}^2)}\Big|_{{\boldsymbol{\xi}}\lawis \mathrm{N}({\mathbf m}_Q,1)}.$$
    Since the moment generating function of a nondegenerate multivariate Gaussian distribution is well-defined everywhere, the conditions \eqref{eq:moment1} and \eqref{eq:moment2} are satisfied. Therefore, Corollary \ref{thm:e-power convergence2} implies that the SHINE construction enjoys exponential convergence when testing $P_1,\dots,P_L$ against $Q$.
\end{example}

\section{Proof of results from Section \ref{sec:composite alt}}\label{app:sec6proof}
 To prove Theorems \ref{thm:iff2} and \ref{thm:pvariable2}, we build upon the ideas from Proposition~\ref{prop:existence}.

\begin{lemma}\label{lemma:subspaces}
    Let $V\subseteq\R^d$ be a subspace containing $\e\in\R^d$ and $\S$ be a collection of affine hyperplanes
    in $\R^d$ containing $\e$ such that whenever $S\in\S$ and an affine subspace $T$ satisfies $T\cap V\subseteq S$, it holds $T\subseteq S'$ for some $S'\in\S$.   Then for each measure $\mu$ centered at $\e$ whose support is not a subset of $S$ for any $S\in\S$, there exists a measure $\nu$ on $V$ such that $\nu\lcx\mu$ and the support of $\nu$ is not a subset of $S$ for any $S\in\S$.
\end{lemma}

\begin{proof}
\sloppy The proof is similar to the proof of Proposition~\ref{prop:existence}. Define $T=\aff\supp\mu$. By Lemma~\ref{lemma:CM}, it suffices to find points $s_1,\dots,s_k\in\supp\mu$  such that $\ri(\conv\{s_1,\dots,s_k\};T)$ contains $\e$ and intersects with $V$ not on a single $S\in\S$.

\sloppy Suppose that the contrary holds. That is,  any $s_1,\dots,s_k\in\supp{\mu}$ satisfies  $\e\not\in\ri(\conv\{s_1,\dots,s_k\};T)$ or $\ri(\conv\{s_1,\dots,s_k\};T)\cap V\subseteq S$ for some $S\in\S$. 
By Lemma~\ref{lemma:elem}(i), any $\ri(\conv\{s_1,\dots,s_k\};T)$ is contained in  $\ri(\conv\{s_1,\dots,s_K\};T)$ for some $k\leq K$ and $s_1,\dots,s_K\in\supp\mu$ such that $\e\in\ri(\conv\{s_1,\dots,s_K\};T)$. This implies for all $s_1,\dots,s_k\in\supp{\mu}$ that $\ri(\conv\{s_1,\dots,s_k\};T)\cap V\subseteq S$ for some $S\in\S$. Consequently, there exists $S\in\S$ such that $\ri(\conv\supp\mu;T)\cap V\subseteq S$. By Lemma~\ref{lemma:elem}(ii), $T\subseteq \aff\ri(\conv\supp\mu;T)$. Since $V,S$ are affine spaces, it holds that 
 $T\cap V\subseteq S$. Moreover, $\supp\mu\subseteq T\subseteq S'$ for some $S'\in\S$ by our assumption. 
Hence, the support of $\mu$ is contained in $S'$, contradicting our assumption.
\end{proof}

\begin{proposition}\label{prop:existence2}
\sloppy Let $L,M\in\N$ and $(P_1,\dots,P_L,Q_1,\dots,Q_M)$ be a jointly atomless tuple of probability measures on $\X$ such that $\spn(P_1,\dots,P_L)\cap\conv(Q_1,\dots,Q_M)=\emptyset$. Then there exist probability measures $F,G_1,\dots,G_M$ on $\R$ such that $F\not\in\conv(G_1,\dots,G_M)$ and
$$\T((P_1,\dots,P_L,Q_1,\dots,Q_M),(F,\dots,F,G_1,\dots,G_M))\neq\emptyset.$$
\end{proposition}

\begin{proof}
\sloppy Let $\mu$ be a dominating measure for $(P_1,\dots,P_L,Q_1,\dots,Q_M)$, say $\mu=(P_1+\dots+P_L+Q_1+\dots+Q_M)/(L+M)$, and $\nu=\Leb$. Then  $\spn(P_1,\dots,P_L)\cap\conv(Q_1,\dots,Q_M)\neq \emptyset$ is equivalent to the existence of $\{\alpha_i\}_{1\leq i\leq L}$  and $\{\beta_j\}_{1\leq j\leq M}$ such that 
$$\sum_{i=1}^L\alpha_i\frac{\d P_i}{\d\mu}= \sum_{j=1}^M\beta_j\frac{\d Q_j}{\d\mu},~\beta_j\geq 0,~\sum_{i=1}^L\alpha_i=\sum_{j=1}^M\beta_j=1.$$Similarly, $F\in\conv(G_1,\dots,G_M)$ is equivalent to the existence of $\{\lambda_j\}_{1\leq j\leq M}$ such that 
$$\frac{\d F}{\d\nu}=\sum_{j=1}^M\lambda_j \frac{\d G_j}{\d \nu},~\lambda_j\geq 0,~\sum_{j=1}^M\lambda_j=1.$$

To this end, we define
$$\S=\left\{S_{\ba,\bb}\mid   \beta_j\geq 0,~\sum_{i=1}^L\alpha_i=\sum_{j=1}^M\beta_j=1\right\},$$where
$$S_{\ba,\bb}:=\left\{(x_1,\dots,x_p,y_1,\dots,y_q) \mid  \sum_{i=1}^L\alpha_ix_i=\sum_{j=1}^M\beta_jy_j\right\}.$$
We now claim that for each measure $\gamma$ such that $\supp\gamma$ is not a subset of some $S\in\S$, there exists $\tau\lcx\gamma$ such that $\tau$ is supported on $V:=\{(\bx,\by)\in\R^{L+M} \mid x_1=\dots=x_L\}$ but not concentrated on a single $V_{\bl}:=\{(\bx,\by)\in\R^{L+M} \mid x_1=\dots=x_L=\sum_{j=1}^M\lambda_jy_j\}$ for all $\lambda_j\geq 0, \sum_{j=1}^M\lambda_j=1$. Provided the claim is true, we construct using Lemma~\ref{lemma:RN}  applied with $d=L+M$ the measures $F,G_1,\dots,G_M$ such that
$$\left(\frac{\d F}{\d\nu},\dots,\frac{\d F}{\d\nu},\frac{\d G_1}{\d\nu},\dots,\frac{\d G_M}{\d\nu}\right)\Big|_\nu =\tau.$$
Since $\nu\lcx\gamma$ and $\gamma$ is supported on the hyperplane $\{(\bx,\by) \mid \sum_{i=1}^Lx_i+\sum_{j=1}^My_j=L+M\}$, our measure $\nu$ will be supported on the same hyperplane, thus $\nu=(LF+\sum_{j=1}^M G_j)/(L+M)$. In other words, $\mu,\nu$ allow the same linear combination of the measure tuples $(P_1,\dots,P_L,Q_1,\dots,Q_M)$ and $(F,\dots,F,G_1,\dots,G_M)$. By Proposition~\ref{prop:ssww} applied with $d=L+M$,  $$\T((P_1,\dots,P_L,Q_1,\dots,Q_M),(F,\dots,F,G_1,\dots,G_M))\neq\emptyset$$ holds as desired.

To prove the above claim, we apply Lemma~\ref{lemma:subspaces} with $d=L+M$, $\mu=\gamma$, $\nu=\tau$, and $V,\S$ defined as above.
Note that if the support of $\tau$ is contained in $V$ but not in a certain $S$, then it cannot be contained in a certain $V_{\bl}$. Thus the conclusion of Lemma~\ref{lemma:subspaces} suffices for our purpose.

It then suffices to check the condition in Lemma~\ref{lemma:subspaces}  that whenever $S\in\S$ and an affine subspace $T$ satisfies $T\cap V\subseteq S$, it holds $T\subseteq S'$ for some $S'\in\S$. To this end, we consider $S_{\ba,\bb}\in\S$ and first assume $T$ is a   hyperplane containing $S_{\ba,\bb}\cap V=\{(\bx,\by)\mid x_1=\dots=x_L=\sum_{j=1}^M\beta_jy_j\}$. In this case, a normal vector to $T$ (which is unique up to a multiplicative constant) must also be a normal vector of $S_{\ba,\bb}\cap V$, and hence must be of the form $(t_1,\dots,t_L,-\beta_1\sum_{i=1}^Lt_i,\dots,-\beta_M\sum_{i=1}^Lt_i)$ for some $t_1,\dots,t_L\in\R$.  If $\sum_{i=1}^Lt_i=0$, then $T\supseteq V$, and thus $V\subseteq S_{\ba,\bb}$, which is impossible. Thus $\sum_{i=1}^Lt_i\neq 0$. It follows that for some $t_1,\dots,t_L$ with $\sum_{i=1}^Lt_i\neq 0$,
$$T=\left\{(x_1,\dots,x_L,y_1,\dots,y_M) \mid  \sum_{i=1}^Lt_ix_i=\sum_{i=1}^Lt_i\sum_{j=1}^M\beta_jy_j\right\}.$$
Therefore, $T\in\S$. More precisely, $T=S_{\mathbf{t}/\sum_{i=1}^L t_i,\bb}$.

Next, we prove the general case of an affine subspace $T$ that satisfies $T\cap V\subseteq S_{\ba,\bb}$. Note that $V\not\subseteq S_{\ba,\bb}$ for each $S_{\ba,\bb}\in\S$, and that $V$ is of dimension $M+1$. Thus $S_{\ba,\bb}\cap V$ is of dimension $M$.  Let $T'=T+(S_{\ba,\bb}\cap V)+V^\perp$. Then  $T'$ is a hyperplane, $T\subseteq T'$, and $$T'\cap V\subseteq (T+(S_{\ba,\bb}\cap V))\cap V= (T\cap V)+(S_{\ba,\bb}\cap V)\subseteq S_{\ba,\bb}.$$ This completes the proof.
\end{proof}

\begin{proof}[Proof of Theorem~\ref{thm:iff2}]
%(i) The ``only if" is clear. The ``if" follows from Hahn-Banach: let $\mu$ be a probability measure dominating $P_i$ and $Q_j$ for all $i,j$. By assumption, $\conv(\frac{\d P_1}{\d \mu},\dots,\frac{\d P_L}{\d \mu})\cap \conv(\frac{\d Q_1}{\d\mu},\dots,\frac{\d Q_M}{\d\mu})=\emptyset$. Since the dual of $L^1(\mu)$ is $L^\infty(\mu)$, there exists a bounded random variable $Y$ with $\E^\mu[\frac{\d P_i}{\d \mu}Y]\leq 0<\E^{\mu}[\frac{\d Q_j}{\d\mu}Y]$ for all $i,j$. Then let $X=1+Y/\sup{|Y|}$. To see that $\E^Q[\log X]>0$ in both settings, we may replace $X$ by $(1-b)+bX$ with $b>0$ chosen  small enough, since our choice of $X$ is bounded.

The direction (a)$\Rightarrow$(b) is  Proposition~\ref{prop:calibration}, (b)$\Rightarrow$(c) is clear, (c)$\Rightarrow$(e) being precisely Proposition~\ref{prop:nec}, and (e)$\Rightarrow$(d) is Proposition~\ref{prop:existence2}. 
To show (d)$\Rightarrow$(a), we let $\phi$ be a nontrivial p-variable with null $\{G_1,\dots,G_M\}$ and alternative $\{F\}$, whose existence is guaranteed by Theorem~\ref{thm:p-variable}.  Then by definition, $\phi\circ X$ has a common law that is $\pst\Leb$ under each $P_i$, and has a law that is $\gst \Leb$ under each $Q_j$. Applying Lemma~\ref{lemma:quantile}(i) then yields a random variable $\Psi$ such that $\Psi\circ \phi\circ X$ is an exact and nontrivial p-variable as desired. Finally, by Proposition~\ref{prop:no ja}, the direction (c)$\Leftrightarrow$(e) also holds without condition \eqref{JA}.
\end{proof}

\begin{proof}[Proof of Theorem \ref{thm:pvariable2}]
    The proofs are similar to  Theorem~\ref{thm:p-variable},
     where in the direction (d)$\Rightarrow$(a) we  replace the linear cone $(-\infty,0)^L\cup(0,\infty)^L\cup\{\z\}$ by the linear cone $(-\infty,0)^L\times (0,\infty)^M\cup (0,\infty)^L\times (-\infty,0)^M\cup\{\z\}$. 
    The last statement is verified by Proposition~\ref{prop:no ja} and Remark \ref{remark:kraft}.
\end{proof}

\begin{proof}[Proof of Theorem~\ref{thm:infinite case}]
We first show the ``only if'' direction. Suppose that $0\in\overline{\overline{\spn}\P+\overline{\conv}\Q}$ and $X$ is an exact and   bounded e-variable satisfying $\inf_{Q\in\Q}\E^Q[\log X]>0$. In particular, since $X$ is bounded, for a sequence of distributions converging in total variation, the expectations of $X$ also converge.
Let $P^{(n)}\in \overline{\spn}\P$ and $Q^{(n)}\in \overline{\conv}\Q$ be such that $P^{(n)}+Q^{(n)}\to 0$ in total variation. It follows that $\E^{P^{(n)}}[X]=1$ and $\E^{Q^{(n)}}[X]\geq \inf_{Q\in\Q}\E^{Q}[X]>1$. But then $\liminf\E^{P^{(n)}+Q^{(n)}}[X]\geq \inf_{Q\in\Q}\E^{Q}[X]-1>0$, a contradiction.

   Next, we show the ``if'' direction.  Let $R\in\Pi(\X)$ be as given. We may abuse notation and identify each $P\in\P$ and $Q\in\Q$ with its density  with respect to  $R$. Clearly, convergence in total variation is equivalent to convergence in $L^1(R)$. Thus, $S:=\overline{\spn}\{P\mid P\in\P\}$ is a closed subspace of $L^1(R)$. By assumption, the set $C:=\overline{\conv}\{Q\mid Q\in\Q\}$ satisfies that  $\overline{C+S}$ is closed, convex, and disjoint from $0$ in the quotient space $L^1(R)/S$. By  the  Hahn-Banach separation theorem, there is $\overline{h}:L^1(R)/S\to\R$ such that $\overline{h}|_{\overline{C+S}}>\ee>0$. Composing with the quotient map we obtain a linear functional $h:L^1(R)\to\R$, and it is easy to check that $h$ vanishes on $S$ and $h|_{C}>\ee$. By duality, we may recognize $h\in L^\infty(R)$. It follows that the bounded random variable $X=h+1$ satisfies $\E^{P}[X]=1+\E^P[h]=1+\int hP\d R=1$ for each $P\in\P$ and $\E^Q[X]=1+\int hQ\d R>1+\ee$ for each $Q\in\Q$. Proposition~\ref{prop:e-var} then concludes the proof. 

   Suppose that $\Q$ is tight and that there exist $P^{(n)}\in \overline{\spn}\P$ and $Q^{(n)}\in \overline{\conv}\Q$  such that $P^{(n)}+Q^{(n)}\to 0$.  By Prokhorov's theorem, $\overline{\conv}\Q$ is weakly compact. This implies for some subsequence $\{n_k\}$, $Q^{(n_k)}$ is convergent. The limit then belongs to $\overline{\spn}\P\cap \overline{\conv}\Q$.
   The other direction is obvious.
\end{proof}

\begin{proof}[Proof of Proposition~\ref{prop:id}]
By \citet[Lemma 7.5]{sato1999levy}, the Fourier transform of  the density $p$ of an infinitely divisible distribution has no real zeros. By Wiener's Tauberian theorem (Theorem 8 of \cite{wiener1988fourier}), the linear span of the set of translates $\{p(\cdot-\theta)\}_{\theta\in\R}$ is dense in $L^1(\R)$. Therefore, there is $P\in\spn\{P_\theta:\theta\in\R\}$ with density $\widetilde{p}$ such that $d_{\mathrm{TV}}(P,Q)=(\int|\widetilde{p}(x)-q(x)|\,\d x)/2<\ee$. In other words, $Q\in\overline{\spn}\{P_\theta:\theta\in\R\}$, say we have $Q=\lim_{k\to\infty} P^{(k)}$, where $P^{(k)}\in{\spn}\{P_\theta:\theta\in\R\} $. 

Suppose that $X$ is an exact  e-variable that is nontrivial for $\{Q\}$. Then there exists a large number $K>0$ such that $\widetilde{X}:=X\bone_{\{X\leq K\}}$ satisfies $\E^Q[\widetilde{X}]>1$. Since $\widetilde{X}$ is bounded, we have $$1<\E^Q[\widetilde{X}]=\lim_{k\to\infty} \E^{P^{(k)}}[\widetilde{X}]\leq \limsup_{k\to\infty} \E^{P^{(k)}}[{X}]= 1.$$This leads to a contradiction. 
\end{proof}

\section{Some technical results}\label{appendix}

\subsection{Useful results from convex analysis}

We start with a few well-known results from convex analysis. We refer the readers to \cite{rockafellar1970convex} and \cite{simon2011convexity} for more background.
\begin{lemma}\label{lemma:elem}
Let $A\subset\R^d$ be a closed set, and $\mu\in\cM(\R^d)$ with $\supp\mu=A$. Then the following statements hold:
\begin{enumerate}[(i)]
    \item $\aff A=\aff\ri(\conv A;\aff A)$;
    \item $\bary(\mu)\in\ri({\conv} A;\aff A)$.
\end{enumerate}
\end{lemma}

\begin{proof}
(i) The $\supseteq$ direction is obvious. To prove $\subseteq$, we may replace $A$ by $\conv A$ and without loss of generality assume $A$ is also convex. Let $a\in A$ and $b\in \ri(A;\aff A)$, then elementary geometric arguments show that $(a+b)/2\in\ri (A;\aff A)$; see Theorem 6.1 of \cite{rockafellar1970convex}. Thus $A\in\aff\ri(A;\aff A)$.

(ii) We may without loss of generality assume $\aff A=\R^d$ and replace the relative interior by interior. An application of the Hahn-Banach separation theorem yields $\bary(\mu)\in \overline{\conv} A$. Suppose $\bary(\mu)\not\in \ri(\overline{\conv} A;\aff A)=(\overline{\conv} A)^\circ$, then the Hahn-Banach separation theorem implies the existence of a closed hyperplane $\H\subseteq\R^d$ such that $\bary(\mu)\in\H$ and $(\overline{\conv} A)^\circ\subseteq\R^d\setminus \H$; see Theorem 11.2 of \cite{rockafellar1970convex}.  Therefore, $A\subseteq \partial\H$, contradicting $\aff A=\R^d$.
     By Theorem 6.3 of \cite{rockafellar1970convex}, $\ri(\overline{\conv} A;\aff A)=\ri({\conv} A;\aff A)$. This completes the proof.
\end{proof}

We also prove the following variant of the Choquet-Meyer theorem.

\begin{lemma}\label{lemma:CM}
 \sloppy Suppose that $\mu$ is a finite measure on $\R^d$, $x_1,\dots,x_k\in\supp\mu$, and $x\in \ri(\conv\{x_1,\dots,x_k\};\aff\supp\mu)$. Then there exists $\delta>0$ such that any measure $\gamma$ with total mass $\gamma(\R^d)\leq \delta$, supported on $B(x;\delta)\cap(\aff\supp\mu)$,  satisfies $\gamma\lcx\widetilde{\mu}$ for some $\widetilde{\mu}\leq\mu$. 
 \end{lemma}
\begin{proof}  First, we may assume without loss of generality that  $\aff\supp\mu=\R^d$, and replace the relative interior by interior. In this case, we  must have $\aff\{x_1,\dots,x_k\}=\aff\supp\mu=\R^d$, otherwise $(\conv\{x_1,\dots,x_k\})^\circ=\emptyset$ and the statement is vacuously true. 

Since $x\in (\conv\{x_1,\dots,x_k\})^\circ$, there exists $\epsilon>0$ such that the distance of $x$ from  $\partial\conv\{x_1,\dots,x_k\}$ is larger than $\epsilon$. Let $\mu_N$ for $N\in \N$ be the conditional distribution of $\mu$ given the $\sigma$-field generated by cubes with coordinates in $\mathbb Z^d/N$. The smallest cubes have size $(1/N)^d$. 
For each $j=1,\dots,k$, pick a cube $D^N_j$ of size $(1/N)^d$ in $\R^d$ containing $x_j$ (possibly on its boundary) that has a positive $\mu$-measure, which is possible since $x_j$ is in the support of $\mu$. Let $y^N_j=\bary(\mu|_{D^N_j})$. 
It is then clear that $\mu_N(\{y^N_j\})>0$ and $\mu_N(\{y^N_j\})\delta_{y^N_j}\lcx \mu|_{D^N_j}$.

 For $N>d^{3/2}/\epsilon$,
$\Vert y^N_j-x_j\Vert<\epsilon/d$.
Therefore, $x\in \ri(\conv\{y^N_1,\dots,y^N_k\};\aff \{y^N_1,\dots,y^N_k\}).$
\sloppy Fix   $N>d^{3/2}/\epsilon$
such that the boxes $\{D^N_j\}_{1\leq j\leq k}$ are disjoint.
 Write $(y_1,\dots,y_k)=(y^N_1,\dots,y^N_k)$ and $D_j=D^N_j.
$ Note that the distance between $x$ and $\partial\conv\{y_1,\dots,y_k\} $ is positive by the triangle inequality,
and hence 
  $\aff\{y_1,\dots,y_k\}=\R^d$, so that $x\in (\conv\{y_1,\dots,y_k\})^\circ$. 

Pick $\delta>0$
small enough such that $B(x;\delta)\subseteq (\conv\{y_1,\dots,y_k\})^\circ$ and that $\delta<\min \{\mu_N(\{y_1\}), \dots,\mu_N(\{y_k\})\}$. By Choquet's theorem (Theorem 10.7(ii) of \cite{simon2011convexity}), for each $y\in B(x;\delta)$, there exists a probability measure $\gamma_y$ supported on $\{y_1,\dots,y_k\}$ such that $\bary(\gamma_y)=y$, and  $\gamma_y$ is continuous in $y$. 

Consider an arbitrary measure $\gamma$ with total mass $\gamma(\R^d)\leq \delta$ and supported on $B(x;\delta)$. Define
$$\widetilde{\gamma}=\int \gamma_y\gamma(\d y).$$Observe that $\gamma\lcx \widetilde{\gamma}$ and $\widetilde{\gamma}$ is supported on $\{y_1,\dots,y_k\}$ with $$\widetilde{\gamma}(\{y_j\})=\int \gamma_y(\{y_j\})\gamma(\d y)\leq \gamma(\R^d)\leq \delta\leq \mu_N(\{y_j\}),\ 1\leq j\leq k.$$  
Define 
$$\widetilde{\mu}=\sum_{j=1}^k\left(\frac{\widetilde{\gamma}(\{y_j\})}{\mu_N(\{y_j\})}\right)\mu|_{D_j}.$$
It follows that
$$\gamma\lcx \widetilde{\gamma}=\sum_{j=1}^k \left(\frac{\widetilde{\gamma}(\{y_j\})}{\mu_N(\{y_j\})}\right)\mu_N(\{y_j\})\delta_{y_j}\lcx \sum_{j=1}^k \left(\frac{\widetilde{\gamma}(\{y_j\})}{\mu_N(\{y_j\})}\right)\mu|_{D_j}=\widetilde{\mu}.$$
Since $\{D_j\}_{1\leq j\leq k}$ are disjoint,
 $$\widetilde{\mu}\leq \sum_{j=1}^k\mu|_{D_j}\leq \mu,$$
as desired.
\end{proof}

\subsection{Useful results about the stochastic order \texorpdfstring{$\lst$}{}}

We next state and prove a few well-known results regarding the stochastic order $\lst$. These are useful when proving the existence of  p-variables.

\begin{lemma}\label{lemma:p-variable}
    Suppose that $F,G\in\Pi(\R)$ are atomless and $F\neq G$. Then there exists a bounded random variable $\phi$ on $\R$ such that its law under $F$ is $\Leb$ and its law under $G$ is $\lst\Leb$ but distinct from $\Leb$.
\end{lemma}

\begin{proof}
We pick a random variable $\phi$ that has law $\Leb$ and is comonotone with $\d G/\d F$ under the law $F$. In particular, $\phi$ and $\d G/\d F$  are positively associated. Therefore, for $\alpha\in[0,1]$,
$$G(\phi\leq\alpha)=\int\frac{\d G}{\d F}\bone_{\{\phi\leq\alpha\}}\d F\geq \int \frac{\d G}{\d F}\d F\int \bone_{\{\phi\leq\alpha\}}\d F=F(\phi\leq\alpha)=\alpha.$$
Since $\d G/\d F$ is not a constant under $F$, there exists $\alpha\in(0,1)$ such that the inequality is strict. 
\end{proof}

\begin{lemma}\label{lemma:quantile}
Suppose that $F_1,\dots,F_L,G$ are atomless probability measures on $[0,1]$.
\begin{enumerate}[(i)]
    \item If $F_i\gst \Leb$ for all $i$ and $G\pst\Leb$, then there exists a random variable $\Psi:[0,1]\to[0,1]$ such that $\Psi|_{F_i}\sst\Leb$ for all $i$ and $\Psi|_G=\Leb$.

    \item If there exists $\beta\in(0,1)$ such that $\d F_i/\d\Leb\leq 1$ on $[0,\beta)$ and $\d F_i/\d\Leb\geq 1$ on $(\beta,1]$, and $F_i\sst\Leb$ for all $i$ and $G=\Leb$, then  there exists a random variable $\Psi:[0,1]\to[0,1]$ such that $\Psi|_{F_i}\gst\Leb$ for all $i$ and $\Psi|_G\pst\Leb$.
\end{enumerate}
\end{lemma}

\begin{proof}
(i) Let $F=\max F_i$ and $\mathrm{Id} $ be the identity on $[0,1]$.
    Let $\widetilde F$, $\widetilde G$ and $\widetilde F_i$ be the corresponding cdfs. 
    $\widetilde G > \mathrm{Id} \ge \widetilde F\ge  \widetilde F_i $ and $\widetilde F_i |_{F_i} \gst  \mathrm{Id} | _{F_i} \laweq \Leb$ for each $i$. Hence $\widetilde G$ follows $\Leb$ under $G$ and it stochastically dominates $\Leb$ under each $F_i$ by Theorem 1.A.3.(a) of \cite{shaked2007stochastic}.

(ii) Denote by $\alpha=\min\{\beta-F_i([0,\beta))\}>0$. Pick $\tau$ such that $\max F_i([\beta,\beta+2\tau))<\alpha$. Define 
\begin{align*}
    \Psi(x)=\begin{cases}
        x&\text{ if }x\in [0,\beta+\tau]\cup[\beta+2\tau,1];\\
        x-\tau&\text{ otherwise}.
    \end{cases}
\end{align*}
 By construction, it is then easy to check that $\Psi|_{\Leb}\pst\Leb$ and $\Psi|_{F_i}\gst\Leb$.
\end{proof}

\end{appendix}

\end{document}